\documentclass[10pt]{article}

\usepackage{amsmath}
\usepackage{graphicx}
\usepackage{latexsym}
\usepackage{amssymb}
\usepackage{amsbsy}
\usepackage{color}
\usepackage{multirow} 
\usepackage{rotating}
\usepackage{bm}   
\usepackage{subfigure}  
\usepackage[title]{appendix}

\usepackage{hyperref} 
\hypersetup{
    colorlinks=true,                          
    linkcolor=blue, % Couleur des liens internes
    citecolor=red, % Couleur des numéros de la biblio dans le corps
    urlcolor=blue  } % Couleur des url
\newcommand{\Int}{\displaystyle\int}
\newcommand{\Sum}{\displaystyle\sum}

\DeclareSymbolFont{msbm}{U}{msb}{m}{n}
\DeclareMathSymbol{\C}{\mathalpha}{msbm}{'103}
\DeclareMathSymbol{\R}{\mathalpha}{msbm}{'122}
\DeclareMathSymbol{\Z}{\mathalpha}{msbm}{'132}
\DeclareMathSymbol{\N}{\mathalpha}{msbm}{'116}

\newtheorem{remark}{Remark}

\def\RR{\mathbb R}
\def\SS{\mathbb S}

\newcommand{\bi}[0]{\begin{itemize}}
\newcommand{\ei}[0]{\end{itemize}}
\newcommand{\beann}[0]{\begin{eqnarray*}}
\newcommand{\eeann}[0]{\end{eqnarray*}}
\def\be{\begin{equation}}
\def\ee{\end{equation}}
\def\bea{\begin{eqnarray}}
\def\ba{\begin{array}{l}\displaystyle}
\def\eea{\end{eqnarray}}
\def\ea{\end{array}}
\def\E{{\cal E}}

\newcommand{\ext}{\text{ext}}

\newfont{\numerikEleven}{ecrm1000}
\newfont{\numerikTen}{cmss10}
\newfont{\numerikNine}{cmss9}
\newfont{\numerikEight}{cmss8}

\usepackage[title]{appendix}

\begin{document}

%\linenumbers

%------------------------------------
% TITLE
%
\title{ \Large
A new class of high order semi-Lagrangian schemes for rarefied gas dynamics.
%\thanks{This work was supported by the French 'Agence Nationale pour la
%Recherche (ANR)' in the frame of the contract JCJC ``ALE INC(ubator) 3D''. } 
}
%------------------------------------

%------------------------------------
% AUTHORS and AFFILIATIONS
%
\author{ 
Giacomo Dimarco\thanks{Department of Mathematics and Computer Science
University of Ferrara, Ferrara, Italy. \& Universit\'{e} de Toulouse; UPS, INSA, UT1, UTM;
CNRS, UMR 5219; Institut de Math\'{e}matiques de Toulouse; F-31062
Toulouse, France. ({\tt
giacomo.dimarco@unife.it})},
\and
Cory Hauck \thanks{Computer Science and Mathematics Division;
Oak Ridge National Laboratory, US. ({\tt hauckc@ornl.gov }~)},
\and
Rapha\"{e}l Loub{\`e}re\thanks{Universit\'{e} de Toulouse; UPS, INSA, UT1, UTM;
CNRS, UMR 5219; Institut de Math\'{e}matiques de Toulouse; F-31062
Toulouse, France. ({\tt raphael.loubere@math.univ-toulouse.fr})}
}
%------------------------------------

%------------------------------------
% DATE
%
\date{\today}
%------------------------------------

\maketitle

%------------------------------------
% ABSTRACT
%
\begin{abstract}
In this paper we genealize the fast semi-Lagrangian scheme developed in 
[J. Comput. Phys., Vol. 255, 2013, pp 680-698] to the case of high order reconstructions
of the distribution function. The original first order accurate semi-Lagrangian scheme is supplemented
with polynomial reconstructions of the distribution 
function and of the collisional operator leading to an effective
high order accurate numerical scheme for all regimes, from extremely rarefied gas to highly collisional siuation. 
The main idea relies on updating at each time step the extreme points of the distribution function for each
velocity of the lattice instead of updating the solution in the cell centers, these extremes points being located at different positions for any fixed velocity of the lattice.
The result is a class of scheme which permits to preserve the structure of the solution over very long times compared to existing schemes from the literature.
We propose a proof of concept of this new approach along with numerical tests and comparisons with 
classical numerical methods.
\end{abstract}
%------------------------------------ 

%------------------------------------
% KEYWORDS
%

\tableofcontents

%------------------------------------
% INTRODUCTION
%

\section{Introduction} \label{sec:introduction}
 
Kinetic models are nowadays a very important tool which is used to describe many different phenomena ranging from
classical fluid mechanics, plasmas and rarefied gases \cite{dimarco-review,cercignani,birsdall} to biology and socio-economy models \cite{PaTo}. 
The advantage of kinetic equations with respect to more classical approaches based on fluid mechanics equations is the possibility to take into account a larger
spectrum of regimes and, consequently, of phenomena. Typically fluid mechanics describes problems which are in the so-called thermodynamics equilibrium state while the kinetic approach
permits to consider both thermodynamics equilibrium as well as far from equilibrium states for the particles constituting the system under analysis. 

Unfortunately, the strong capability in terms of modelization are 
%countered 
mitigated 
by several numerical problems of different nature \cite{dimarco-review}. 
One typically would like to preserve the physical conserved quantities also at the numerical level because they characterize the steady states and, in addition, to be able to 
deal with different spatial and temporal scales which are often encountered in practical applications.
%By the way, for sure, 
A more fundamental problem related to kinetic numerical simulation 
is represented by 
% the high dimensionality
a light version of 'the curse of dimensionality'
of these models. In general
the distribution function which defines the state of the system and which gives the probability for a particle to be at a given position with a given speed at a fixed time  
depends on seven independent variables: $(x,v,t) \in \mathbb{R}^7$. On the contrary 
fluid mechanics unknowns live in a $\mathbb{R}^4$ space: three physical space components 
plus the time. 
This makes numerical simulations of realistic multi dimensional problems 
%highly prohibitive.
demanding in terms of computational time and memory consumption.

Historically, there exists two different approaches which are typically used to tackle kinetic equations: deterministic numerical methods such as
finite volume, semi-Lagrangian and spectral schemes \cite{dimarco-review} and 
probabilistic numerical methods such as Direct Simulation Monte Carlo (DSMC) \cite{bird, Cf}.
Both methodologies have strengths and weaknesses. While the first could normally reach high order of accuracy, the second are often faster, especially
for solving steady problems but, typically, they exhibit lower convergence rate and difficulties in describing non stationary and slow motion flows.

For this reason, many different research works have been dedicated in the recent past to reduce
some of the disadvantages of these classes of methods. 
In particular, concerning the drawback of the excessive computational cost for deterministic methods we recall the fast
spectral approaches developed in \cite{Filbet2,Pareschi,Pareschi2}, while for tackling the multiple scales problem we remind the work in \cite{Jin, dimarco-IMEX,dimarco-Exp} and
the references therein to have 
%an insight on the problem and the relative numerical techniques. 
a broad view on the problematic and the associated numerical techniques.

The recent review on deterministic numerical methods
for kinetic equation \cite{dimarco-review} gives a panorama of the state of the art of numerical discretization of kinetic equation. 
Concerning the issues related to the excessive noise of Monte Carlo methods we quote \cite{Cf}
for an overview on efficient and low variance Monte Carlo methods while for applications of variance reduction techniques to kinetic
equation we remind to the recent works \cite{Hadji, Hadji1} and \cite{dimarco2, dimarco3, dimarco4}. 

In this work, we specifically tackle the issue related to the construction of efficient numerical methods for the linear transport part of a kinetic model.
In details, we consider the development of a new class of semi-Lagrangian methods with focus on the efficiency of the method and on the possibility to
extend the method developed in \cite{FKS} to arbitrarily high order polynomial reconstructions of the distribution function. Thus, the class of methods described in this paper is the natural extension
of the fast kinetic scheme (FKS) described in \cite{FKS,FKS_HO} which was based on the passage from the continuous kinetic model 
to discrete velocity models (DVM) \cite{bobylev, Pal, Pal1, Mieussens} and on a semi Lagrangian approach \cite{Chacon, Cheng, CrSon,CrSon1, Gu, Gu2, Filbet, Shoucri} 
to solve the transport part. The FKS approach was shown to be an efficient way to solve linear transport equations which 
%leads to full six dimensional simulations
has permitted to simulation the
full six dimensions
on a single processor machine. Unfortunately the solutions computed with this method were limited to a first order in space and time precision. 
This approach was, in fact, based on a piece-wise constant reconstruction of the distribution function for each velocity of the lattice.
In this paper, contrarily, we consider 
%the extension of such a scheme to 
arbitrary polynomial reconstructions which permit to get a theoretical arbitrary high order 
accurate 
method if coupled with
suitable high order splitting techniques \cite{Strang, Des} for solving the full problem given by the transport plus the collision operator. The main idea is to use a semi-Lagrangian technique
to compute the solution of the transport equation and then to compute the solution of the collision operator on the maxima and the minima of the distribution function instead of 
computing it on the grid points. This approach, as will be shown later, 
%permits 
allows 
to diminish the numerical diffusion with respect to classical semi-Lagrangian
methods 
%which employ 
employing
the same polynomial reconstruction.
Moreover this approach is well suited to preserve the structure of the solution for very long times compared to existing schemes.

The class of methods proposed can be used for discretizing different kinetic equations with different type of collision operators as it will be clarified later.
Since the focus is on the way in which the transport part is solved in the numerical test section a simple relaxation type collision operator, i.e 
the BGK (Bhatnagar-Gross-Krook) operator \cite{Gross}, has been chosen. Extensions to the full Boltzmann operator and to multidimensional settings is a work in progress.

The article is organized as follows. In the first section \ref{sec:FKS}, we present the model 
and the Fast Kinetic numerical Scheme
(FKS). In section \ref{sec:HO} we generalize the FKS to arbitrary 
high order reconstructions. A full example is described in the case of BGK kernel and 
piece-wise linear reconstructions.
Several test problem results are presented in section \ref{sec:numerics} to
assess the efficiency of this new high accurate numerical method exposed in this
paper by comparison with classical schemes. 
Conclusion and future developments are finally drawn in section
\ref{sec:conclusion}.

%------------------------------------

%------------------------------------
% WENO

%
% SECTION 2 : numerical scheme on sequential machine
%
\section{Kinetic equations and Fast Kinetic Scheme} \label{sec:FKS}

In this section we briefly present the model and the numerical method designed in \cite{FKS}. 
However we refer the reader to the cited papers for 
%more details
an exhaustive description. 

\subsection{The kinetic model}
\label{ssec:Boltzmann}

In the kinetic theory of rarefied gases, the non-negative function $f=f(x,v,t)$ characterizes the state of the system and it defines 
the density of particles having velocity $v\in\R^3$ in position $x\in\R^3$ at time $t\in\R^+$. The time evolution of the system is obtained through the equation
\be
\frac{\partial f}{\partial t} + v\cdot \nabla_x f= Q(f).
\label{eq:i1}
\ee
The operator $Q(f)$, on the right hand side in equation (\ref{eq:i1}), 
describes the effects of particle interactions and its form depends on the details of the microscopic dynamic. 
%The most well-known example is represented by the nonlinear Boltzmann integral~\cite{dimarco-review, cercignani}. 
Independently on the type of microscopic interactions considered, typically the operator characterizes the conservation properties of the physical system
\be
\int_{\R^3} Q(f)\phi(v)\,dv=0,
\label{eq:i5}
\ee
where $\phi(v)=(1,v,|v|^2)$ are commonly
called the collision invariants. We denote by
 $U(x,t)=\Int_{\R^3} f(x,v,t) \phi(v)\,dv\in\R^5$ the first three moments of the distribution function $f$.
Integrating (\ref{eq:i1}) against $\phi(v)$ yields a system of macroscopic conservation laws
\be
\frac{\partial}{\partial t}\int_{\R^3}f\phi(v)\,dv + \int_{\R^3} v \cdot \nabla_x f\,\phi(v)\,dv =0.
\label{eq:i4}
\ee
The above moment system is not closed since the second term involves higher order moments of the distribution function $f$. 
However, using the additional property of the operator $Q(f)$ that the functions belonging to its kernel satisfy \be Q(f)=0\quad\hbox{iff}\quad f=M[f],\label{eq:i6}\ee 
where the Maxwellian distributions $M[f]=M[f](x,v,t)$ can be expressed in terms of the set of moments $U(x,t)$, one can get a closed system
by replacing $f$ with $M[f]$ in (\ref{eq:i4}).
% which 
This corresponds to the set of compressible Euler equations
which can be written as
\be
\frac{\partial U}{\partial t} + \nabla_x \cdot F(U) =0,
\label{eq:i7}
\ee
with $F(U)=\Int_{\R^3} M[f]v\phi(v)\,dv$.
Note that the simplest operator satisfying (\ref{eq:i5}) and (\ref{eq:i6}) is the linear relaxation operator~\cite{Gross}
\be
Q_{BGK}(f)=\nu(M[f]-f),
\label{eq:ibgk}
\ee 
where $\nu=\nu(x,t) > 0$ defines the so called collision frequency. On the other hand, the classical Boltzmann operator reads
\be
Q_B(f)=\int_{\R^3}
\int_{\SS^2} B(|v-v_*|,\omega) \left( f(v')f(v'_*)-f(v)f(v_*) \right) dv_* d \omega,
\label{eq:bolt}
\ee
where $\omega$ is a vector of the unitary sphere
$\SS^2 \subset \R^3$. The post collisional velocities $(v',v'_*)$ are given by the
relations \be v'=\frac{1}{2}(v+v_*+|q|\omega),\quad
v'_*=\frac{1}{2}(v+v_*+|q|\omega), \ee
where $q=v-v_*$ is the
relative velocity, the kernel $B$ characterizes the
details of the binary interactions. 
In this work, the construction of the numerical scheme is first done in the general case of a collision operator which only satisfies conditions (\ref{eq:i5}) and (\ref{eq:i6}).
%Successively, 
For the numerical simulations we consider the BGK operator 
%above defined 
which has been shown to be able to model with sufficient accuracy 
many different rarefied regimes \cite{cercignani}.

\subsection{The first order Fast Kinetic Scheme}
\label{ssec:FKScheme}

The Fast Kinetic Scheme (FKS) \cite{FKS} belongs to the family 
of so-called semi-Lagrangian schemes \cite{CrSon, CrSon1, Filbet}
which are typically applied to a Discrete Velocity Model (DVM) \cite{bobylev, Mieussens} approximation of the original kinetic equation.
We briefly recall the basics of the scheme in the one dimensional setting and we remind to \cite{FKS} for details.

Let us truncate the velocity space fixing some given bounds and set a grid in velocity space of $N$ points with $\Delta v$ the
grid step. The continuous distribution function $f$ is then replaced by a vector whose components are assumed to be approximations of the 
distribution function $f$ at locations $v_k$:
\be
f_{k}(x,t) \approx f(x,v_{k},t).
\ee

The discrete velocity kinetic model consists then of a set of $N$ evolution
equations in the velocity space for $f_k$, $1\leq k \leq N$, of the form
\be
\partial_t f_{k} + v_{k} \cdot\nabla_{x} f_{k} = Q(f_k), % (\E_{k}[F]-f_{k}),
\label{eq:DM1}
\ee
where $Q(f_k)$ is a suitable approximation of the collision operator $Q(f)$ at location $k$. 
Observe that, due to the truncation of the velocity space and to the finite number of points with which $f$ is discretized, 
the moments of the discrete distribution function $f_k$ are such that
\be
\widetilde{U}(x,t) = \sum_{k}  \phi_{k} \, f_k(x,t)\, \Delta v \ne U(x,t),
\label{eq:DM}
\ee
with $\phi_k=(1,v_k,v_k^2)$ the discrete collision invariants, are no longer the ones given by the continuous distribution $f$. This problem concerns all 
numerical methods based on the discrete velocity models and different strategies can be adopted to restore the correct macroscopic physical quantities 
\cite{Mieussens,Pal,Pal1,Gamba}. In order to solve this problem we adopt 
a $L_2$ projection technique for the discretized distribution which permits
to project the discretized $f_k$ in
the space of the distributions for which the moments are exactly the continuous macroscopic quantities we aim to preserve. 
We do not enter into more details and we refer the reader to \cite{Gamba} for a description of this technique.

%======================================================
% FIG :transport
\begin{figure}
  \hspace{-2cm}
  \begin{center}
    \includegraphics[width=1.0\textwidth]{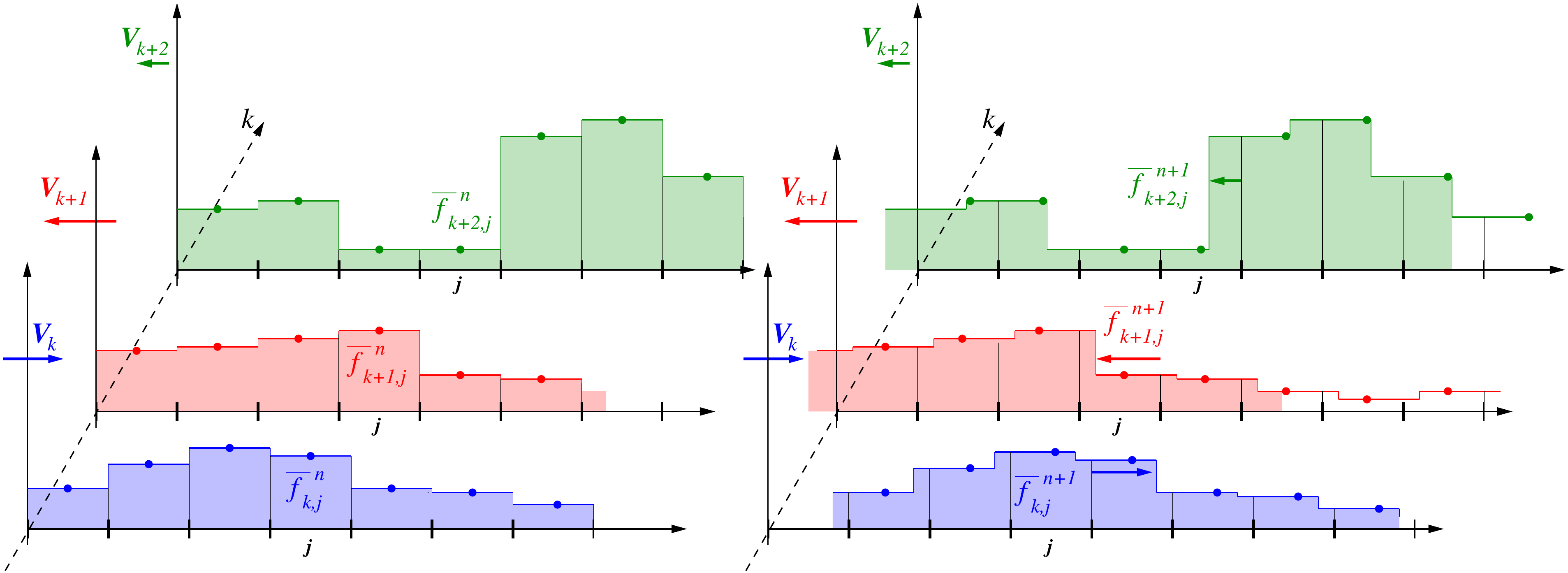}
    \caption{ \label{fig:transport}
      Illustration of the transport scheme for the first order FKS scheme.
      Left panel before transport step, right panel after transport step. Each discrete velocity (index $k$) drives its own
      transport equation with velocity $v_k$. The representation
      of $f$ is made by means of a piece-wise constant function.
      The shape of the entire function has not changed during the transport
      but the cell-centered values (bullets) may have.}
  \end{center}
\end{figure}
%===========================================================
Let us now introduce a Cartesian uniform grid in 
physical space made of $M$ points with $\Delta x$ a scalar which represents the grid step in the
physical space.
% and 
Further we define
a time discretization $t^{n+1}= t^n + \Delta t$ starting at $t^0$,
where $\Delta t$ is the time step defined by an opportune CFL condition discussed next.
We restrict ourselves to uniform meshes in space. 
Investigations are in progress to extend the present method to more complex geometries.

Each equation of system (\ref{eq:DM1}) is solved by a time splitting procedure. 
We recall here a first order splitting approach: 
first a transport step exactly solves the left-hand side, whereas a collision stage solves 
the right-hand side using the solution from the transport step as initial data:
\bea
\label{eq:transport} \textit{Transport stage}  & \longrightarrow & \partial_t f_{k} + v_{k} \cdot\nabla_{x} f_{k} = 0, \\
\label{eq:collision}\textit{Collisions stage} & \longrightarrow & \partial_t f_{k}  = Q(f_k).
\eea
\paragraph{Transport step.}
Let $f^{0}_{j,k}$ be the point-wise values at time $t^0$ of the distribution $f$,
${f}^{0}_{j,k}=f(x_{j},v_{k},t^0)$.
The idea behind the fast kinetic scheme is to solve the transport stage (\ref{eq:transport}) continuously in space, see
Fig.~\ref{fig:transport} for a sketch. To this aim we define at the initial time the function $\overline{f}_k^0(x)$ as a piece-wise continuous function for all
$x\in \Omega_j$, where spacial cell $\Omega_{j} = [x_{j-1/2}; x_{j+1/2}]$ and $\Omega=\bigcup_j \Omega_{j}$.
Hence starting from data $\overline{f}_{k,j}^0$ at time index $0$, the exact solution of (\ref{eq:transport})
is simply
\bea
\overline{f}_k^{*,1} = \overline{f}_k^0(x -v_{k} \Delta t), \quad \quad \forall x \in \Omega.
\eea
In other words, the entire function $\overline{f}_k^0$ is advected with velocity $v_k$ 
during $\Delta t$ unit of time and the $*$ superscript indicates that only the transport step has been solved so far. 
The extension of this procedure to the generic time step $n$ gives 
\bea
\overline{f}_k^{*,n+1} = \overline{f}_k^n(x -v_{k} \Delta t), \quad \quad \forall x \in \Omega , \label{eq:f_bar}
\eea
where now, the key observation is that the discontinuities of the piece-wise function $\overline{f}_k^n(x)$ do not lie on the interfaces
of two different cells. Instead, the positions of the discontinuities depend entirely on the previous advection step and thus they may be 
located anywhere in the physical space. 
This means that if only the linear transport equation has to be solved, 
this approach gives the exact solution to the equation if the initial data is effectively a piece-wise constant function 
initially 
centered on the spacial mesh.

\paragraph{Collision step.}
%We discuss now the collisions step. 
The effect of the collisional step is to change the amplitude of $\overline{f}_k(x)$.
The idea is to solve the collision operator locally on the grid points and successively extend these computed values to the full domain $\Omega$. 
Thus we need to solve the following ordinary differential equation
\be
\partial_t f_{j,k} =Q(f_{j,k}), %\frac{1}{\tau}(\E_{j,k}[F]-f_{j,k}), 
%\quad \ k=1,\ldots,N, \quad \ j=1,\ldots,M,
\label{eq:relax}
\ee
for all velocities of the lattice $k=1,\ldots,N$ and grid points $j=1,\ldots,M$.
The initial data for solving this equation is furnished by the result of the
transport step obtained by (\ref{eq:f_bar}) at points $x_j$ of the mesh at time $t^{n+1}=t^{n}+\Delta t$, i.e.
$\overline{f}^{*,n+1}_{k}(x_{j})$, for all $k=1,\ldots,N$, and $j=1,\ldots,M$.
Then, the solution of (\ref{eq:relax}), locally on the grid points, reads if, for simplicity, a forward Euler scheme is used
\be \label{eq:f_coll}
f^{n+1}_{j,k}= f^{*,n+1}_{j,k} \, + \, \Delta t \, Q(f^{*,n+1}_{j,k}),
\ee
where $f^{*,n+1}_{j,k}=\overline{f}^{*,n+1}_{k}(x_{j})$. Of course many different time integrators can be employed for solving this equation.
In particular special care is needed in the case in which the equation becomes stiff, refer to \cite{dimarco-IMEX,dimarco-Exp} for 
%different possible 
alternative strategies. Since the time integration of the collision term is not the issue considered in this paper, we considered the simplest possible scheme, the
FKS technique remains the same when other time integrators are employed. Equation (\ref{eq:f_coll}) furnishes a new value for the distribution $f$ at time
$t^{n+1}$ \underline{only in the cell centers of the spacial cells} for each velocity $v_k$. However, one needs also the value of the distribution $f$ in all points of
the domain in order to perform the transport step at the next time step. 
Therefore, we define a new piece-wise constant function $\overline{Q}_k$ for each velocity of the lattice $v_k$ as
\be
\overline{Q}^{n+1}_{k}(x)=%\overline{\E}_{k}(\bm{X},t^{n+1})=
Q(f^{*,n+1}_{j,k}),
\; \; \forall x \; \mbox{such that} \; \; \overline{f}^{*,n+1}_{k}(x)=\overline{f}^{*,n+1}_{k}(x_j). %\ j=1,\ldots,M.
\ee
% Raph
Said differently we make the fundamental assumption that the discontinuities of $Q$ coincide
with the ones of $f$.
Thanks to the above choice one can rewrite the collision step in term of
spatially reconstructed functions as
\be \label{eq:expdt}
\overline{f}^{n+1}_{k}(x) = \overline{f}_{k}(x, t^n+\Delta t) =
  \overline{f}^{*,n+1}_{k}(x)
  + \, \Delta t \, \overline{Q}^{n+1}_{k}(x).
\ee
This ends one time step of the FKS scheme.

Concerning the transport part of the scheme, the time step $\Delta t$ is constrained by a CFL condition of type
\be
\Delta t \max_{k} \left( \frac{|v_{k}|}{\Delta x} \right) \leq 1 = \text{CFL}.
\label{eq:Time}
\ee
while the time step constraint for the collision step depends on the choice of the operator $Q$. Since in the numerical test section
we used a BGK operator, the time step constraint for this part has been chosen as $\Delta t\leqslant \nu\rho$.
As observed in \cite{FKS} the transport scheme is stable for every choice of the time step, being the solution
for a given fixed reconstruction performed exactly. 
%However, being 
Nonetheless 
the full scheme being based on a time splitting technique, the error is of the
order of $\Delta t$ in the case of first order splitting or of order $(\Delta t)^q$ for a splitting of order $q$. 
This suggests to take the usual CFL condition in order to maintain the time splitting error small enough.

To conclude this section let us observe that time accuracy can be increased by high order time splitting methods,
while spacial accuracy can be increased close to the fluid limit to a nominally second-order accurate scheme by the use
of piece-wise linear reconstructions of state variables, see the details in \cite{FKS_HO}.
In the next section, we will introduce a procedure which allows to increase the spacial accuracy 
by using high order polynomial reconstructions.

%
%      H I G H     A C C U R A C Y
%
\section{High order polynomial reconstruction for the FKS scheme} \label{sec:HO}

%======================================================
% FIG :transport
\begin{figure}
  \hspace{-2cm}
  \begin{center}
    \includegraphics[width=1.0\textwidth]{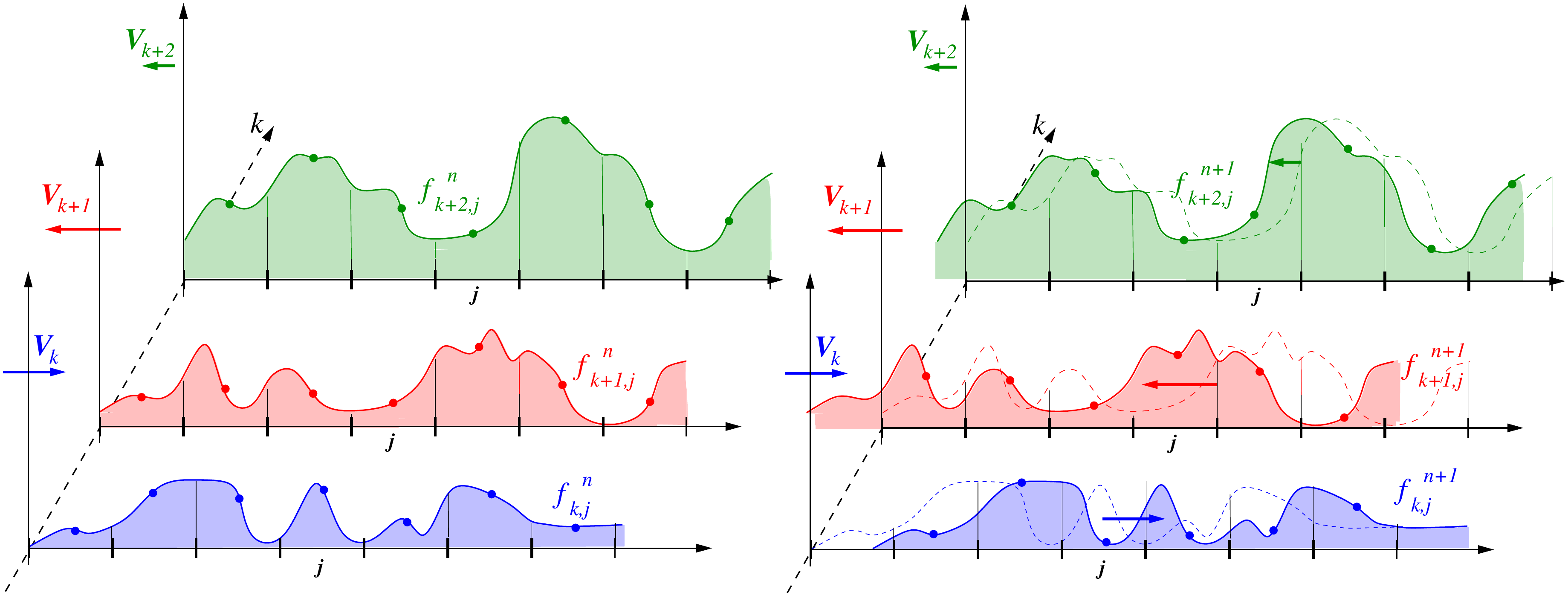}
    \caption{ \label{fig:transport2}
      Illustration of the transport scheme in 1D for the high order FKS scheme.
      Left panel before transport step, 
      right panel after transport step. Each discrete velocity (index $k$) drives its own
      transport equation with velocity $v_k$. The representation
      of $f$ is made by means of continuous function.
      The shape of the entire function has not changed during the transport
      but the cell-centered values (bullets) may have.}
  \end{center}
\end{figure}
%===========================================================

\subsection{General context for high accuracy extension}
The purpose of this section is to provide design principles to generalize the 
family of FKS schemes to higher accuracy in space. 
The FKS method employs simple piece-wise constant functions while our aim is 
to take into account arbitrary polynomial reconstruction of the distribution $f$.
The goal is to obtain a class of highly accurate in space schemes and to design a method which
will be viable in 3D$\times$3D \cite{FKS_HO2}, relatively inexpensive and independent
of the collisional operator employed: BGK, ES-BGK, Boltzmann, etc.
In this section we first present the general framework and then we introduce three possible choices which permit to design different schemes.
We conclude with a critical discussion which leads naturally to choose only one of the listed possibilities for the numerical experiments of the numerical section \ref{sec:numerics}. 
While, in the first part, the method is presented in the general context valid for any reconstruction of the distribution function, in the second part of the section 
we will give the details for the special case of piece-wise linear reconstructions. 

The main idea behind such an extension of the original FKS method 
is to use a specific local piece-wise representation of the underlying data, for instance using polynomial reconstructions following classical MUSCL like technique \cite{leveque} and
to update the reconstructed distribution in the extremes points instead of updating it as usual in the cell centers.
%By local we explicitly express the fact that we avoid the use of classical global spline
%reconstruction because of their cost in multi-dimension computations.
% knowing that
% more accurate reconstruction would be possible, only the 
% associated limiting may be more complex to design.
The starting point is, as before, the Discrete Velocity Model approximation of the original kinetic equation and a split into collision and transport terms.
Only a first order splitting is considered in this work and
the extension to higher order splitting is under study.

\paragraph{Initial data representation.}
Let, as before, $f^{0}_{j,k}$ be the point-wise values at initial 
time $t^0$ of the distribution $f$,
${f}^{0}_{j,k}=f(x_{j},v_{k},t^0)$ at cell center $x_j$ for the lattice velocity $v_k$. 
Now, we define a space piece-wise polynomial reconstruction $\widetilde{f}^0_{k}$ of a given order $p$
starting from point-wise values ${f}^{0}_{j,k}$ for each velocity of the lattice $v_k$
in cell center of $\Omega_j$.
%We will discuss later the construction and the notion of limiting for such reconstructions.
%This reconstruction is such that the Discrete Maximum Principle (DMP) whose precise definition will be  given next is satisfied. This is enough to avoid spurious oscillations in the reconstruction.
Then, as in the first order case starting from data $\widetilde{f}_{k}^0$ at time index $0$, the exact solution of (\ref{eq:transport}) is simply
\bea
f_{k}^{*,1}(x) = \widetilde{f}_{k}^0(x -v_{k} \Delta t), \quad \quad \forall x \in \Omega .
\eea
In other words, the entire piece-wise polynomial function defined by $\widetilde{f}_{k}^0$, is advected with velocity $v_k$ during $\Delta t$. 

\paragraph{Reconstruction in the collisionless case.}
If for one moment we simply consider the collisionless case, i.e. $Q(f)=0 \ \forall f$,
then the extension of this procedure to the generic time step index $n$ yields 
\bea
f_{k}^{*,n+1}(x) = \widetilde{f}_{k}^n(x -v_{k} \Delta t), \quad \quad \forall x \in \Omega, \label{eq:f_barho}
\eea
where the key point is that the new transport step restarts from where the previous one was finished. As a consequence there is no re-projection
of the distribution function $f$ on the mesh points $x_j$ to compute the new polynomial reconstruction as it is the case for classical semi-Lagrangian schemes.
Therefore, in this pure convective case, this method is exact if the initial data are exactly represented by the initial polynomial reconstruction 
employed (the one at time $0$), see figure \ref{fig:transport2} for a sketch of the transport phase in the one dimensional case when a continuous distribution function is considered. \\
  \begin{remark}
  In order to practically compute the solution of the transport step, one simply needs to store at the beginning of the computation for each
  point of the mesh a vector of size $p+1$ where $p$ is the order of the polynomial reconstruction containing the coefficients of the polynomial. 
  Successively, in order to know the value of the
  distribution $f$ everywhere, for each velocity of the lattice one needs to keep track of a single fixed point in the physical space (assuming that the transport
  velocity is constant and the mesh is uniform). This means that the computational cost associated to the solution of this step is very low compared to any kind
  of finite volume technique.
  \end{remark}

\paragraph{Reconstruction in the collisional case.}
In the case where collisions occur, their effect is to change the shape of the transported distribution in velocity space. 
Let $\widetilde{f}_{k}^{*,n+1}(x)$ be the solution known for each point $x$ in the physical space after the $(n+1)$th transport step. 
The idea behind the construction of the schemes is to replace the classical solution computed in the center of the cell (first order forward Euler scheme has been employed 
for simplicity but no restriction on the time integrators is implied), that is
\be \label{eq:f_collho}
f^{n+1}_{k}(x_j)= \widetilde{f}^{*,n+1}_{k}(x_j) \, + \, \Delta t \, Q \left( \widetilde{f}^{*,n+1}_{k}(x_j) \right),
\ee
by
\be \label{eq:f_collho1}
f^{n+1}_{k}(x_{\ext_k})= \widetilde{f}^{*,n+1}_{k}(x_{\ext_k}) \, + \, \Delta t \, \widetilde{Q}^{n+1}_{k}(x_{\ext_k}),
\ee
where $x_{\ext_k}$ are the extreme points of $\widetilde{f}^{*,n+1}_{k}(x)$ which are in general different for each $k$. 
Let observe that there exists only $M$ extremes points for each value $k$ of the velocity mesh since the reconstruction 
is performed in such a way that the the Discrete Maximum Principle is satisfied and then successively preserved since the transport phases keep
the structure of the solution unchanged. This choice permits to minimize the loss of
accuracy due to the numerical diffusion since in general the extreme points do not lie on the cell centers and thus they are lost at each time step.
Notice also that the effect of (\ref{eq:f_collho}) as well as of (\ref{eq:f_collho1}) is to compute the new value of the distribution in only one point $x_j$ or alternatively on $x_{\ext_k}$. 
This means that only the solution in the cell centers or in the extreme points have been updated with one of the two choices. However, the principle of the FKS method is to know the solution
in each point of the domain in order to continue to the next time step.
As a consequence
% and thus 
in order to compute the next transport step and advance in time, we need to define a new continuous distribution 
function, that is for all point $x$. For this we employ the same initial polynomial reconstruction used for the distribution $f$ at time $0$ using these new $M$ extrema as interpolation points. 
This ends one time step of the scheme. The only point that remains to specify is how, in practice,
are computed the functions $\widetilde{Q}^{n+1}_{k}(x_{\ext_k})$?
Three different possible choices as anticipated will discussed in the next paragraph.
They lead to three possible schemes. Only one choice is retained for the numerical simulation.

\paragraph{Computation of the collision operator in the extremes points.}
We identified three different possibilities in order to give an estimation of $\widetilde{Q}^{n+1}_{k}(x_{\ext_k})$ for the different lattice $k$ and extreme points.
\begin{enumerate}
 \item The first possibility is to observe that the distribution function is known everywhere after the transport step and so especially in $x_{\ext_k}$. 
 Thus one could simply get $\widetilde{Q}^{n+1}_{k}(x_{\ext_k})$ by direct computation from the values $\widetilde{f}^{*,n+1}_{k}(x_{\ext_k})$. 
  However, this first simple and natural choice may be very expensive when realistic
 collision operators have to be used such as the Boltzmann operator (\ref{eq:bolt}). In fact, this choice implies the collision operator to be computed $N\times M$ times, i.e. 
 the number of extreme points, at each time step. Unfortunately, this will be impracticable in multidimensional settings and with a Boltzmann like collision operator.
 \item The second possibility is to compute at the $M$ locations corresponding to the $M$ cell centers the collision operator that is 
 \be \label{eq:f_collholocal}
Q \left( \widetilde{f}^{*,n+1}_{k}(x_j) \right),\ee
 and
 then to interpolate the known values in the cell centers in order to get an estimation of $\widetilde{Q}^{n+1}_{k}(x_{\ext_k})$. This second choice will be largely
 computationally less expensive compared to the first proposed solution. However, the interpolation procedure will cause the same loss of accuracy that one gets with a classical
 semi-Lagrangian discretization since the information about the extreme points of the distribution function will not be used in this case to update the solutions; only the cell
 centers values will be used. The resulting scheme
 would perform very well in the non collisional or almost non collisional case since in this case extreme values are preserved (or almost preserved)
 by the transport phase. However this scheme would perform no better than a classical semi-Lagrangian scheme in the limit in which the number of collision becomes very large.
 \item The third possibility is to compute the collision operator 
at the $M$ locations corresponding to the $M$ cell centers as before, i.e. $Q \left( \widetilde{f}^{*,n+1}_{k}(x_j) \right)$.
Next we
% and then to 
perform the same polynomial reconstruction employed for 
%to reconstruct 
the distribution $f$ for each velocity of the lattice $v_k$. This gives 
the function $\widetilde{Q}^{n+1}_{k}(x)$, starting from the 
 point-wise values $Q(\widetilde{f}^{*,n+1}_{k}(x_j))$. The fundamental assumption for using such an approach lies on the fact that $f$, which gives the distribution of particles before 
 the collision and $Q(f)$,  which gives the distribution of particles after the collision, are supposed to evolve in a close relationship so that to share the same space topology. 
More precisely, our 
 statement is that the collision operator modifies the values $f_k$ for each velocity of the lattice, in fact the collisions redistribute the relative weights of the lattice velocity, 
 but not their spacial distribution.
 %Let observe that we define a reconstructed continuous 
 %function $\widetilde{Q}^{n+1}_{k}$ for each velocity of the lattice. This means that $\widetilde{Q}^{n+1}_{k}\ne \widetilde{Q}^{n+1}_{k'}$ with $k\ne k'$ even if the same polynomial reconstruction is
 %employed, since the points used for the reconstruction are different for each velocity of the lattice as it will be made clear in the following. This means that since $f_k$ can be 
 %interpreted as the number of particles with velocity $v_k$, and since the collision only changes the number of particles with a given velocity: in fact in general 
 %$f^{n+1}_{k}(x_j)\ne \widetilde{f}^{*,n+1}_{k}(x_j)$; then our hypothesis can be simply interpreted by saying that the particles with a given velocity maintain the same space 
 %topology over the collision. 
 This last solution, which will be the one 
%employed next
implemented and tested, 
%it permits to 
preserves the structure of the solution as in the first possibility but at a cost comparable to that
 of the second possibility. 
In fact only $M$ evaluations of the collision operator are used in this case.
\end{enumerate}
In order to detail the proposed approach (the third case) we need to define the polynomial reconstruction and 
%consequently to define a 
an associated 
limiting procedure to ensure that new extrema are not generating spurious oscillations. 
In the following, 
for the sake of clarity, 
we choose a simple piece-wise linear reconstruction of the distribution function along with 
\textit{ad hoc} limiters. 
Moreover we restrict ourselves to the simple case of BGK collisional operator. 
Note that the detailed construction of FKS schemes in multidimensions 
with higher order polynomial reconstructions 
for Boltzmann type of operator is the subject of a future work.

\subsection{Linear polynomial reconstruction case}
%======================================================
% FIG :transport
\begin{figure}
  \hspace{-2cm}
  \begin{center}
    \includegraphics[width=1.0\textwidth]{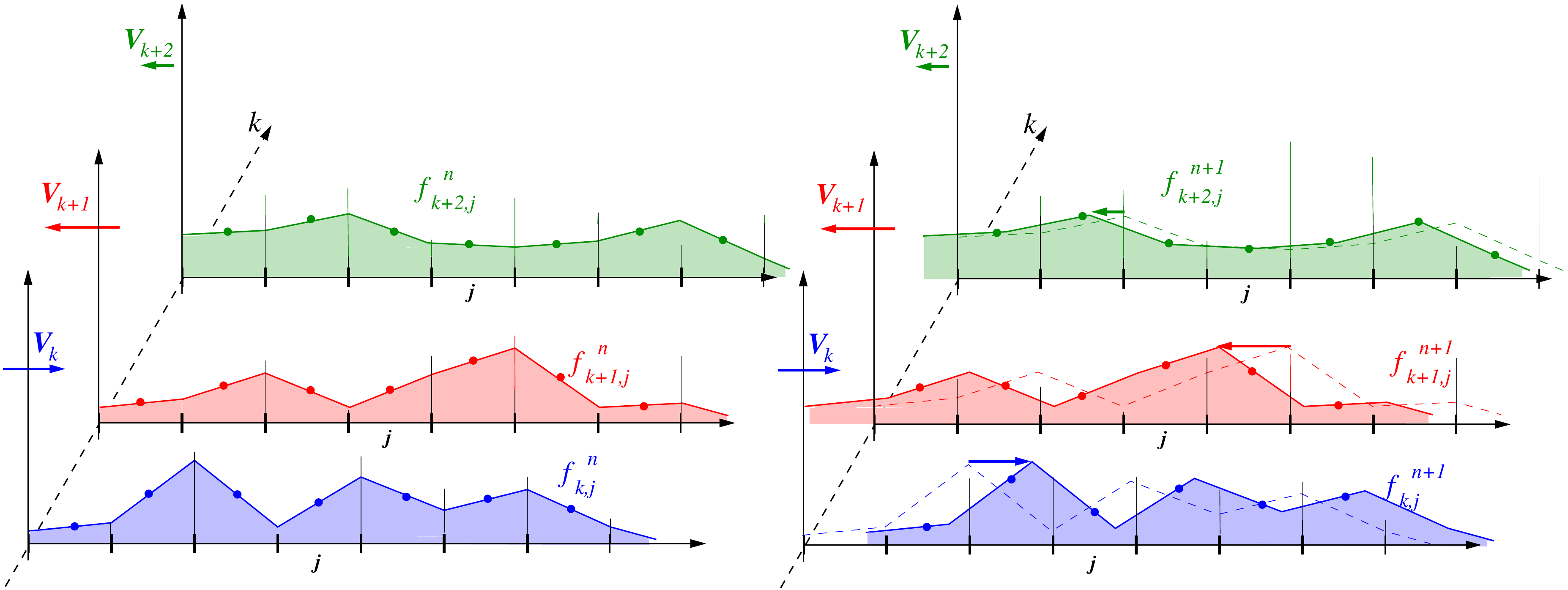}
    \caption{ \label{fig:transport3}
      Illustration of the transport scheme in 1D for the high order FKS scheme with linear reconstruction.
      Left panel before transport step, 
      right panel after transport step. Each discrete velocity (index $k$) drives its own
      transport equation with velocity $v_k$. The representation
      of $f$ is made by means of continuous function.
      The shape of the entire function has not changed during the transport
      but the cell-centered values (bullets) may have.}
  \end{center}
\end{figure}
%===========================================================

\begin{figure}
  \begin{center}
  \begin{tabular}{c}
  \hspace{-1.cm} 
    \includegraphics[width=1.0\textwidth]{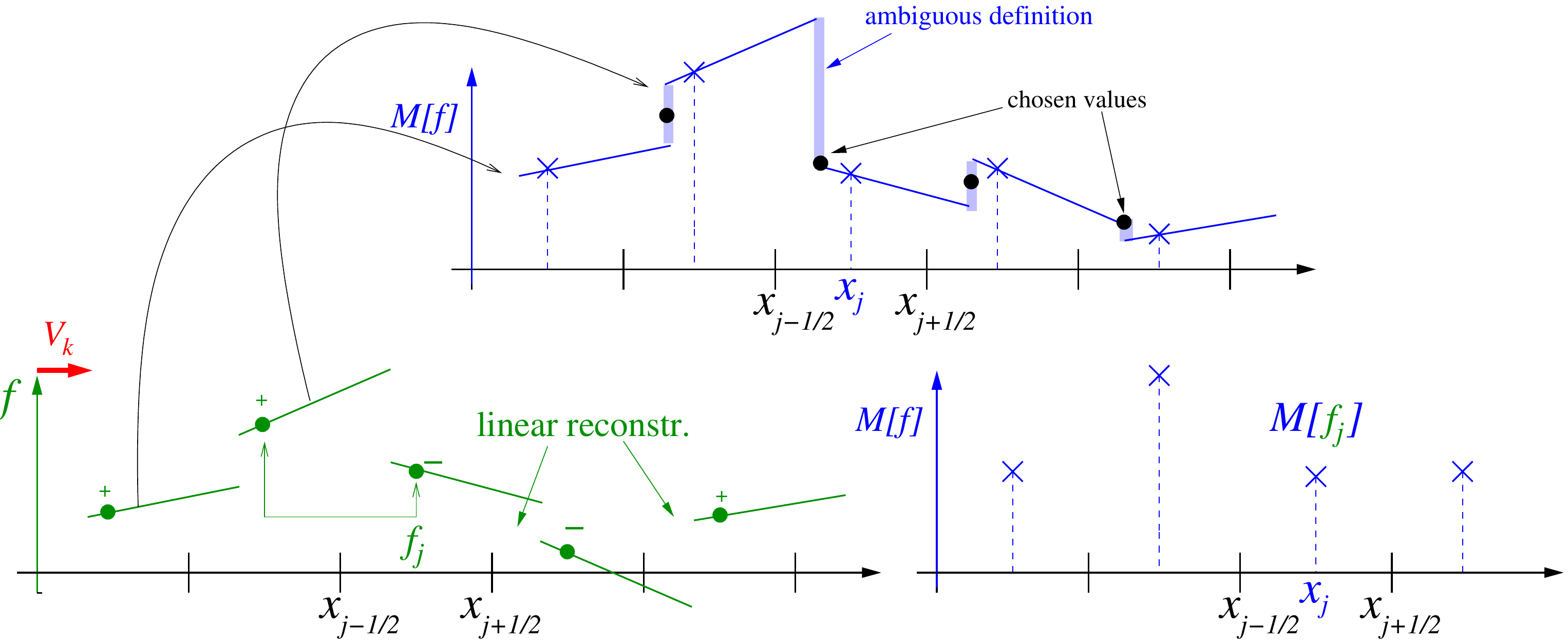}
  \end{tabular}
    \caption{ \label{fig:limiting}
      Illustration of the reconstruction technique employed in the R-FKS ---
    From the distribution function data one reconstructs an upwind piece-wise linear polynomial 
    (green on the left panel). 
    From the distribution functions known at the cell centers we evaluate the 
    Maxwellian also at cell centers (blue cross). 
    Last on the top panel, we employ the
    same slope from the distribution function which passes through the Maxwellian cell center
    value. An ambiguous definition at each node is clarified by choosing the less extreme values
    in case the slopes are of different signs and the mean value otherwise, see the black bullets.}
  \end{center} 
\end{figure}

In this section we will 
focus on a FKS method which employs a polynomial reconstructions 
of degree one.  The transport phase has been detailed in the previous paragraphs, it can
however  
be summarized by the sketch in figure \ref{fig:transport3} for the specific case under study.
The collisions are modeled by a relaxation towards the thermodynamics equilibrium $Q_{BGK}=\nu(M[f]-f)$.
In this specific case, being the transport phase unaltered, the collision phase becomes on the grid points 
\be \label{eq:collBGK}
f^{n+1}_{k}(x_j)= \left[ \exp(-\nu\Delta t) \right] \, f^{*,n+1}_{k}(x_j)+
\left[1-\exp(-\nu\Delta t) \right] \, \E_{k}[U_j^{*,n+1}],
\ee
where $\E_{k}[U_j^{*,n+1}]$ is a consistent discretization of the Maxwellian distribution at the grid points. 
% which
This consistent discretization only depends on the macroscopic 
%values
moments $U(x_j,t^{*,n+1})=U_j^{*,n+1}$ also defined on the grid points and 
satisfies the conservation of macroscopic quantities %property
\be
U(x_j,t^{*,n+1}) 
\, = \, \sum_{k}  \phi_{k} \,\E_{k}[U_j^{*,n+1}] \Delta v
\, = \, \sum_{k}  \phi_{k} \, f^{*,n+1}_{k}(x_j)\Delta v .
\ee
In the previous equation we have 
%and where we 
used the exact solution of equation $\partial_t f = Q_{BGK}$ in this case instead of the 
solution of the forward Euler scheme as before. 
%, since it was available.
% the pointwise macroscopic values being computed by summing the local value of the discrete
% distribution $f$ over the velocity set after the transport step
% \be
% \sum_{k}  \phi_{k} \, f^{*,n+1}_{k}(x_j)\Delta v=U(x_j,t^{n+1})=U(x_j,t^{*,n+1}),
% \ee
% where the last equality holds since the relaxation phase conserves the macroscopic quantities.
Equation (\ref{eq:collBGK}) furnishes the new value of the distribution $f$ at time
$t^{n+1}$ only in the cell centers of the spacial cells for each velocity $v_k$,
see the green bullets in figure~\ref{fig:limiting}.
Then, since our aim is to update the distribution function at the extreme points $x_{\ext_k}$, we need to evaluate the Maxwellian distribution in these points.
In order to do that we use the same first order polynomial interpolation used for $f$ but for the Maxwellian distribution, let observe the piece-wise linear green line in figure~\ref{fig:limiting}.
These slopes are then transferred onto the Maxwellian distribution to create a continuous version of $M$ from cell centered data 
(blue crosses in figure~\ref{fig:limiting}). However, this procedure causes loss of uniqueness in the definition of the Maxwellian state at the extreme points, in fact for each of the
extreme points we get a left state $\E_k^{-}(x_{\ext_k})$ and a right $\E_k^{+}(x_{\ext_k})$ state, see figure~\ref{fig:limiting}.
In order to restore uniqueness and to avoid the creation of extrema which may be possibly premises of spurious oscillations, the value of the Maxwellian function at these points is defined
by taking the average between the left and the right states when 
the slopes on both sides have the same sign. This average is weighted by the distance from the cell center.
On the other hand, when the slopes of the left and right states are of opposite
signs, then the value of the equilibrium function is fixed as the minimum between the left and the right states if the left slope is positive (presence of a maximum), or as the maximum 
%between the left and the right state
if the left slope is negative (presence of a minimum). For all cases, this value is denoted by $\E_k[U_j^{*,n+1}](x_{ext_{k}})$.
%======================================================
% FIG :transport
%===========================================================
Finally the new distribution $f^{n+1}_{k}$ at any extreme point $x_{\ext_k}$ is defined by summing the values of 
the transported distribution and of the reconstructed Maxwellian distribution at the $M$ points at which the transported
distribution function reaches its extrema, that is
\be \label{eq:f_BGKho}
f^{n+1}_{k}(x_{\ext_k})= \left[  \exp(-\nu\Delta t)\right] \, f^{*,n+1}_{k}(x_{\ext_k})+
\left[ 1-\exp(-\nu\Delta t) \right] \, \E[U^{n+1}]_{k}(x_{\ext_k}).
\ee
From the knowledge of these new extreme values, we redefine the new reconstructed distribution $f$ by the same first order polynomial interpolation employed
at time $t^n$ using the points $x_{\ext_k}$ as interpolation nodes. 
Consequently we are ready to compute another transport phase and 
this ends one time step of the scheme.

%------------------------------------

%------------------------------------
% NUMERICS

\section{Numerical experiments}
\label{sec:numerics}
In this section we test the FKS scheme with piece-wise linear reconstructions against the original FKS schemes 
(referred to as R-FKS and FKS respectively) and some
classical semi-Lagrangian (SL-) schemes of first order (referred to as 'SL-Upwind') and 
second order of accuracy (referred to as 'SL-MUSCL').
The SL-MUSCL scheme employs a simple piece-wise linear reconstruction with van Leer flux limiter 
\cite{leveque}. Note that most of semi-Lagrangian second-order limited schemes would produce rather equivalent results.

The methodology of verification is revolving around three test cases:
\begin{itemize}
\item Numerical convergence test. Given a smooth solution for different value of $\nu$,
  we observe that the numerical method is nominally second order accurate.
\item Riemann problem using Sod like initial data
  \begin{itemize}
  \item for $\nu=10^{4}$ which consists of an almost pure hydrodynamics situation where genuine shock wave
    and discontinuous solutions may appear;
  \item for $\nu=10^{3}$ which is a transition regime where discontinuities are no more present but
    steep fronts are still present;
  \item for $\nu=10^{2}$ which is the entrance door towards the full collisional regime for which 
    physical diffusion dominates any other non-linear waves.
  \end{itemize}    
\item Highly oscillating initial data which are generating a large number of discontinuities 
  even in the kinetic regime $\nu=10^{-2}$. This problem
  requires a low dissipative scheme to maintain accurate results. For this test we will 
  present mesh convergence and efficiency studies.
\end{itemize}
One expects, and we will show that it is indeed the case, that the new FKS scheme with reconstruction (R-FKS) is more accurate
 than FKS and the SL-Upwind schemes, and has the same -or possibly better- accuracy than SL-MUSCL scheme.
At the same time one expects that the R-FKS scheme be computationally less expensive than the SL-MUSCL
scheme.
In an ideal situation the cost of R-FKS should be close to the FKS scheme, which is already known
to be an inexpensive scheme compared to classical SL strategy.
Finally, as expected, we will show that for a fixed error the R-FKS demands less computational resources than SL-MUSCL, 
FKS and SL-Upwind in 1D$\times$1D assuming that, if so, 
this situation will be vastly more interesting in 
2D$\times$2D and 3D$\times$3D computations. 
In all simulations we set the CFL number to one for all schemes.

\subsection{Convergence test problem} \label{sec:convergence}
\begin{figure}[ht!]
\begin{center}
  \hspace{-0.5cm}
  \includegraphics[width=0.36\textwidth]{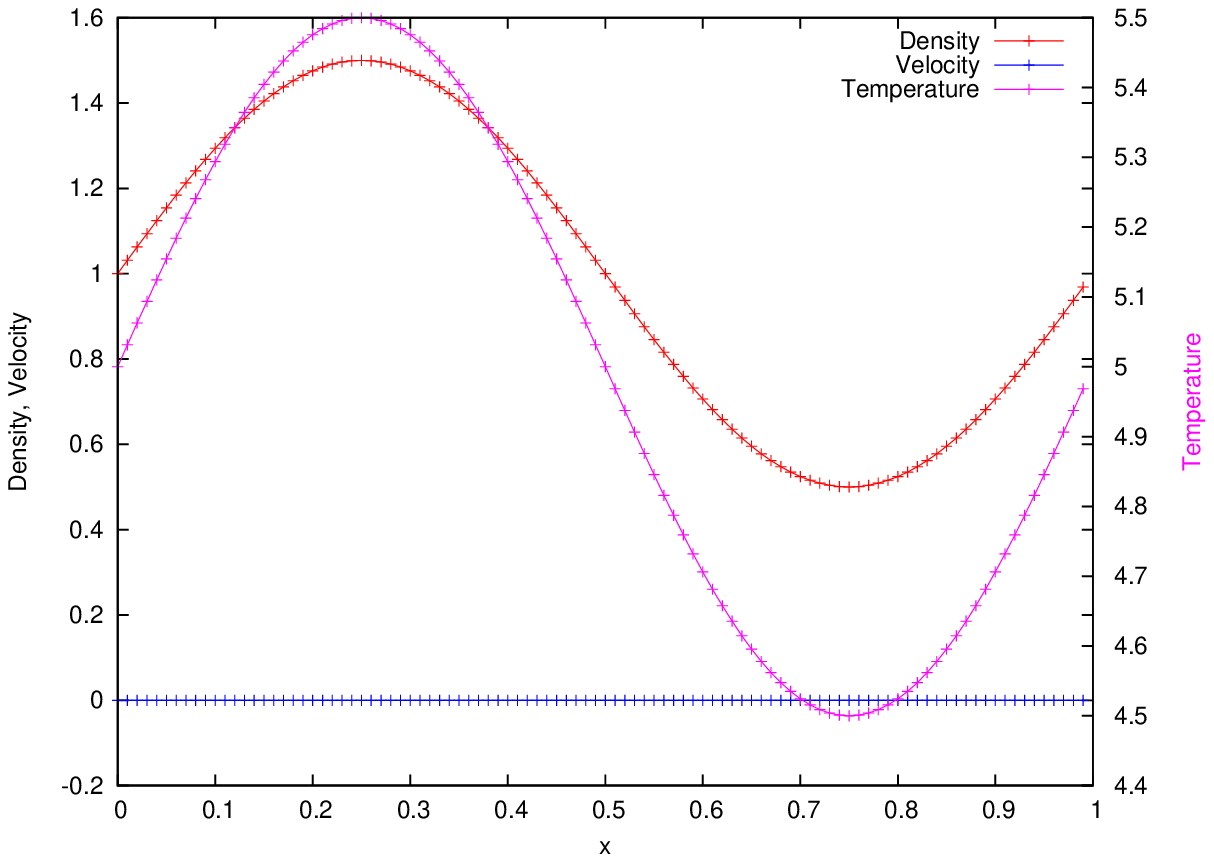} 
  \begin{tabular}{ccc}
    \hspace{-0.75cm}
    \includegraphics[width=0.36\textwidth]{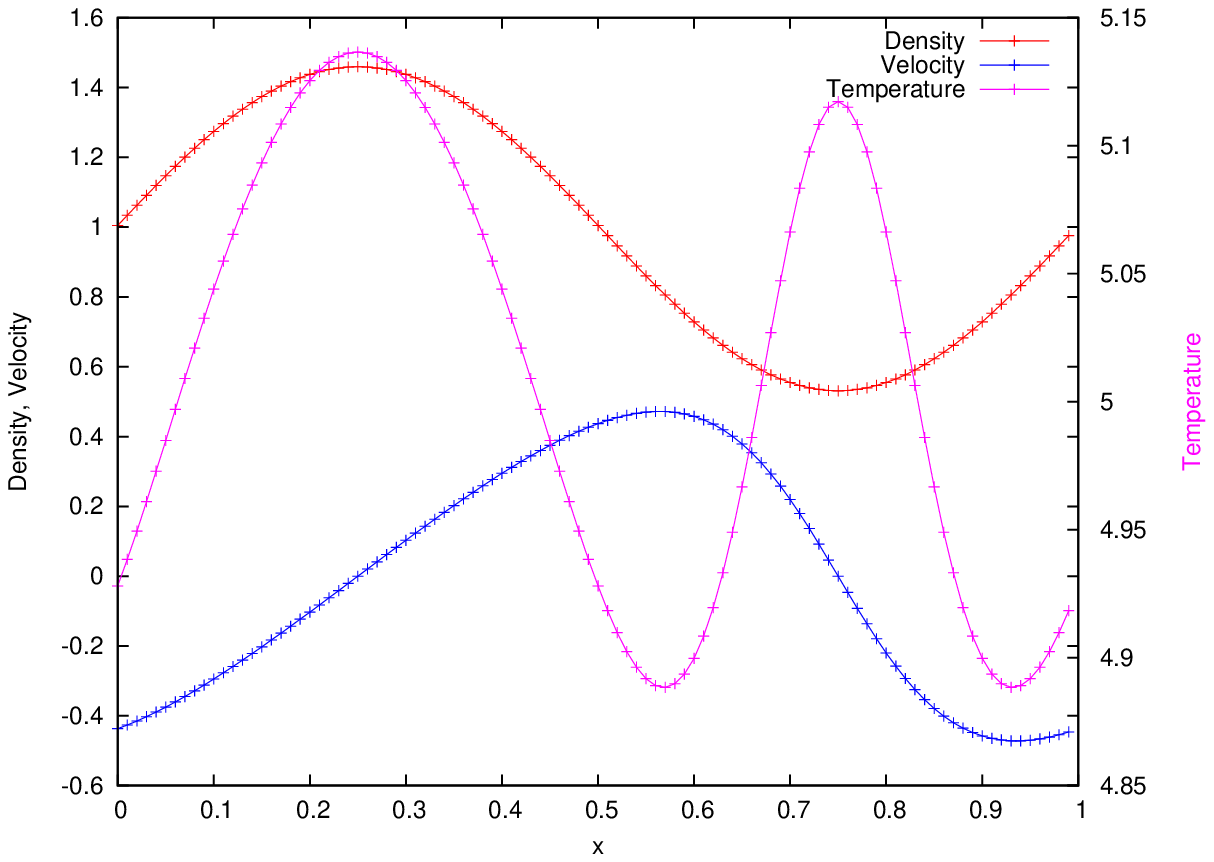}&
    \hspace{-0.5cm}
    \includegraphics[width=0.36\textwidth]{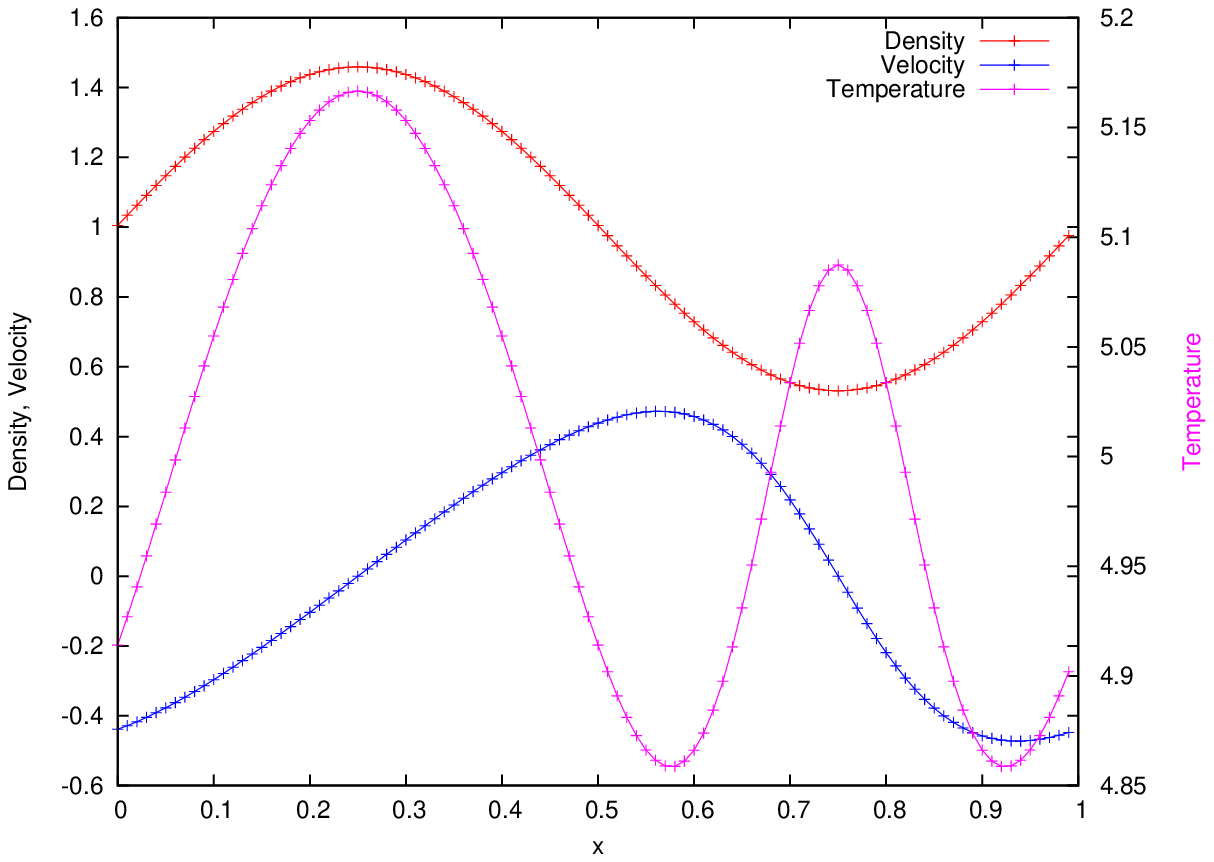}&
    \hspace{-0.5cm}
    \includegraphics[width=0.36\textwidth]{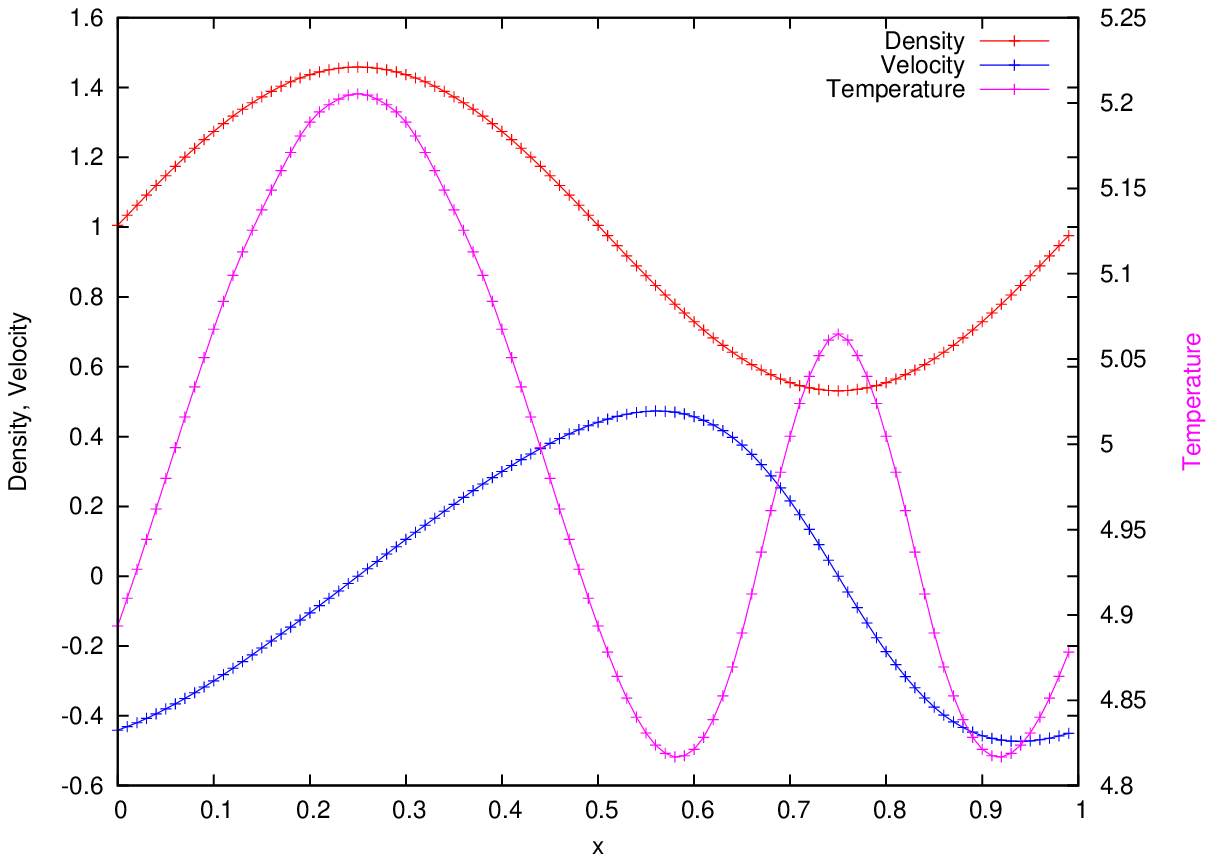}
  \end{tabular}
\caption{
    Convergence test case results for $M_1=100$ cells --- 
    Top panel: initial data ---
    Left to right panels: final solution for $\nu=10^1, 10^2, 10^4$ .}
\label{fig:cvg}
\end{center}
\end{figure}

% Problem description
This problem is initialized with a smooth macroscopic density, velocity and temperature data
\bea
\rho(x,t=0) = 1 + \frac12\sin( 2 \pi x), \quad
u(x,t=0) = 0, \quad
T(x,t=0) = 5 + \frac12\sin( 2 \pi x).
\eea
where $\Omega=[0,1]$ and an initial thermodynamics equilibrium is chosen, i.e. $f(x,v,t=0)=M[f](x,v,t=0)$.
The final time is set to $t_{\text{final}}=0.025$.
The initial truncation of the velocity space is $L_v=[-15,15]$. 
The meshes are uniform both in velocity and physical spaces with 
$N=50$ and $M_j=100 j$ cells with $j=1,2,4,8,16,32$ and $64$ respectively. \\
% Figures and results
In figure~\ref{fig:cvg} are presented the initial data (top panel) and
the final solution for $\nu=10^1, 10^2, 10^4$ from left to right.
Density, velocity and temperature are displayed when $M_1=100$ cells are considered.
%=== B E G I N    F I G U R E ================================================
%=== E N D   F I G U R E ================================================
% Error computation description
According to our initialization the cell centers for mesh $M_1$ are also cell centers
for any mesh $M_i$, $1 \leq i\leq 64$. Consequently we can compute 
an estimation of the numerical order of convergence using three consecutive meshes 
of indices $j,j+1, j+2$ as
\bea \label{eq:order_cvg}
p_j =\log_2\left( \frac{ \Sum_{i=1}^{M_1} | \rho_i^{j} - \rho_{2i}^{j+1} | }{\Sum_{i=1}^{M_1} | \rho_{2i}^{j+1} - \rho_{4i}^{j+2} |}  \right) ,
\quad
j=1,2,4,8,16.
\eea
In order to properly measure the error due to the space discretization, 
we make the time discretization error very small by fixing the time step $\Delta t$
to the value needed for the finest mesh, that is 
$\Delta t=(\Delta x)_{64}/\max(|v|) \simeq 10^{-5}$, where $(\Delta x)_{64}=1/6400$.
Moreover, in order to observe the convergence rate for different regimes, we consider three different collision frequencies :  $\nu=10^{1}$, $10^{2}$ and $10^{4}$. 
%== T A B L E =========================================================
\begin{table}
\begin{center}
  \begin{tabular}{|c||c|c|c|}
    \hline    
    \textbf{Meshes}  & \multicolumn{3}{|c|}{\textbf{Order of convergence}} \\ 
    \cline{2-4}
    $M_j,M_{j+1},M_{j+2}$  & $\nu=10^1$  &  $\nu=10^2$ & $\nu=10^4$ \\
%    \cline{2-4}
%    $j$  & $p_j$   &  $p_j$ & $p_j$ \\
    \hline   
    $1,2,4$    & $1.010$   & $1.056$  & $1.776$      \\
    $2,4,8$    & $9.404$   & $6.695$  & $4.569$      \\
    $4,8,16$   & $2.711$   & $2.012$  & $2.584$      \\
    $8,16,32$  & $2.658$   & $2.043$  & $2.556$      \\
    $16,32,63$ & $1.968$   & $1.948$  & $2.357$      \\
    \hline
  \end{tabular}     
  \caption{ \label{tab:convergence}
    Convergence test case for a smooth solution computed from
    three successively refined meshes $M_j,M_{j+1},M_{j+2}$ via equation (\ref{eq:order_cvg})
    and for three collisional frequencies.    
    % Results obtained by the code in Raphael's machine
    % CODE_RFKS_CONVERGENCE_SINE_17_02  
  }  
\end{center}
\end{table}
%=== E N D   T A B L E ================================================
% Comments
The results are gathered in table~\ref{tab:convergence} and we observe that, 
for the three collisional frequencies considered, the nominal orders of accuracy are
approaching $2$. The first two orders of accuracy that can be computed
present some eratic behaviors because the numerical solution is not yet mesh-converged.

\subsection{Riemann problem using Sod like initial data}
\label{sec:sod-shock}

% Problem description
In this paragraph we study the classical Riemann problem with Sod like initial data
\bea
\nonumber
 \int_\RR f dv=1,    \; \int_\RR fv dv=0, \; \int_\RR fv^2dv=2.5, & \;\;  \; \text{if} \;\;& x\leq L/2,  \\
\nonumber
 \int_\RR f dv=0.125 \; \int_\RR fv dv=0  \; \int_\RR fv^2dv=0.25,& \;\;  \; \text{if} \;\;& x> L/2,  
\eea
with $\Omega=[0,1]$ the domain length and initial thermodynamics equilibrium, i.e. $f(x,v,t=0)=M[f](x,v,t=0)$.
The CFL condition is fixed for all different tests to $\Delta t=\Delta x/\max(|v|)$,
where $\max(|v|)$ corresponds to the maximum speed of the particles which is fixed by the initial truncation of the velocity space $L_v=[-20,20]$. 
The same initial data are used to run three different problems corresponding to three different collision frequencies :  $\nu=10^{2}$, $10^{3}$ and $10^{4}$. 
The meshes are uniform both in velocity and in the physical space and the number of points is chosen respectively as $N=50$ and $M=300$. 

\begin{figure}[ht!]
\begin{center}
  \begin{tabular}{cc}
    \hspace{-1.0cm}
    \includegraphics[width=0.5\textwidth]{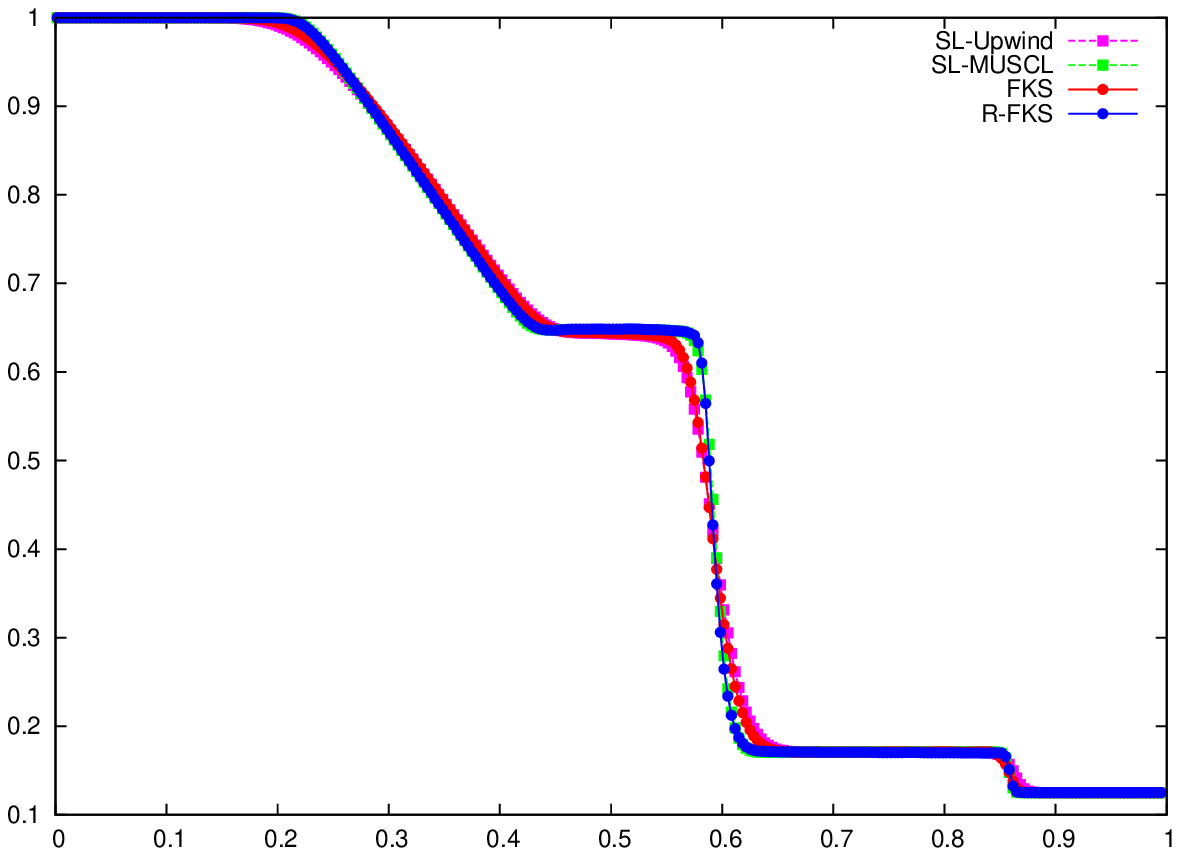}&
    \hspace{-0.35cm}
    \includegraphics[width=0.5\textwidth]{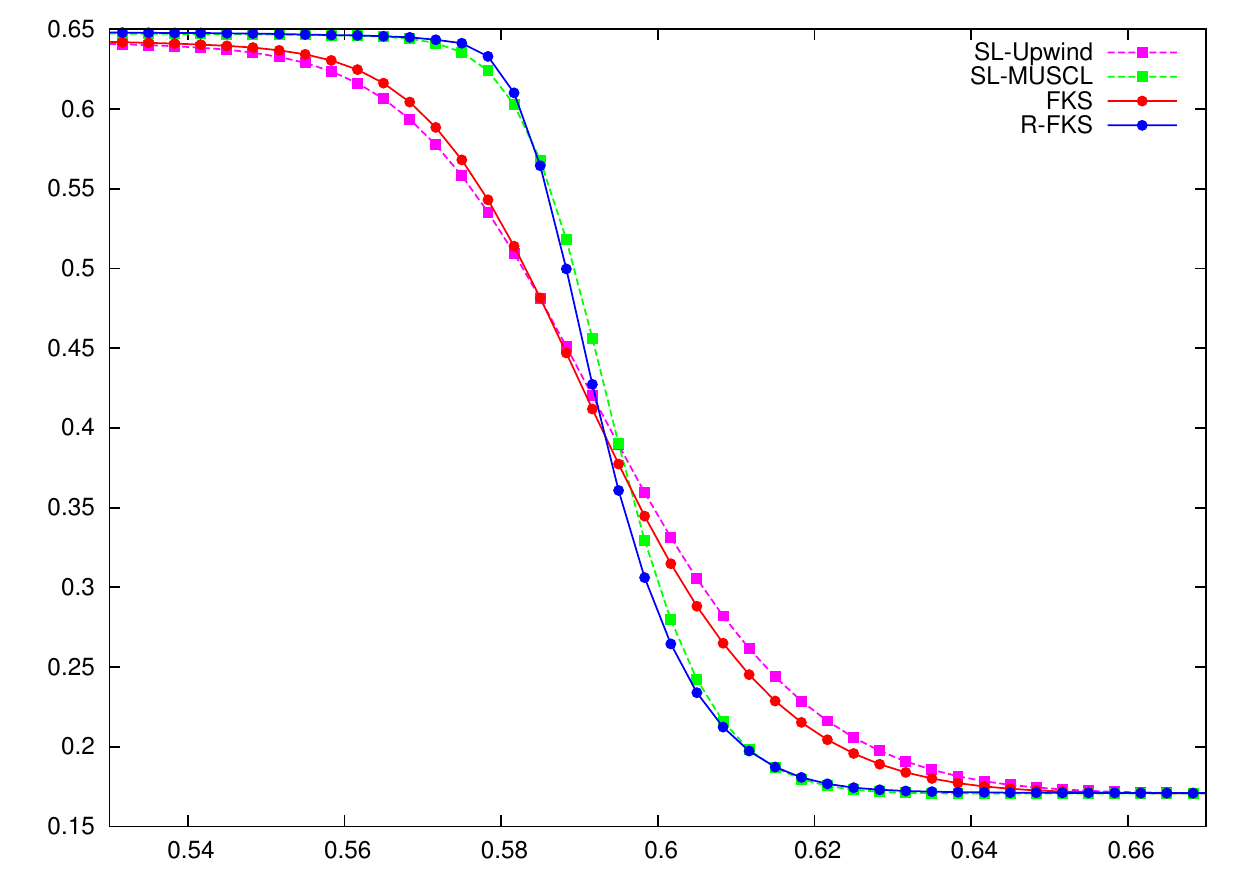} \\
    \hspace{-1.0cm}
    \includegraphics[width=0.5\textwidth]{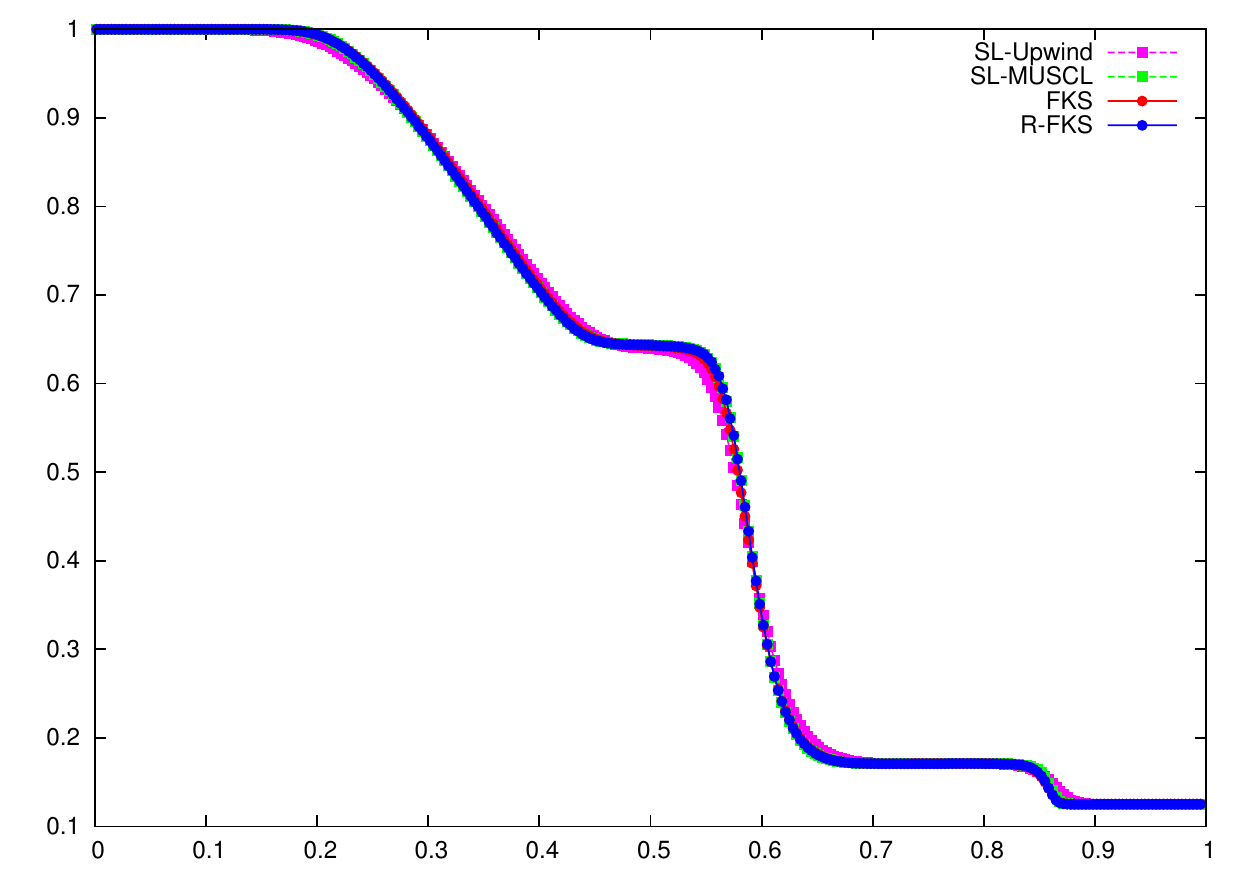}&
    \hspace{-0.35cm}
    \includegraphics[width=0.5\textwidth]{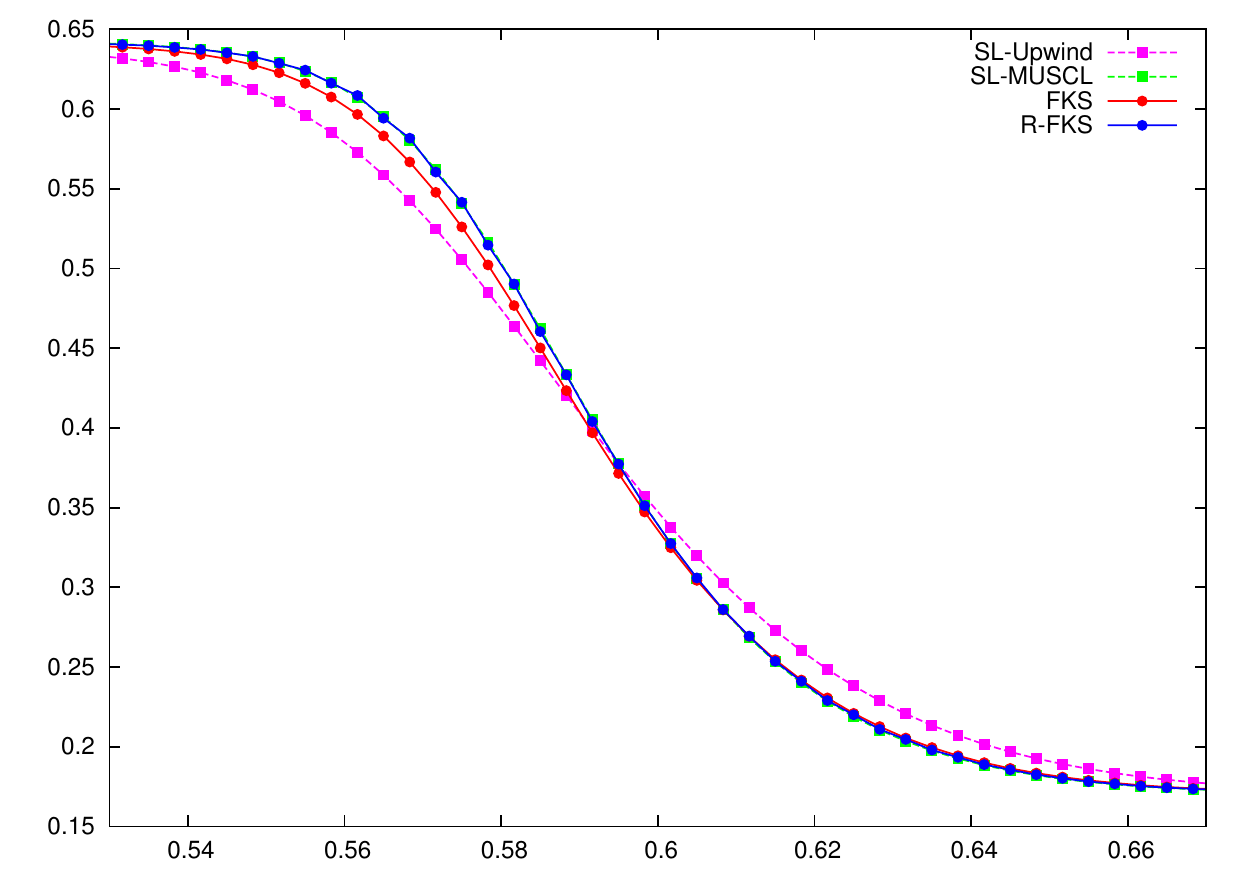} \\
    \hspace{-1.0cm}
    \includegraphics[width=0.5\textwidth]{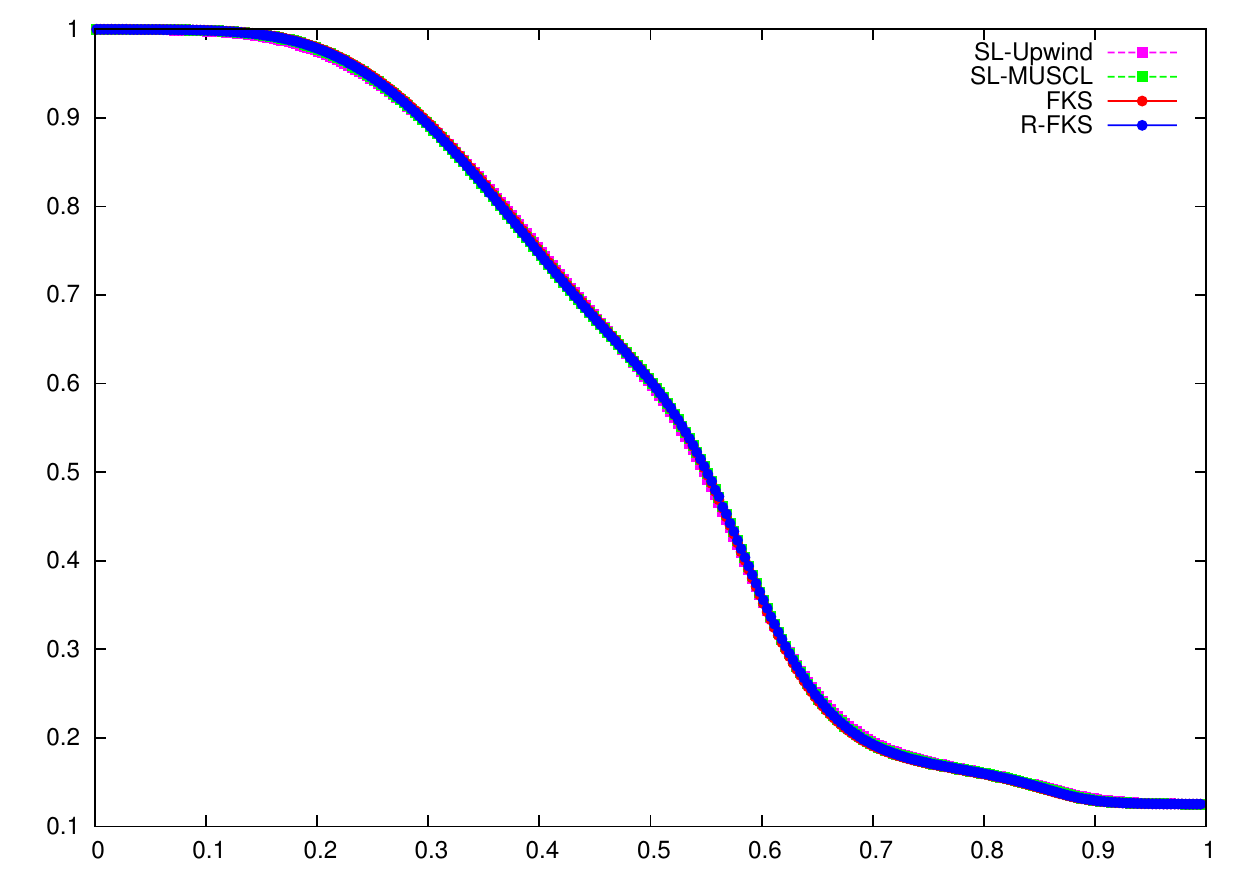}&
    \hspace{-0.35cm}
    \includegraphics[width=0.5\textwidth]{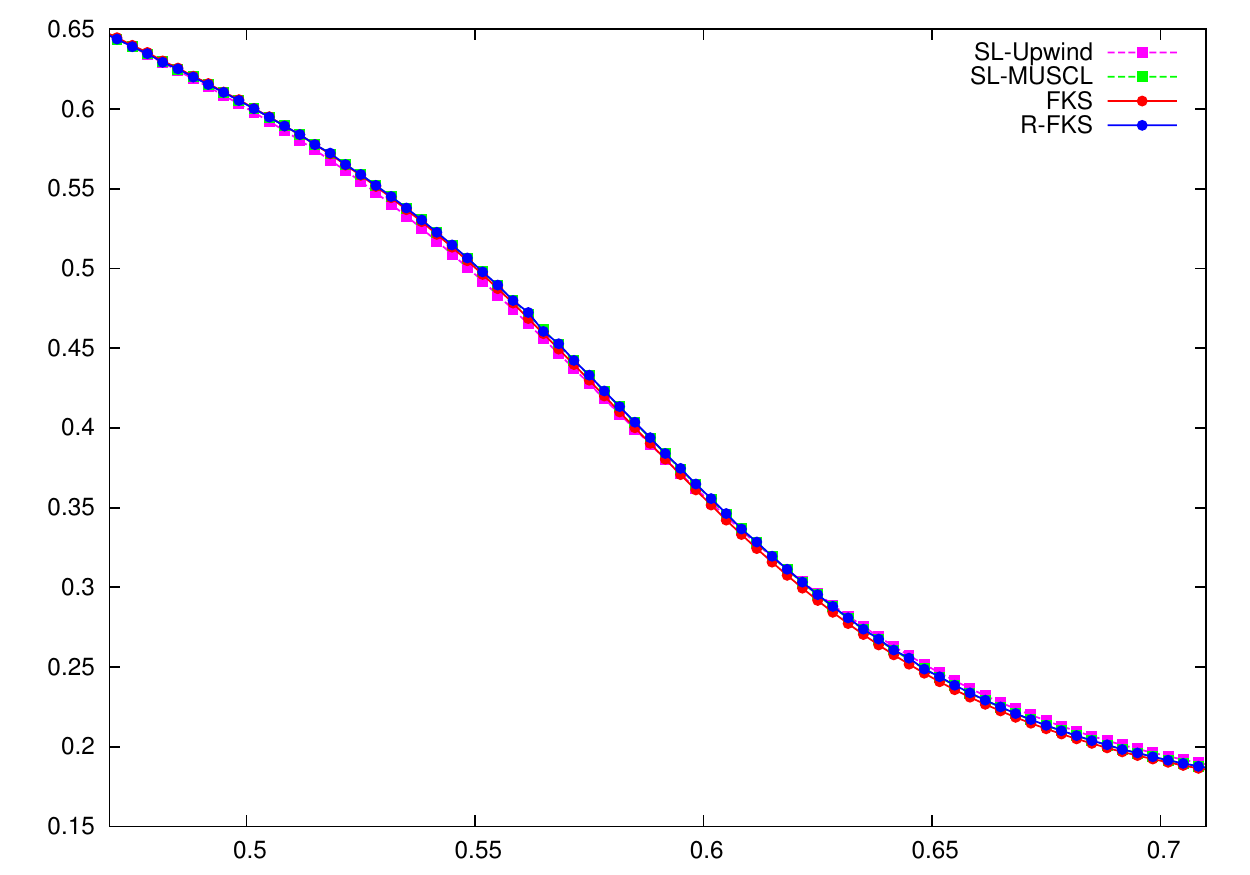} \\
  \end{tabular}
\caption{Riemann problem --- Solution at $t_{\text{final}}=0.07$ for the
density. Left full solution, right zoom close to the region in which the solution develops a contact discontinuity in the limit of infinite
collisions. Collision frequency $\nu=10^{4}$ (top), $10^{3}$ (middle), $10^{2}$ (bottom).}
\label{fig_Sod1}
\end{center}
\end{figure}

\begin{figure}[ht!]
\begin{center}
  \begin{tabular}{cc}
    \hspace{-1.0cm}
    \includegraphics[width=0.5\textwidth]{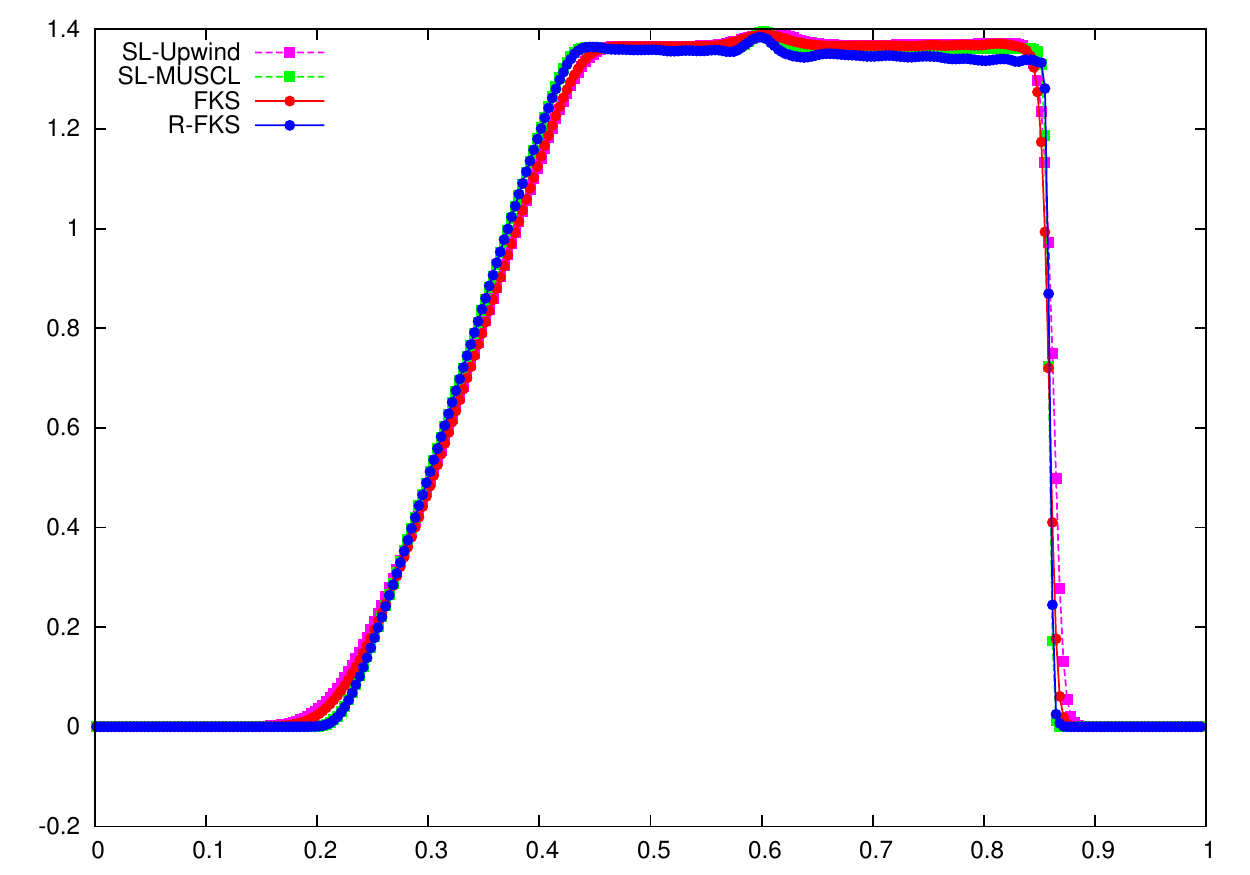}&
    \hspace{-0.35cm}
    \includegraphics[width=0.5\textwidth]{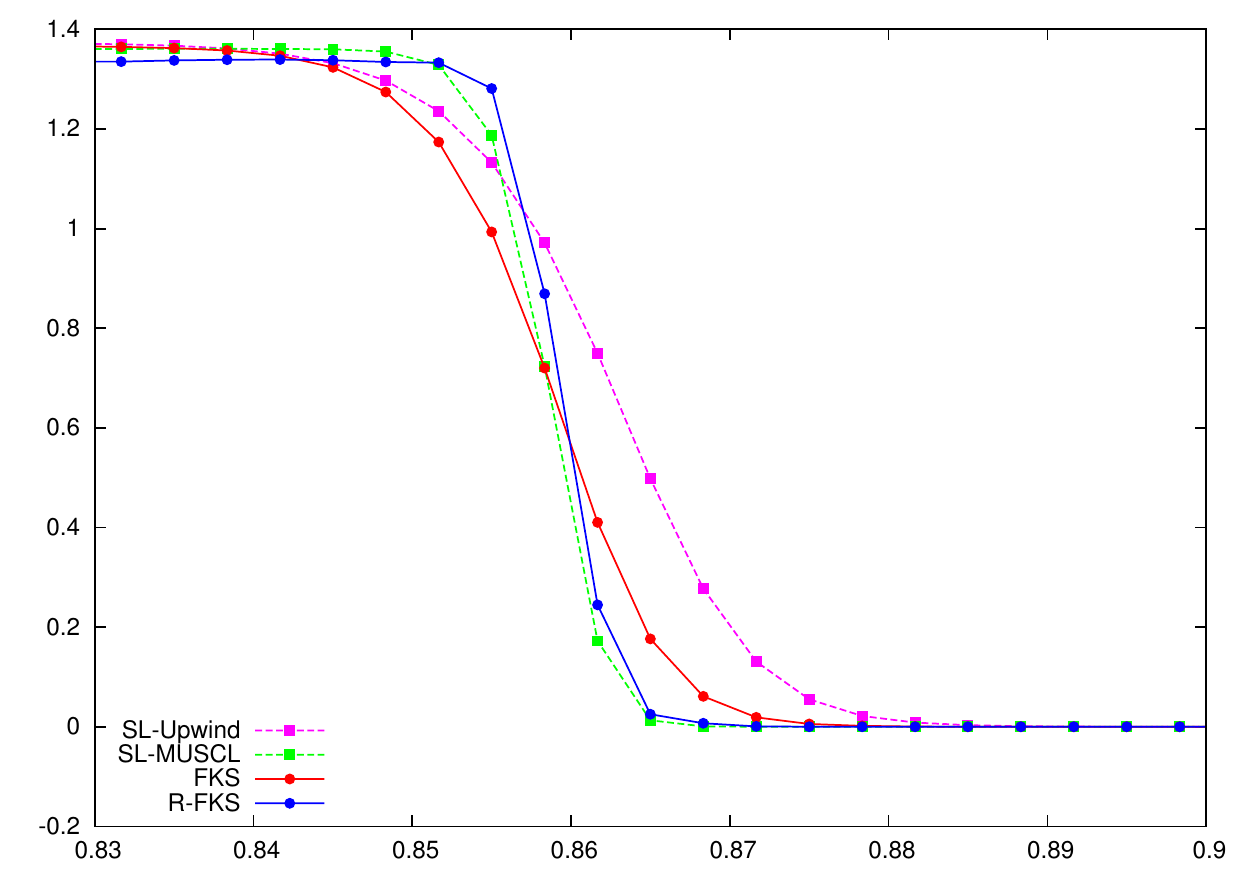} \\
    \hspace{-1.0cm}
    \includegraphics[width=0.5\textwidth]{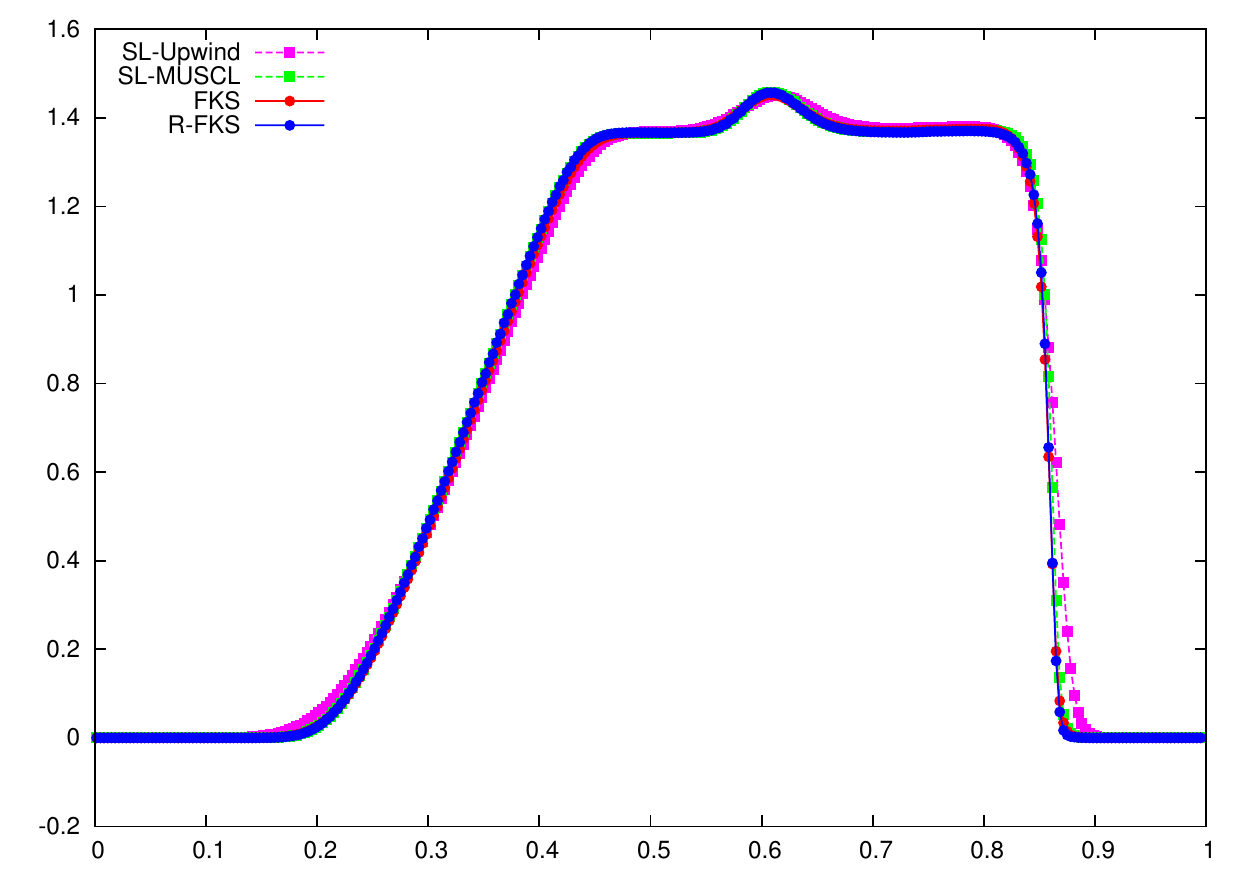}&
    \hspace{-0.35cm}
    \includegraphics[width=0.5\textwidth]{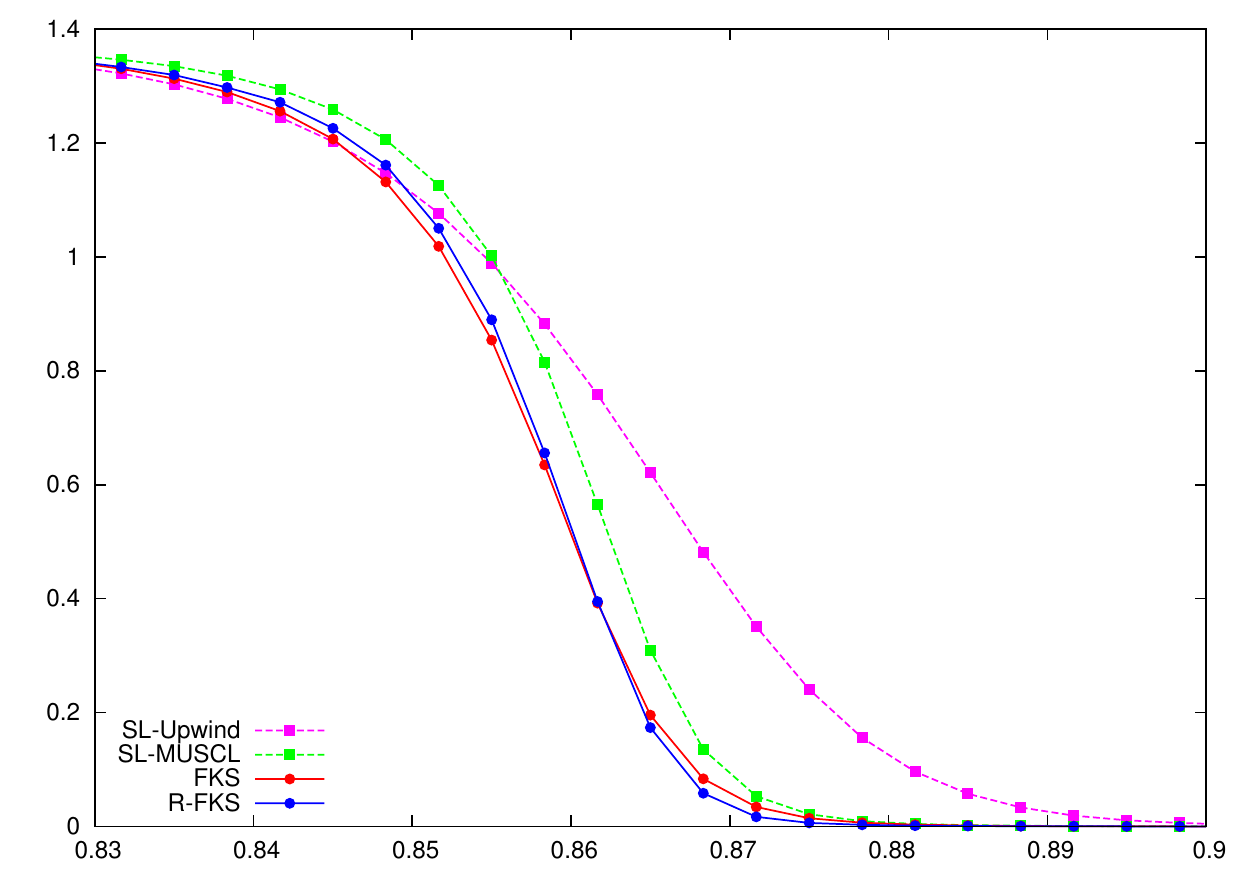} \\
    \hspace{-1.0cm}
    \includegraphics[width=0.5\textwidth]{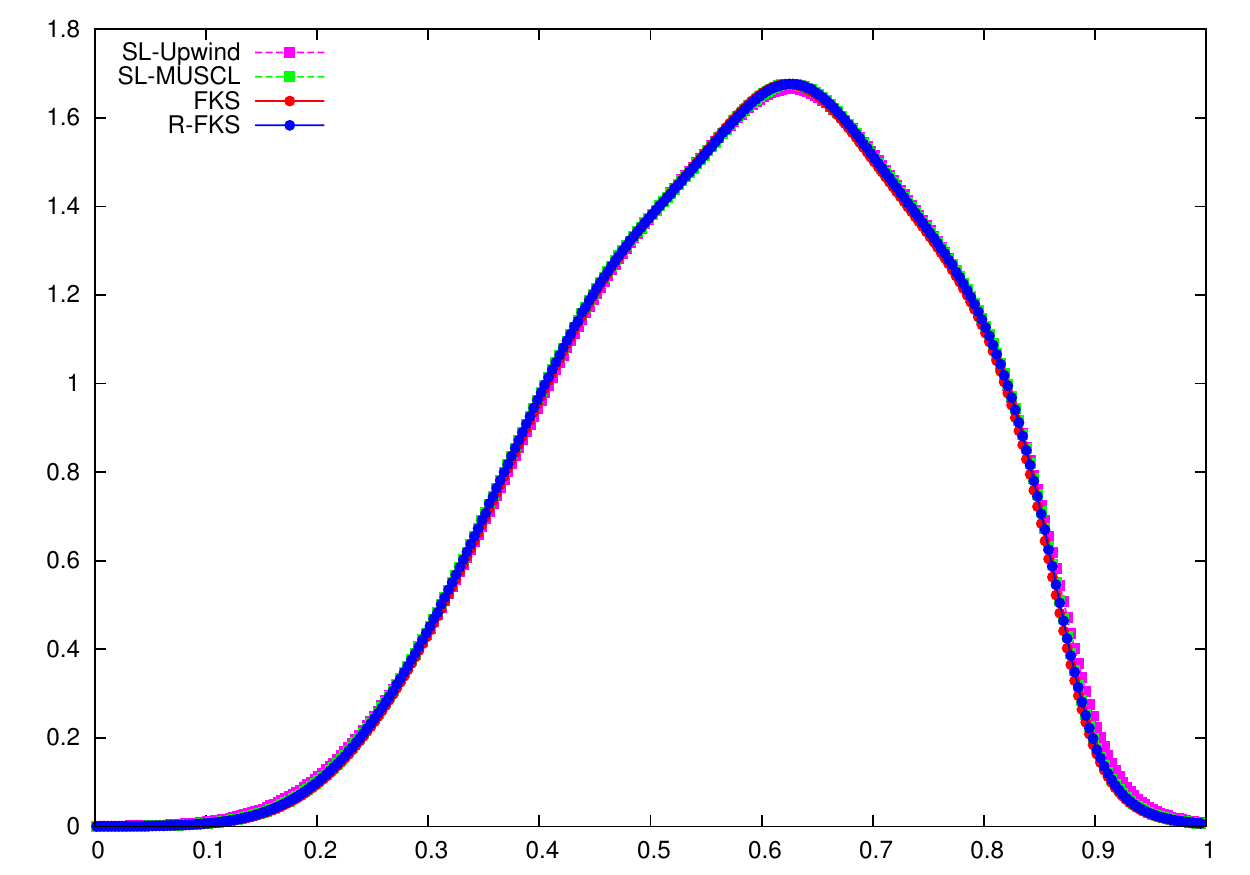}&
    \hspace{-0.35cm}
    \includegraphics[width=0.5\textwidth]{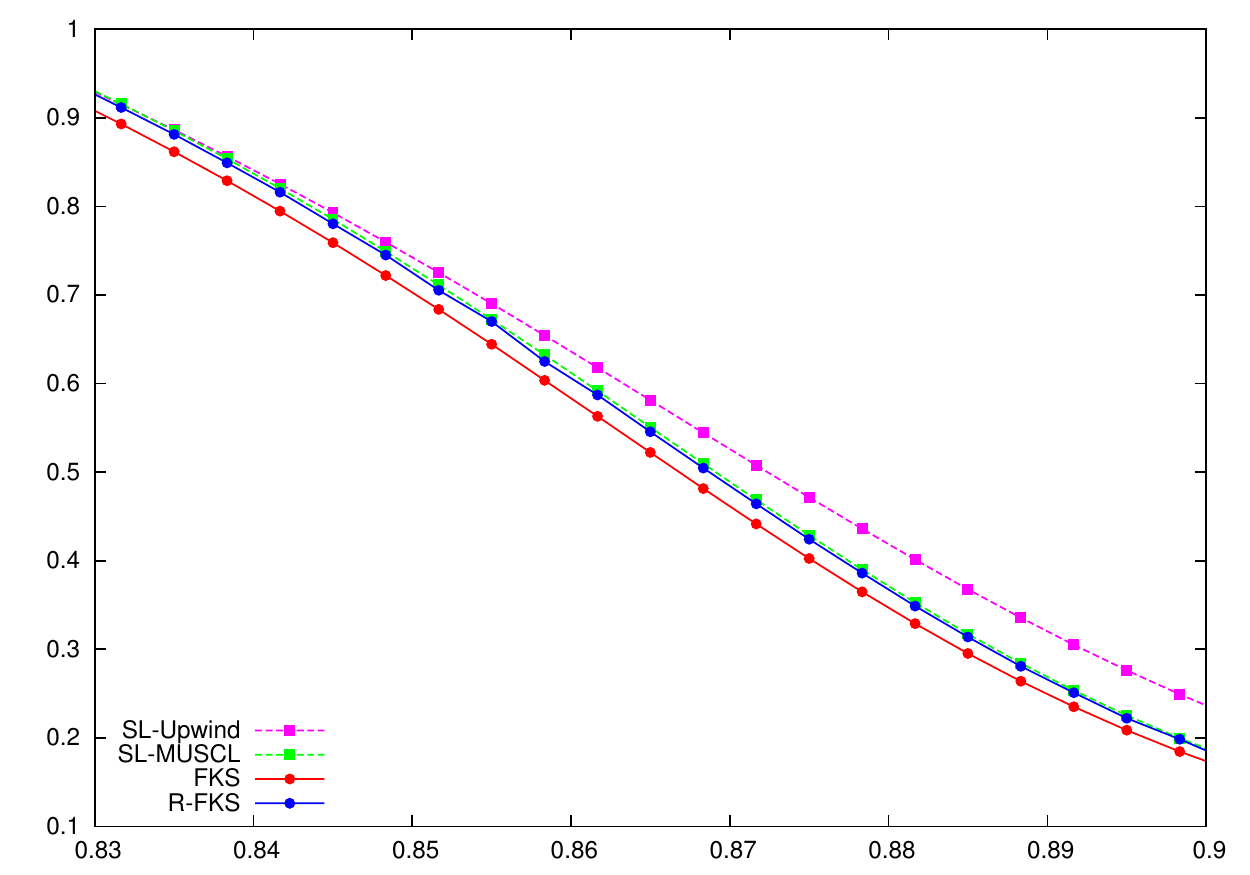} \\
  \end{tabular}
\caption{Riemann problem --- Solution at $t_{\text{final}}=0.07$ for the
velocity. Left full solution, right zoom close to the region in which the solution develops a shock wave in the limit of infinite
collisions. Collision frequency $\nu=10^{4}$ (top), $10^{3}$ (middle), $10^{2}$ (bottom).}
\label{fig_Sod2}
\end{center}
\end{figure}

\begin{figure}[ht!]
\begin{center}
  \begin{tabular}{cc}
    \hspace{-1.0cm}
    \includegraphics[width=0.5\textwidth]{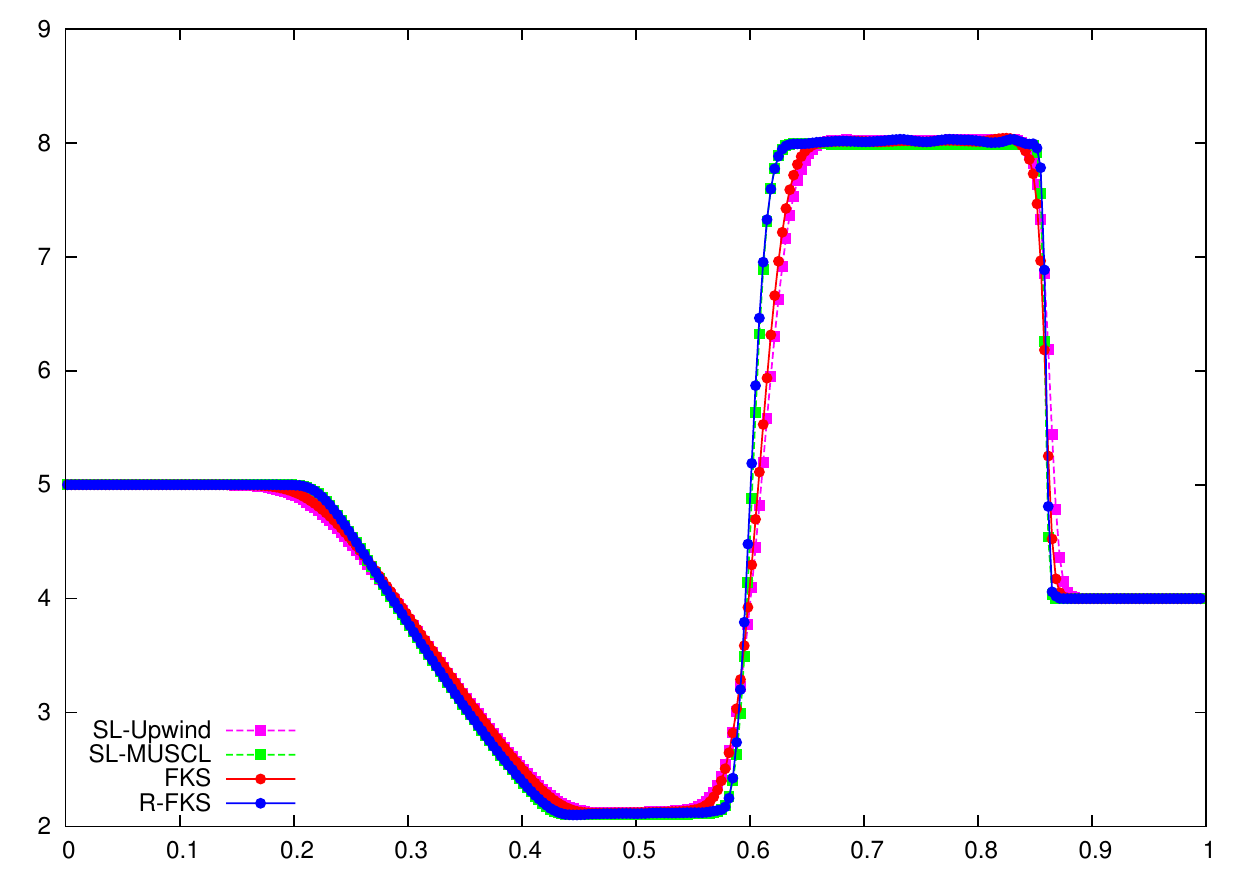}&
    \hspace{-0.35cm}
    \includegraphics[width=0.5\textwidth]{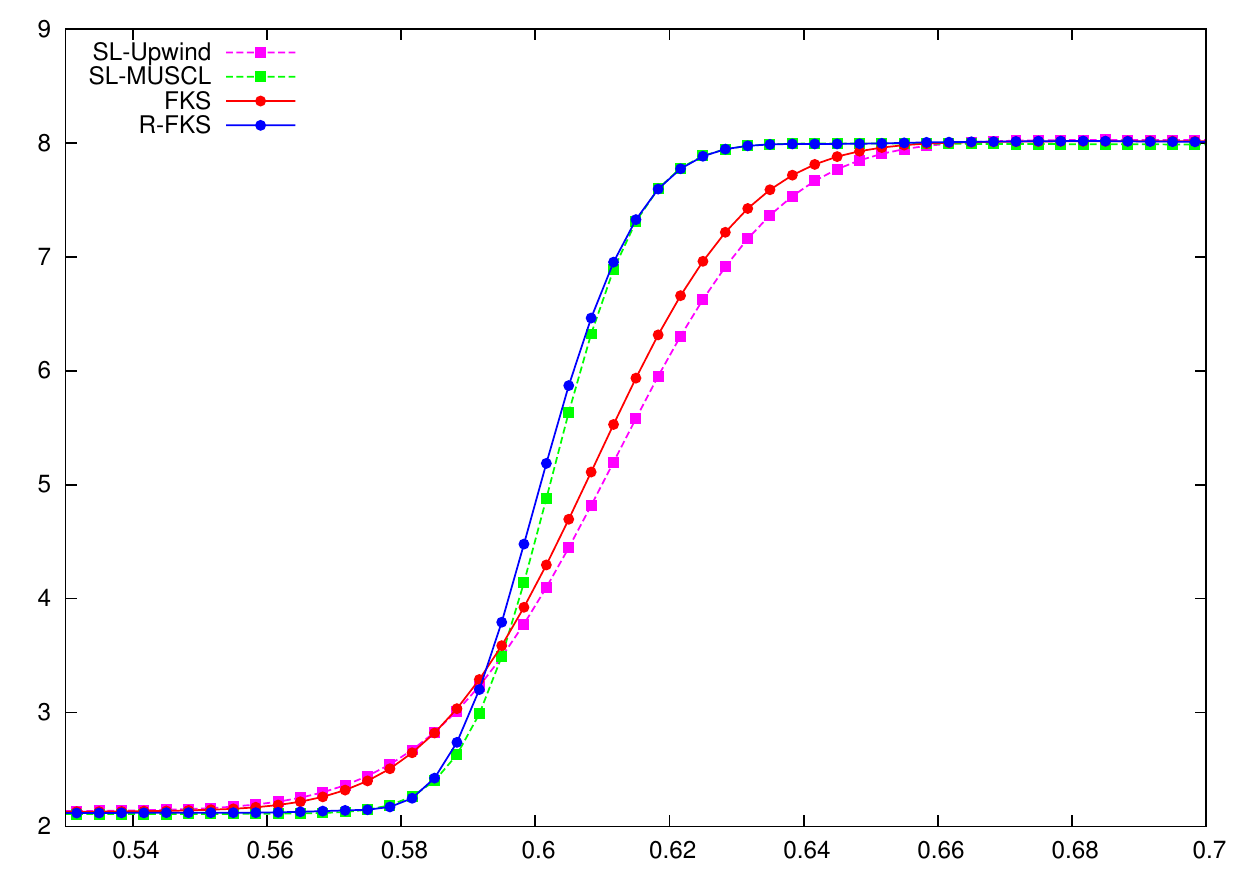} \\
    \hspace{-1.0cm}
    \includegraphics[width=0.5\textwidth]{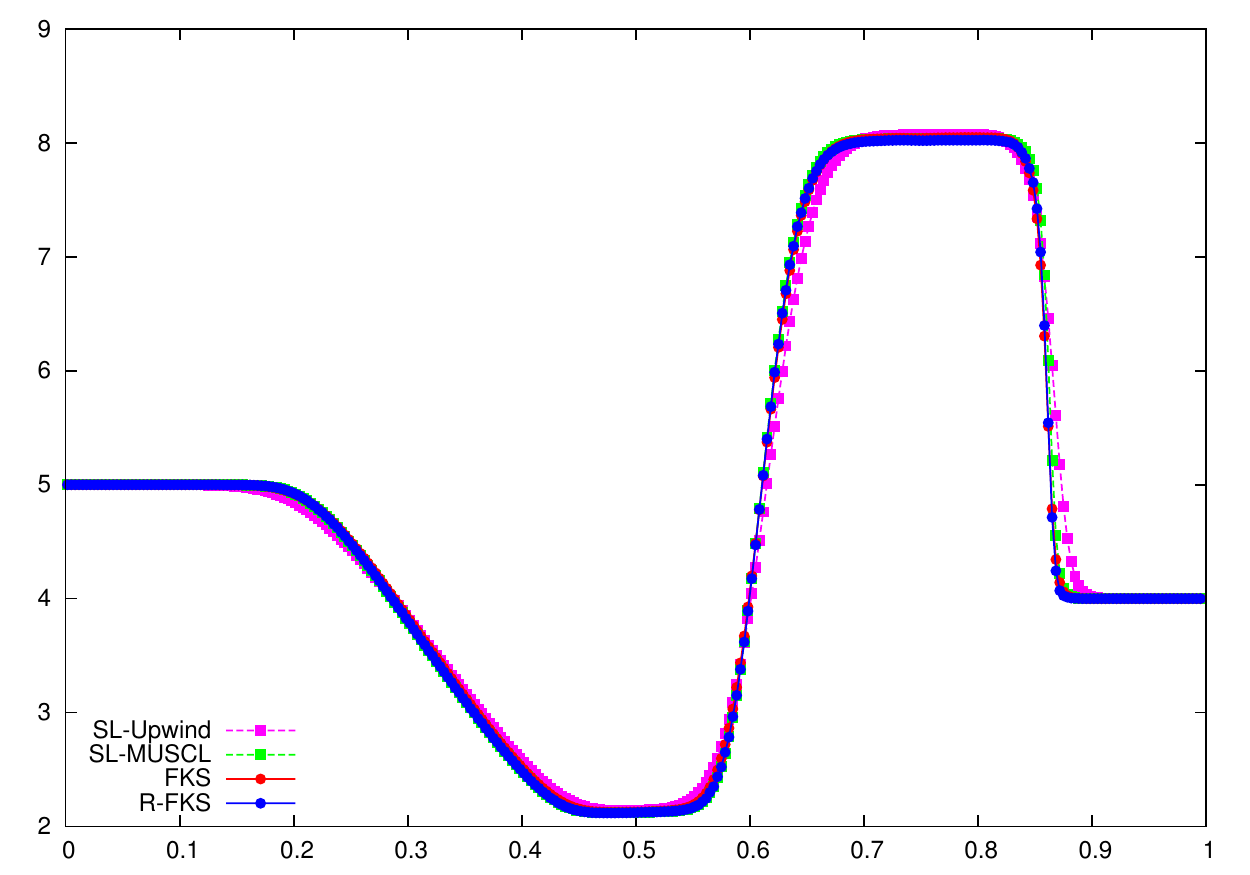}&
    \hspace{-0.35cm}
    \includegraphics[width=0.5\textwidth]{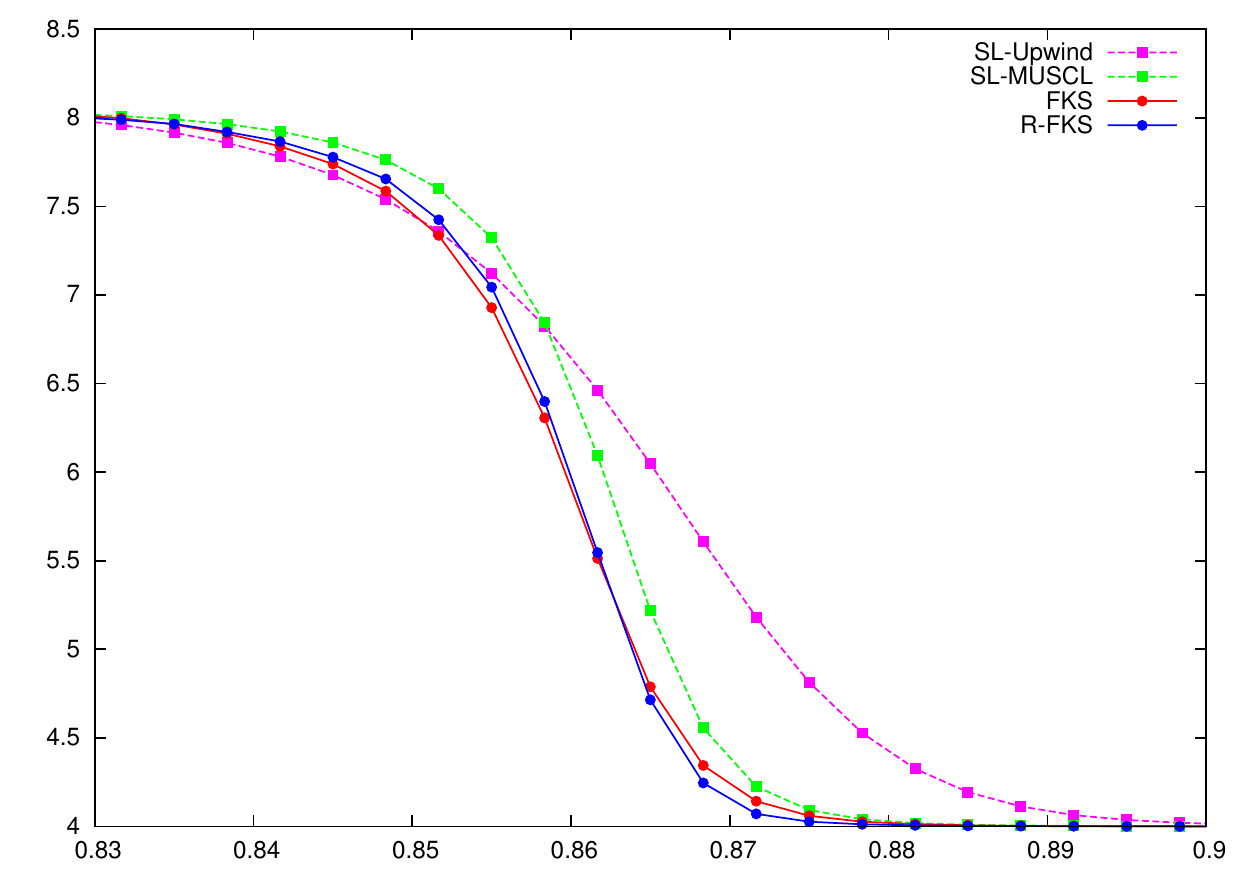} \\
    \hspace{-1.0cm}
    \includegraphics[width=0.5\textwidth]{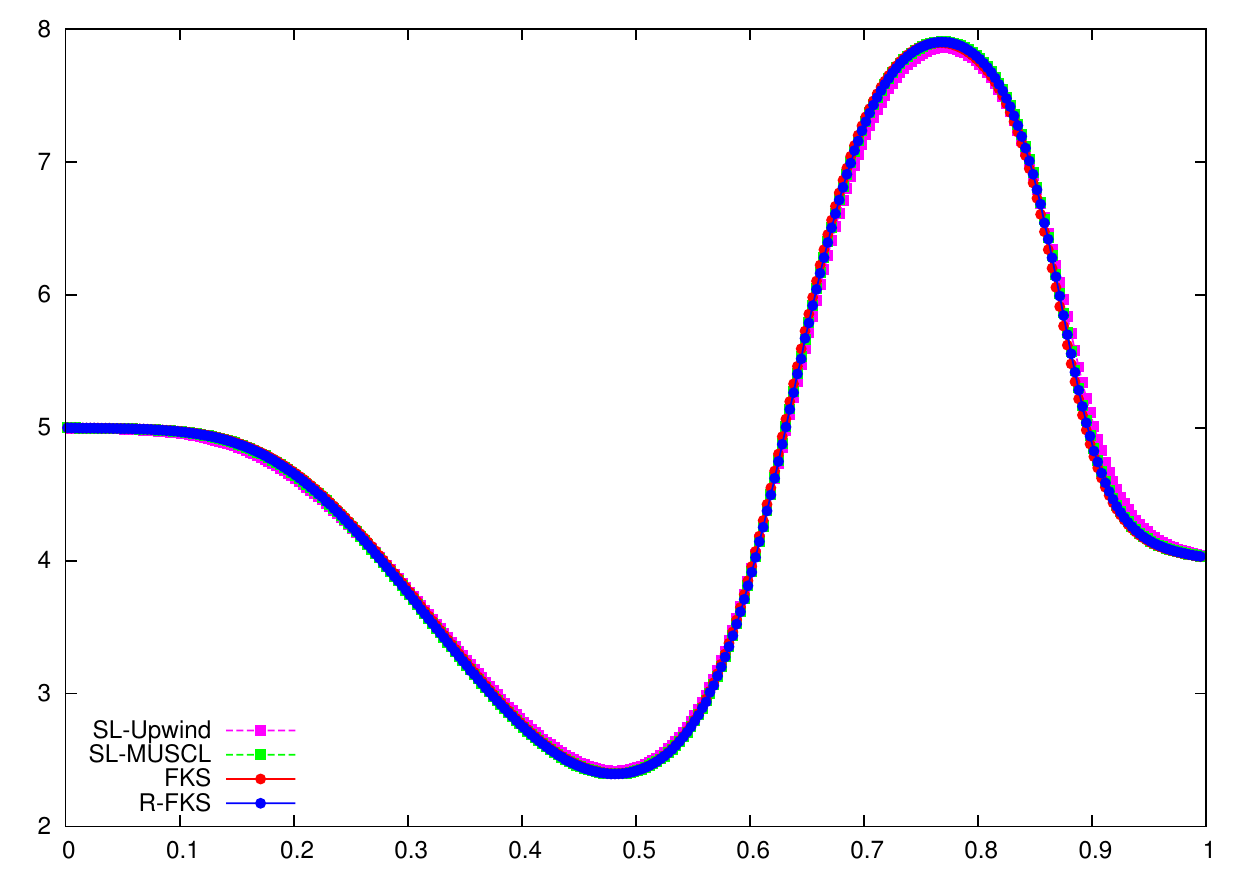}&
    \hspace{-0.35cm}
    \includegraphics[width=0.5\textwidth]{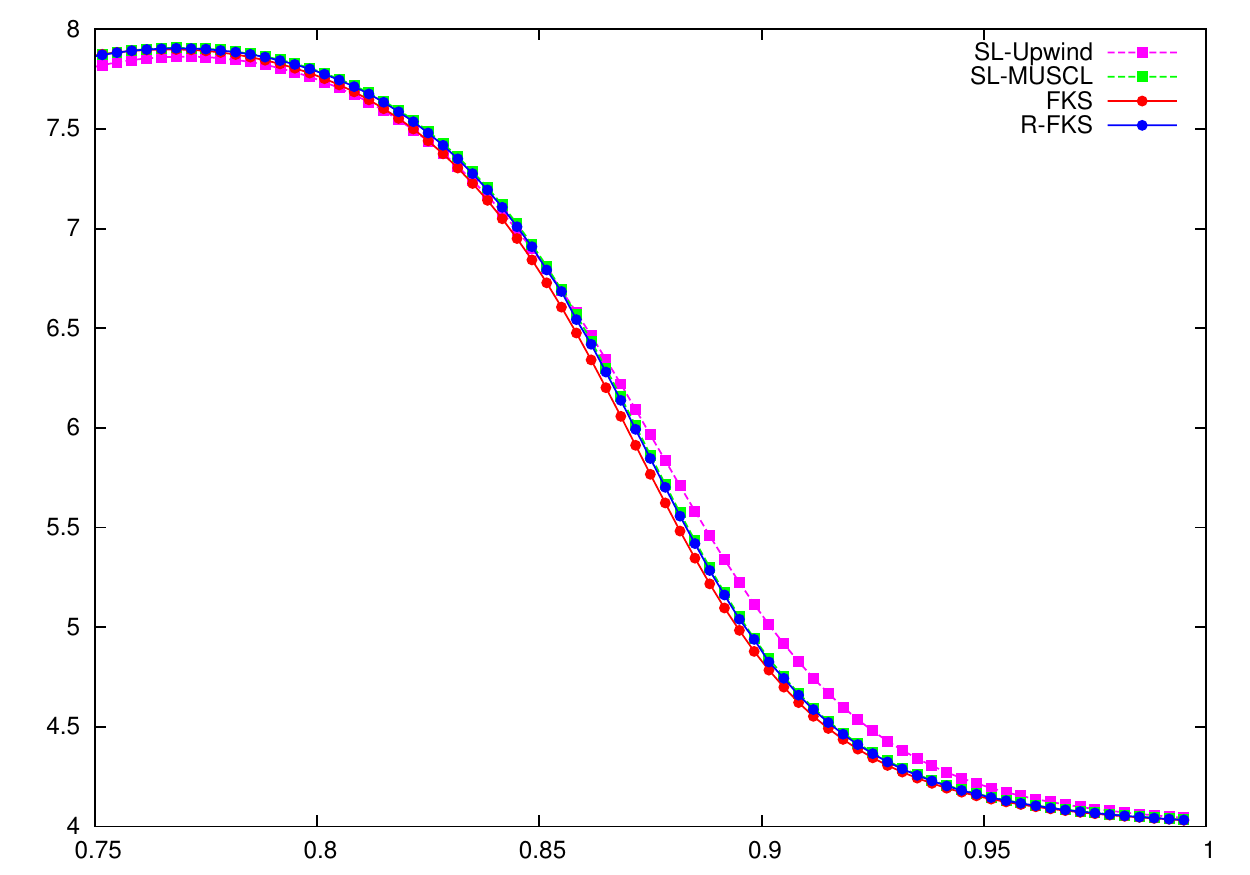} \\
  \end{tabular}
\caption{Riemann problem --- Solution at $t_{\text{final}}=0.07$ for the
temperature. Left full solution, right zoom close to the region in which the solution develops a contact discontinuity (top) or a shock wave (middle and bottom) in the limit of infinite
collisions. Collision frequency $\nu=10^{4}$ (top), $10^{3}$ (middle), $10^{2}$ (bottom).}
\label{fig_Sod3}
\end{center}
\end{figure}

% Description of the plots
In figure~\ref{fig_Sod1} is reported the density of the computed solution for a final time $t_{\text{final}}=0.07$, in figure~\ref{fig_Sod2} the velocity and in 
figure~\ref{fig_Sod3} the temperature for the three different collision frequencies $\nu=10^{4}$ (top panels), $10^{3}$ (middle panels), $10^{2}$ (bottom panels).
On the left side the solution on the entire domain is plotted while on the right side a magnification of the solutions in some interesting locations is reported. 
In particular for the density the zoom is around the region in which the solution develops a contact discontinuity in the fluid dynamic limit, 
for the velocity around the region in which the solution forms a shock wave in the limit, while for the temperature the zooms are 
around the region of the contact discontinuity for $\nu=10^{4}$  and around the region of the shock wave for $\nu=10^{3}$ 
and $\nu=10^{2}$. 
The results obtained by the four schemes are reported on each plot: FKS in red, R-FKS in blue, SL-Upwind in magenta and SL-MUSCL in green. \\
The results show that SL-MUSCL and R-FKS are clearly more accurate than SL-Upwind and FKS
in the highly collisional case $\nu=10^{4}$. In particular, for the density we observe that the contact discontinuity is better captured by the R-FKS scheme compared to the SL-MUSCL, 
while for the temperature and the velocity the solutions are almost superimposed. The same behavior is observed for the shock wave in the fluid dynamic limit, 
the solution is better captured (less points on the shock wave are present) by the R-FKS scheme in comparison with the three other schemes employed, while the complexity 
of the scheme is the one of a first order scheme (only first order polynomial are employed). 
On the other hand, when $\nu$ decreases, the differences between the schemes are less pronounced as the physical diffusion takes over on the waves. However, for $\nu=10^3$ we still observe
a difference between the four solutions, in particular the two semi-Lagrangian schemes SL-Upwind and SL-MUSCL are more diffusive than the two FKS schemes (a convergence test, not reported, 
shows that the position of the waves converges to the one of the FKS method when the number of points increases). In addition, comparing FKS and R-FKS we observe that R-FKS produces slightly 
less diffusive results. Finally, in the case in which the kinetic effects are large, i.e. $\nu=10^{2}$ the four solutions are almost superimposed and it is 
difficult to appreciate the differences between the methods, only the SL-Upwind is clearly slightly more diffusive than the others. \\
In order to better appreciate the behavior of our R-FKS method in the next paragraph we consider a case test in which in the kinetic regime many oscillations are present, this problem
permits to enhance specific characteristics of the schemes.

%=== B E G I N    F I G U R E ================================================
%=== E N D   F I G U R E ================================================

%=== B E G I N    F I G U R E ================================================
%=== E N D   F I G U R E ================================================

%=== B E G I N    F I G U R E ================================================

\subsection{Highly oscillating kinetic problem} \label{ssec:oscill}
\begin{figure}
  \begin{center}
    \includegraphics[width=0.49\textwidth]{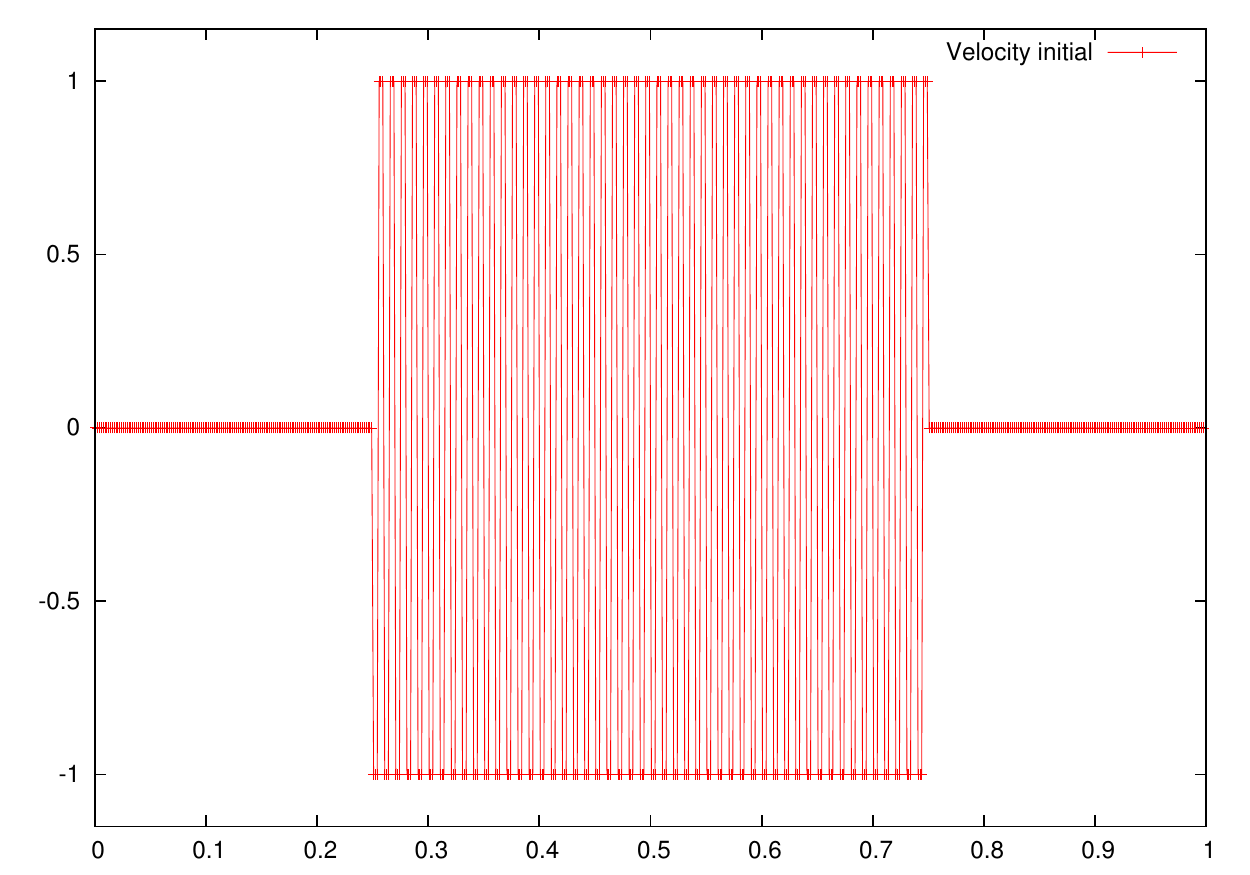}    
    \includegraphics[width=0.49\textwidth]{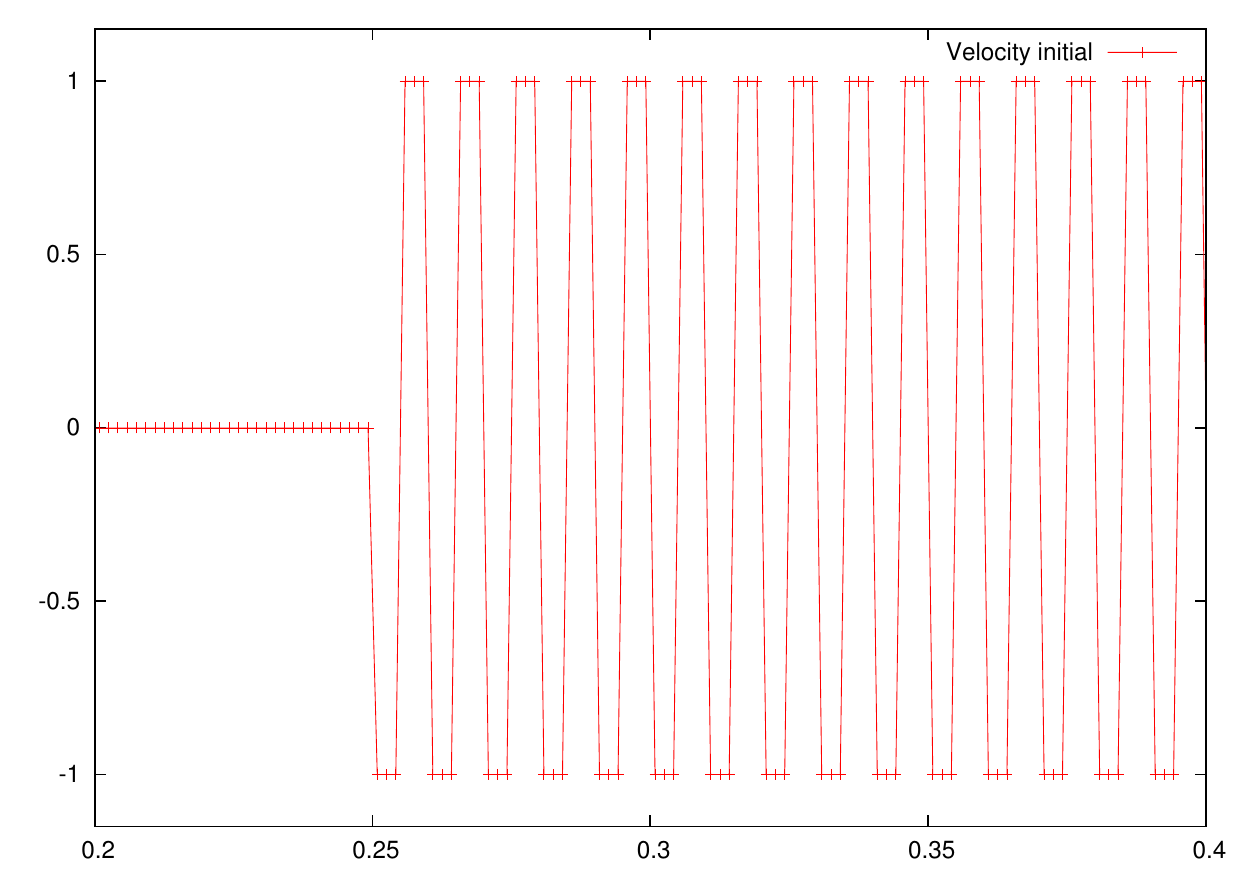}    
    \caption{ \label{fig:oscill} Initial velocity field for the highly oscillating 
      kinetic test case --- $\delta=0.02$ and $600$ mesh cells are used ---
      Left: full view on $\Omega$, right: zoom on $[0.2,0.4]$.}
  \end{center}
\end{figure}

The previous test case was clear in illustrating the gain of the R-FKS scheme compared to the other schemes close to the fluid limit and for intermediate regimes, i.e. $\nu=10^3$, while
the gain was less obvious when the regime was far from the thermodynamics equilibrium.
Contrarily, the test case in this section is designed to illustrate the anti-diffusive behavior of 
the FKS methods compared to classical semi-Lagrangian
schemes in a kinetic regime. In order to do that, let us consider an homogeneous gas, $\rho=1$, at constant temperature, $T=5$ and at rest $u=0$ on a domain $\Omega=[0,1]$. 
This gas is further animated by a stair case oscillating velocity field ($u=\pm 1$) over a period $0<\delta\ll L$ as
\bea
u(x) = \left\{ \begin{array}{lll}
+1 & \text{if} &  x \in [x_{2k}-\delta/2:x_{2k}+\delta/2] \cap [0.25;0.75] ,\\
-1 & \text{if} &  x \in [x_{2k+1}-\delta/2:x_{2k+1}+\delta/2] \cap  [0.25;0.75],\\
0 & \text{else} &  
\end{array}
\right. \qquad \forall 0 \leq k \leq M/2 ,
\eea
where $\{ x_i \}_{i=1,M}$ is the spacial mesh, see figure~\ref{fig:oscill} for an illustration 
of the initial velocity field.
%=== B E G I N    F I G U R E ================================================
%=== E N D   F I G U R E ================================================

We then consider a kinetic regime, $\nu=10^{2}$, for which kinetic effects are non negligible.
In this situation, the highly oscillating velocity field generates waves emanating from
 each discontinuity. These waves
 later interact, creating secondary waves further in interaction and so on. 
The physical relaxation will damp some of these waves, their amplitude and their structure,
and, any low dissipative numerical method must reproduce such behaviors. 
Associated with this physical dissipation, due to the embedded numerical diffusion of the schemes,
 some small structures may 
also be lost during the time evolution, or, at least, they may be damped.
% Goal
The goal of this test is to show that the class of FKS methods are capable to maintain in time these small structures whereas they are, for a large part, lost using
a semi-Lagrangian approach due to the intrinsic numerical diffusion. 

% "Raph : This looks already like the conclusion of the test"
%Moreover, in order for the semi-Lagrangian schemes to obtain similar behaviors of the FKS methods, we
%need a very large number of points in space and consequently we need an investment of computer %resources which is much larger compared to that of the FKS methods.%
%This is of particular interest in the multidimensional settings for which precision is mandatory in order to avoid very fine meshes.
In the following we consider two kind of tests which measure the performances of the schemes. The first one is a mesh convergence test while the second one is an efficiency test.

\subsubsection{Mesh convergence study}
\begin{figure}[ht!]
\begin{center}
  \begin{tabular}{ccc}
    \hspace{-1.2cm}
    \includegraphics[scale=0.44]{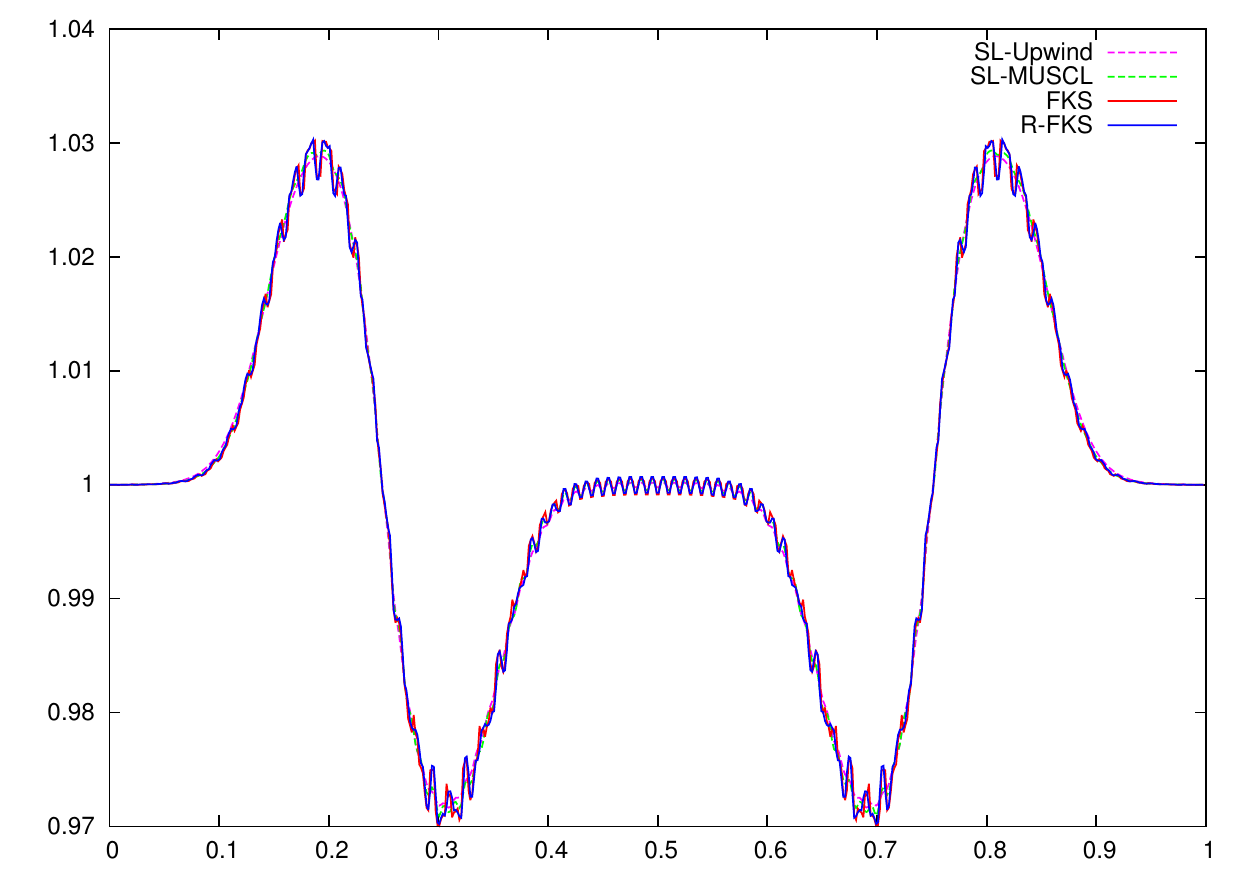}&
    \hspace{-0.75cm}
    \includegraphics[scale=0.44]{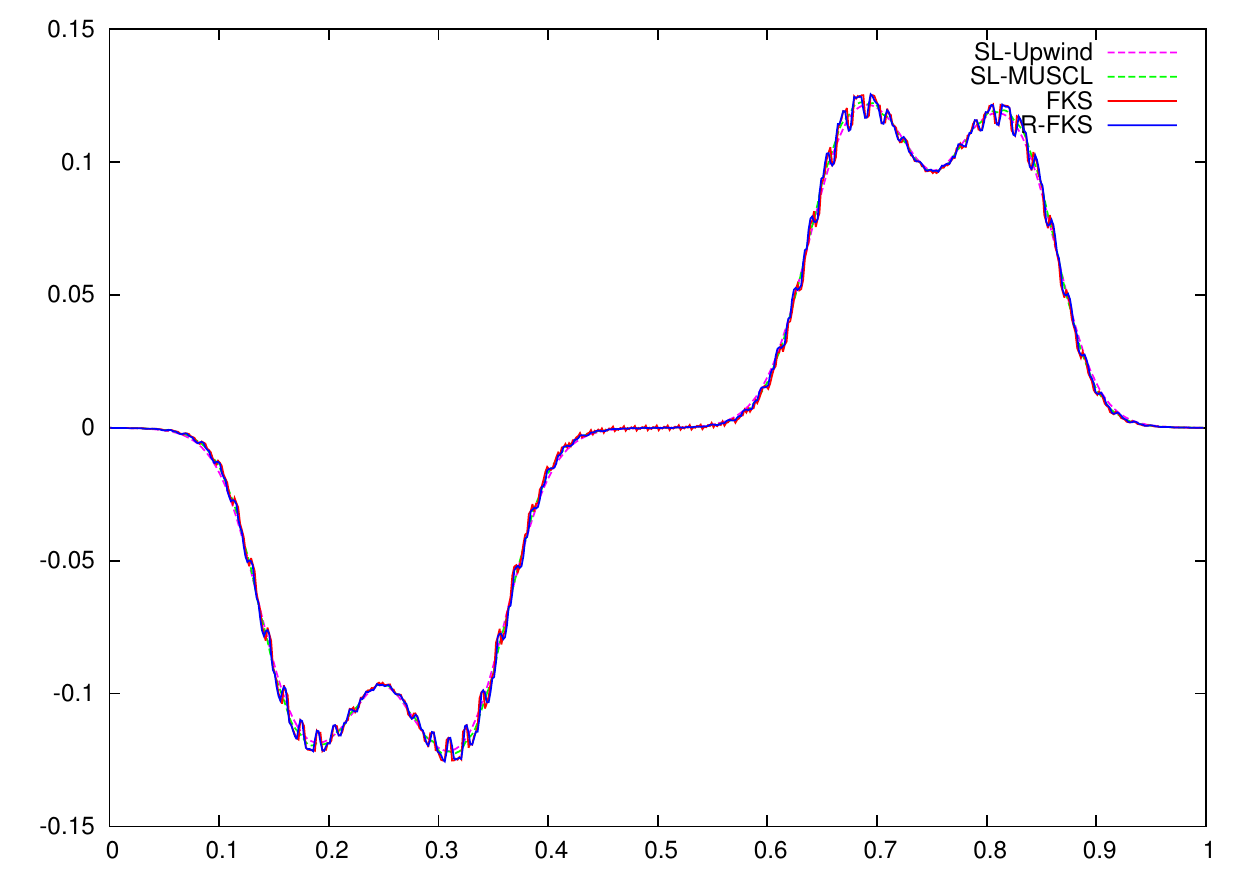}&
    \hspace{-0.75cm}
    \includegraphics[scale=0.44]{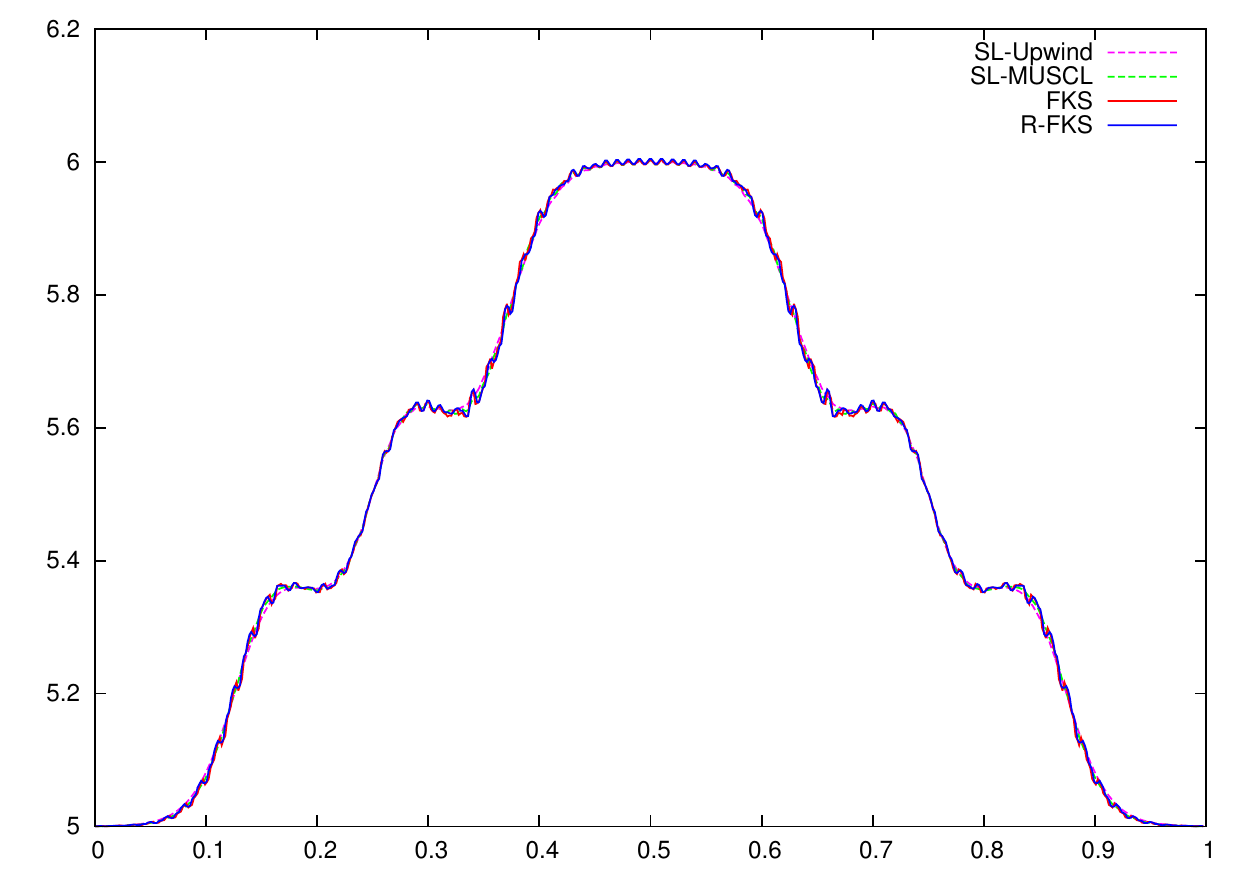} \\
  \end{tabular}
\caption{Oscillating problem --- $600$ cells --- 
  Solution at $t_{\text{final}}=0.025$ for the
density (left), mean velocity (middle) and temperature (right) with $\nu=10^{2}$
and the SL-Upwind, SL-MUSCL, FKS and R-FKS schemes.
}
\label{fig:oscill600}
\end{center}
\end{figure}

\begin{figure}[ht!]
\begin{center}
  \begin{tabular}{ccc}
    \hspace{-1.2cm}
    \includegraphics[scale=0.6]{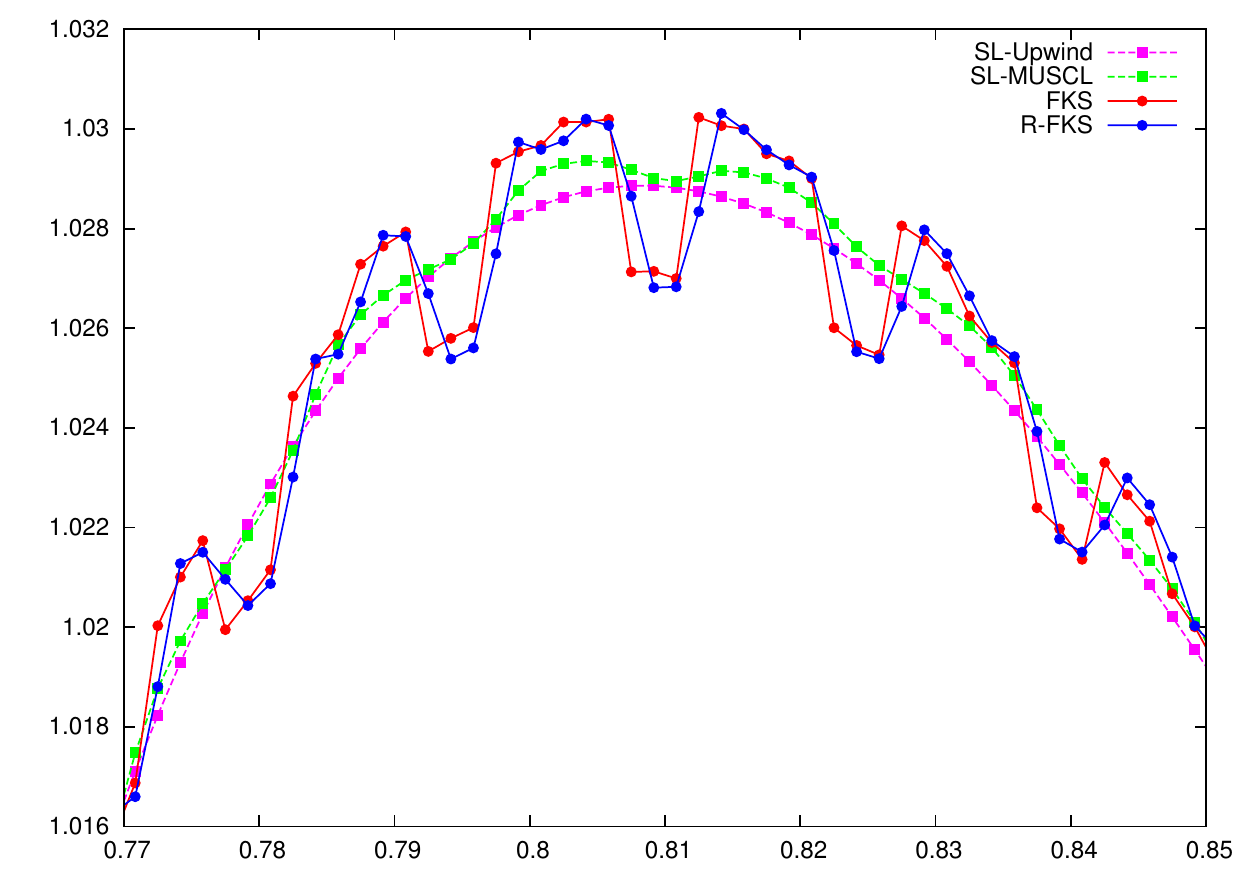}&
    \hspace{-0.75cm}
    \includegraphics[scale=0.6]{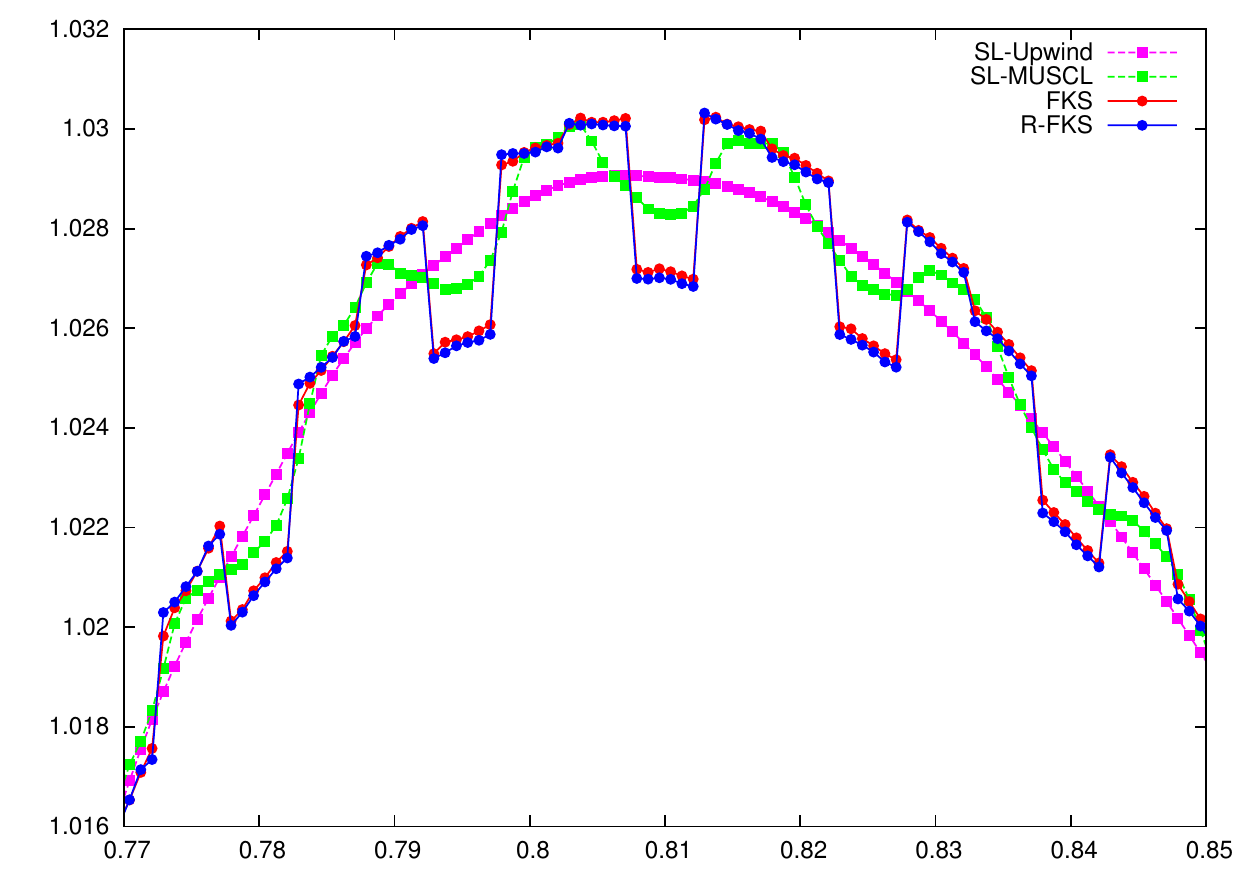}\\
    \hspace{-1.2cm}
    \includegraphics[scale=0.6]{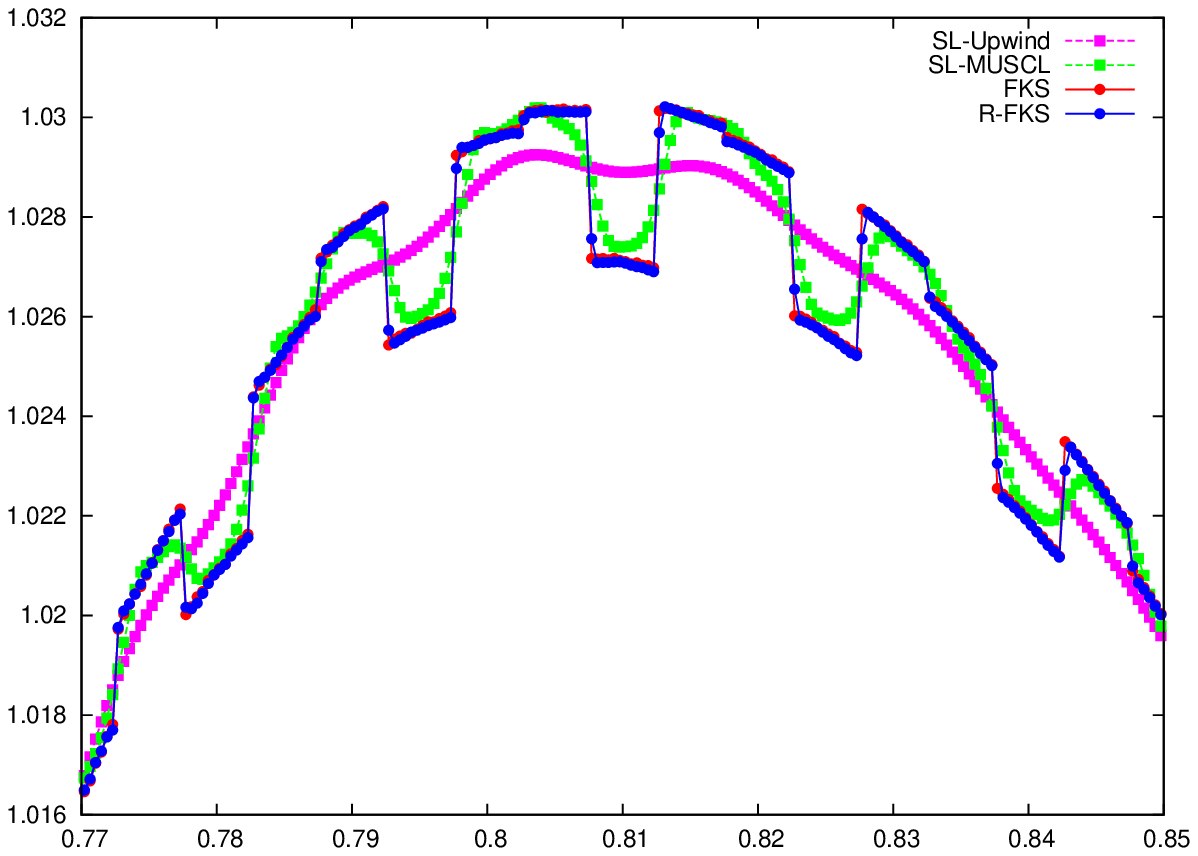}&
    \hspace{-0.75cm}
    \includegraphics[scale=0.6]{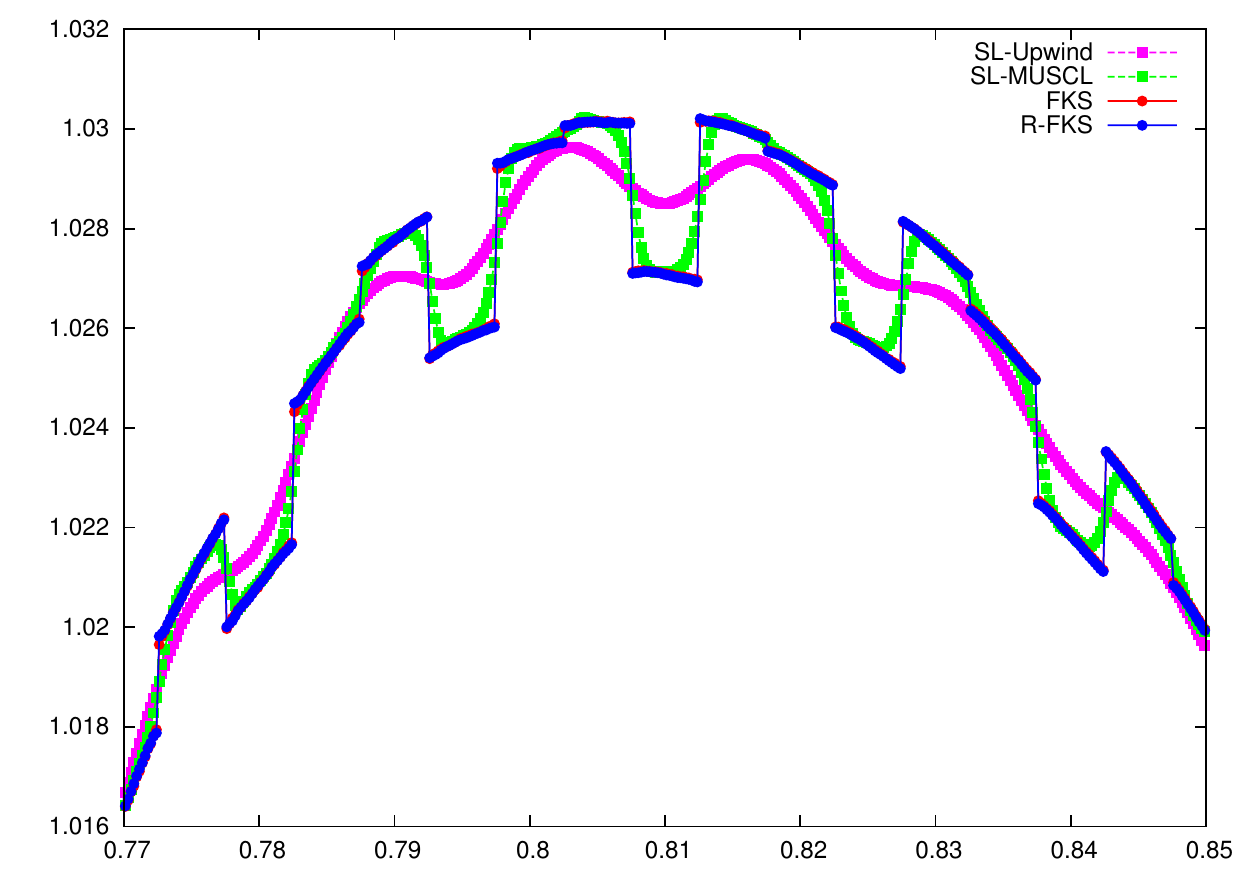}
  \end{tabular}
\caption{Oscillating problem --- Zooms on
  density for 
  $600$, $1200$, $2400$ and $4800$ cells (from top-left to bottom-right) --- 
  Solution at $t_{\text{final}}=0.025$  with $\nu=10^{2}$
  for the SL-Upwind, SL-MUSCL, FKS and R-FKS schemes.
}
\label{fig:oscillcomp_density}
\end{center}
\end{figure}

\begin{figure}[ht!]
\begin{center}
  \begin{tabular}{ccc}
    \hspace{-1.2cm}
    \includegraphics[scale=0.6]{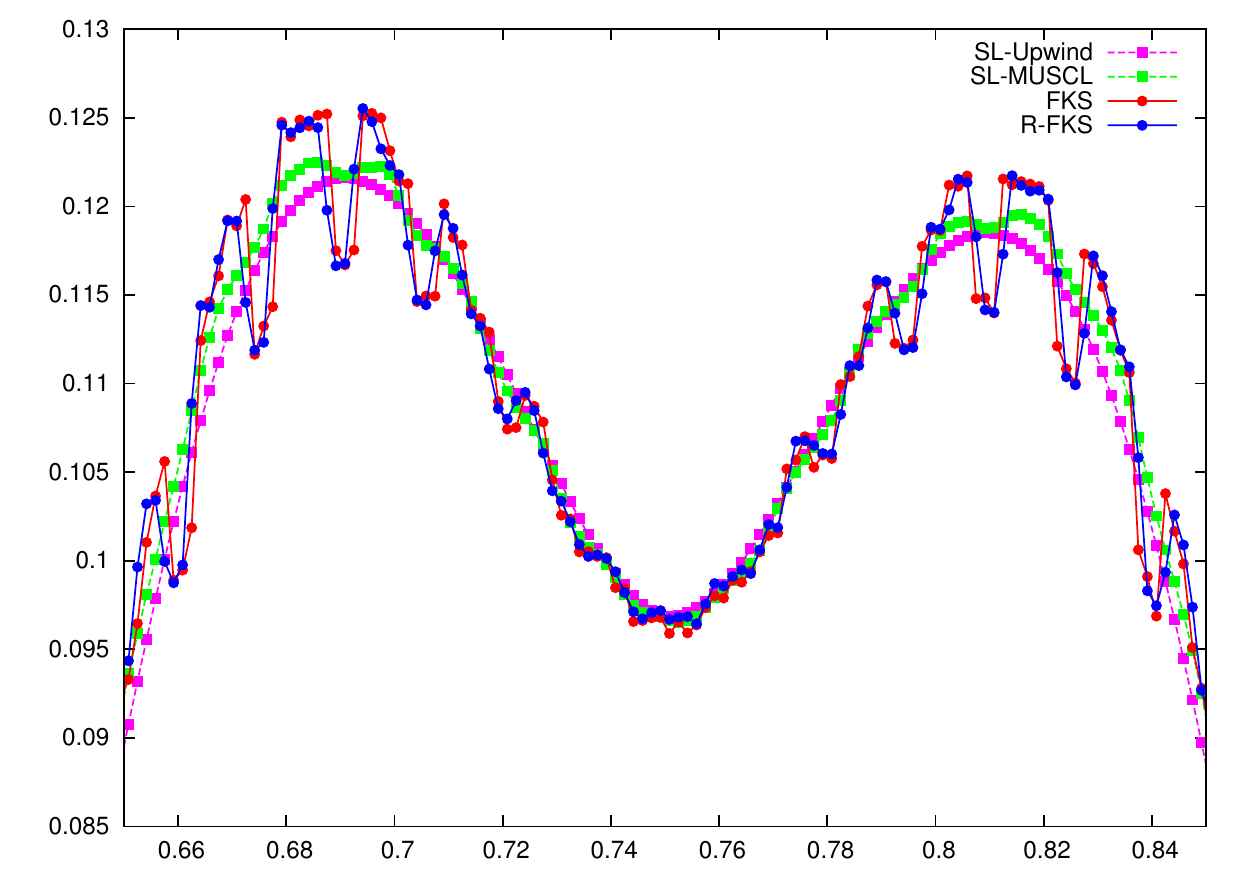}&
    \hspace{-0.75cm}
    \includegraphics[scale=0.6]{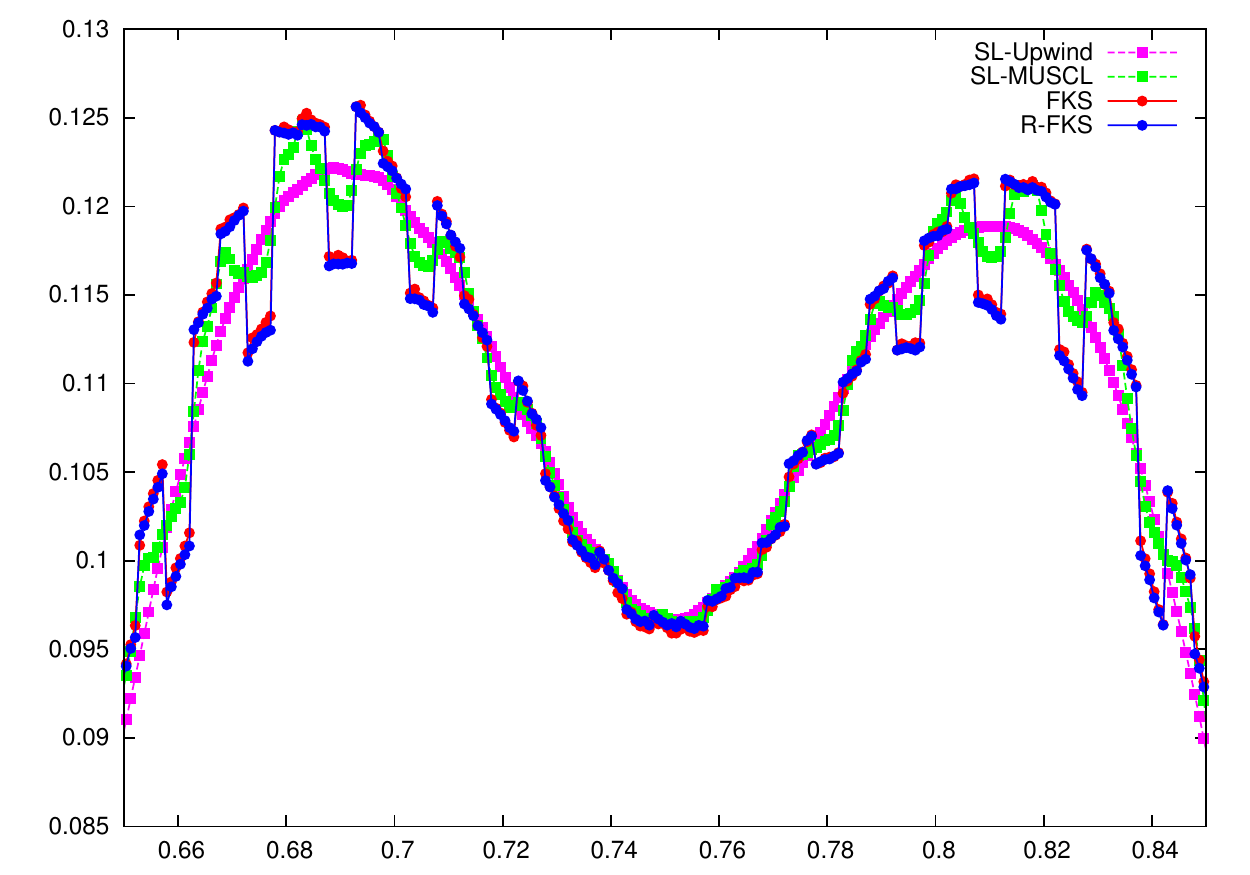}\\
    \hspace{-1.2cm}
    \includegraphics[scale=0.6]{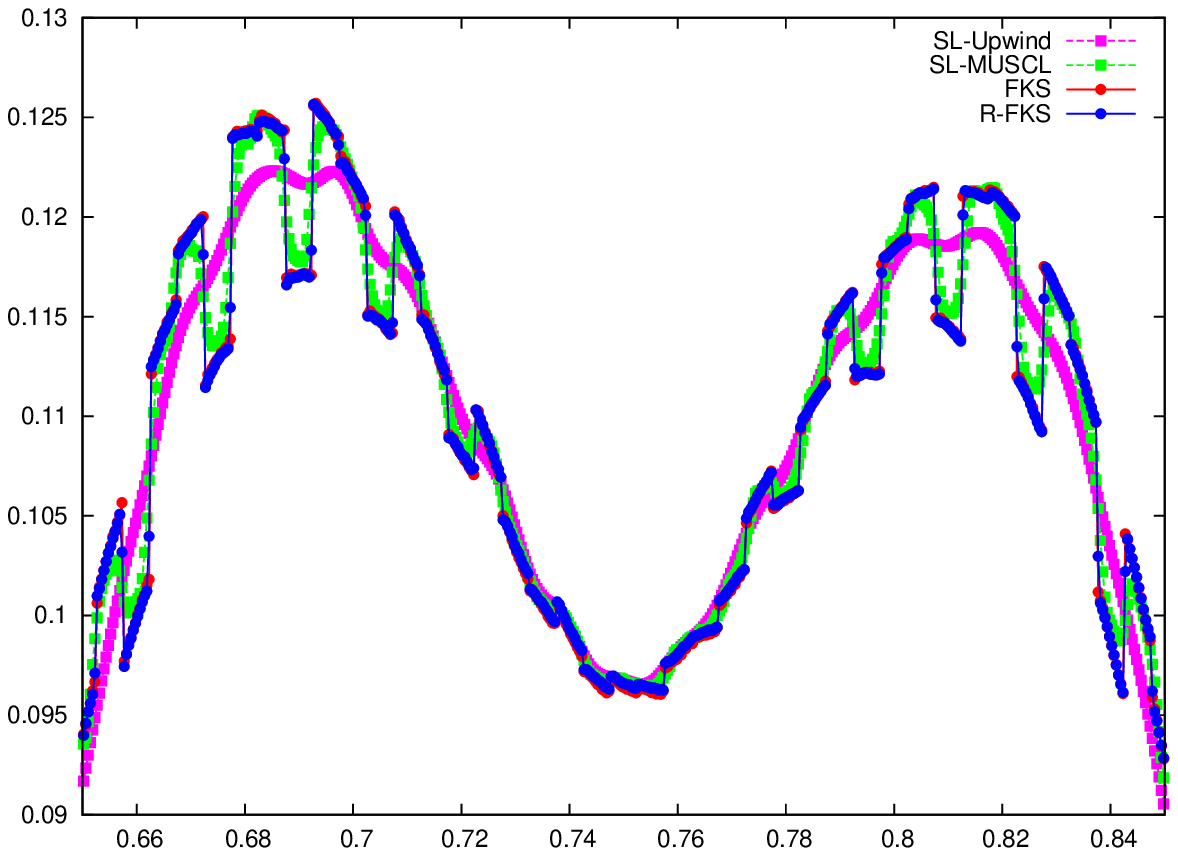}&
    \hspace{-0.75cm}
    \includegraphics[scale=0.6]{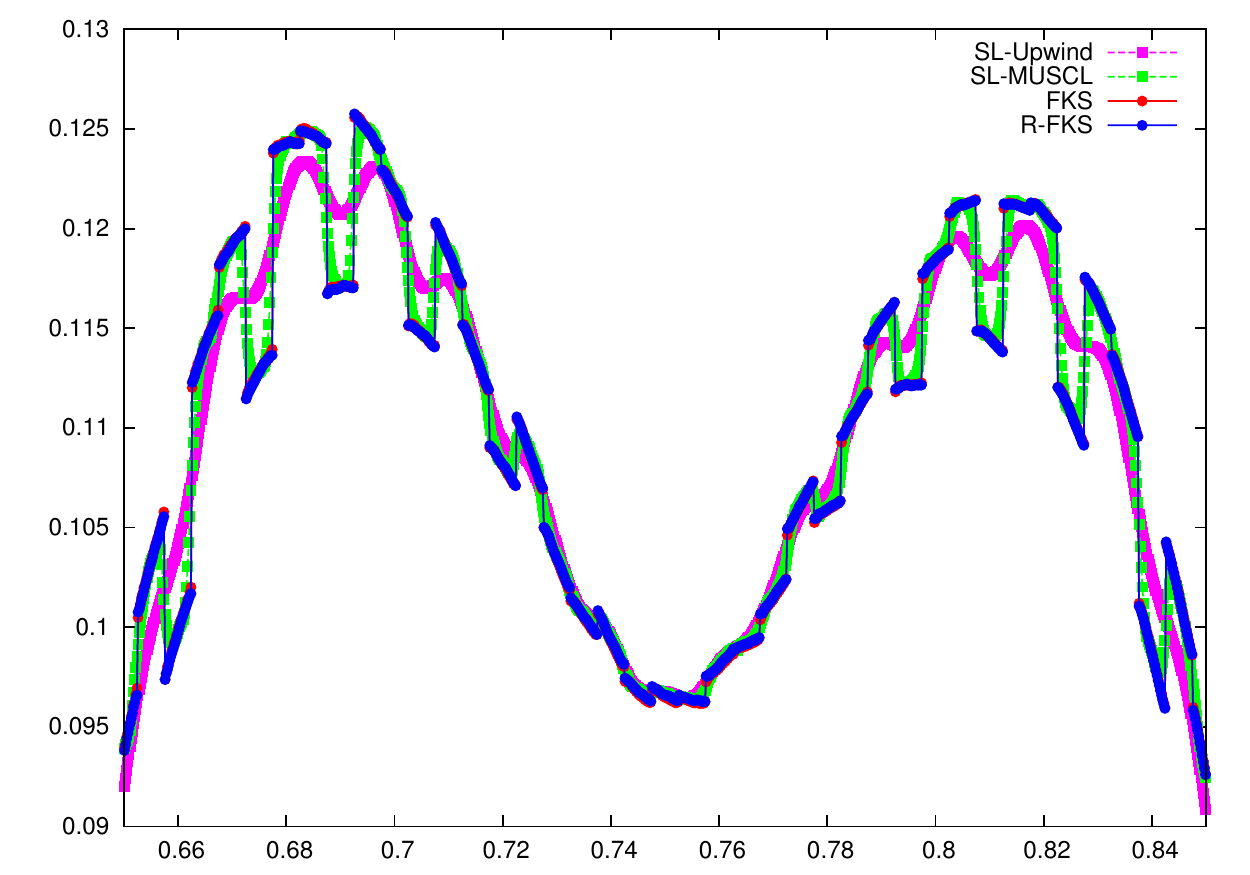}
  \end{tabular}
\caption{Oscillating problem --- Zooms on
  velocity for 
  $600$, $1200$, $2400$ and $4800$ cells (from top-left to bottom-right) --- 
  Solution at $t_{\text{final}}=0.025$  with $\nu=10^{2}$
  for the SL-Upwind, SL-MUSCL, FKS and R-FKS schemes.
}
\label{fig:oscillcomp_velocity}
\end{center}
\end{figure}

\begin{figure}[ht!]
\begin{center}
  \begin{tabular}{ccc}
    \hspace{-1.2cm}
    \includegraphics[scale=0.6]{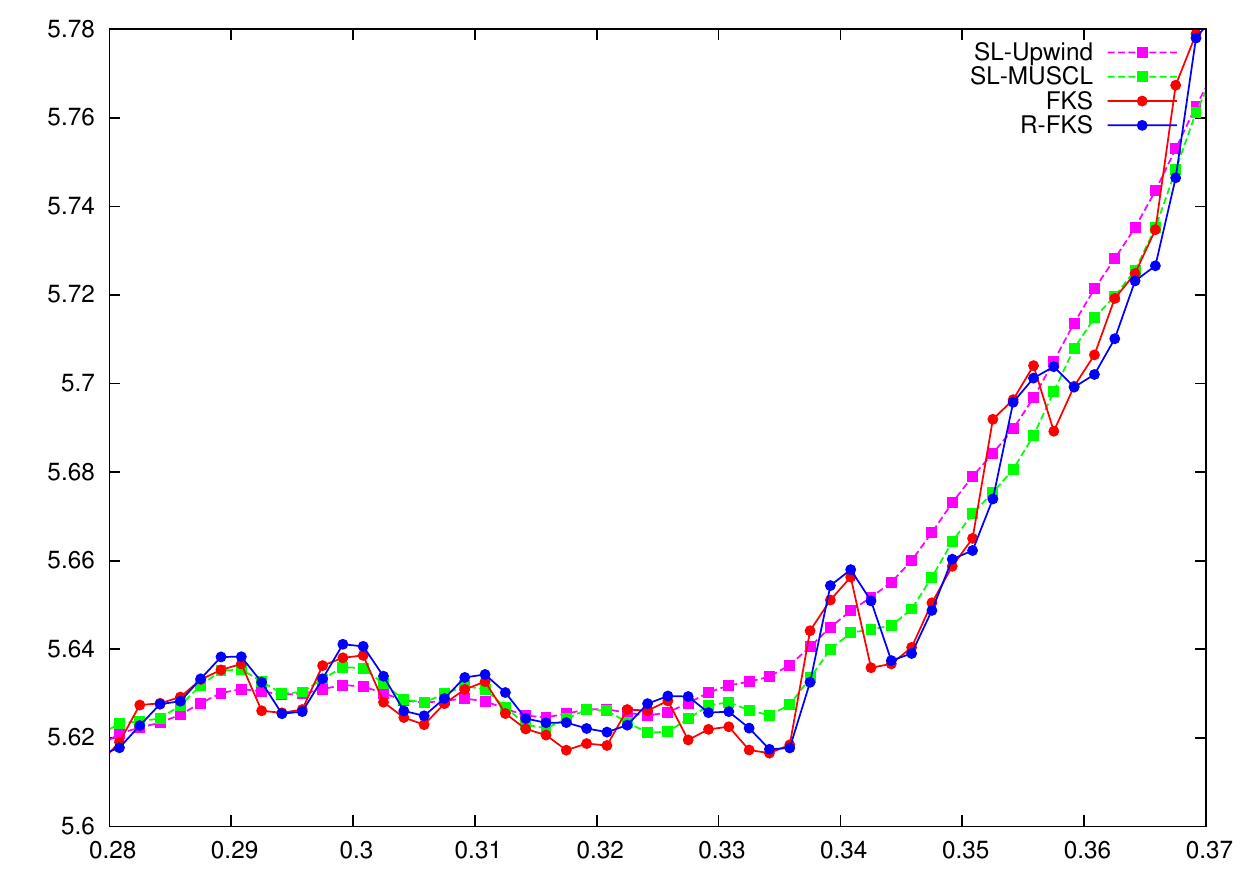}&
    \hspace{-0.75cm}
    \includegraphics[scale=0.6]{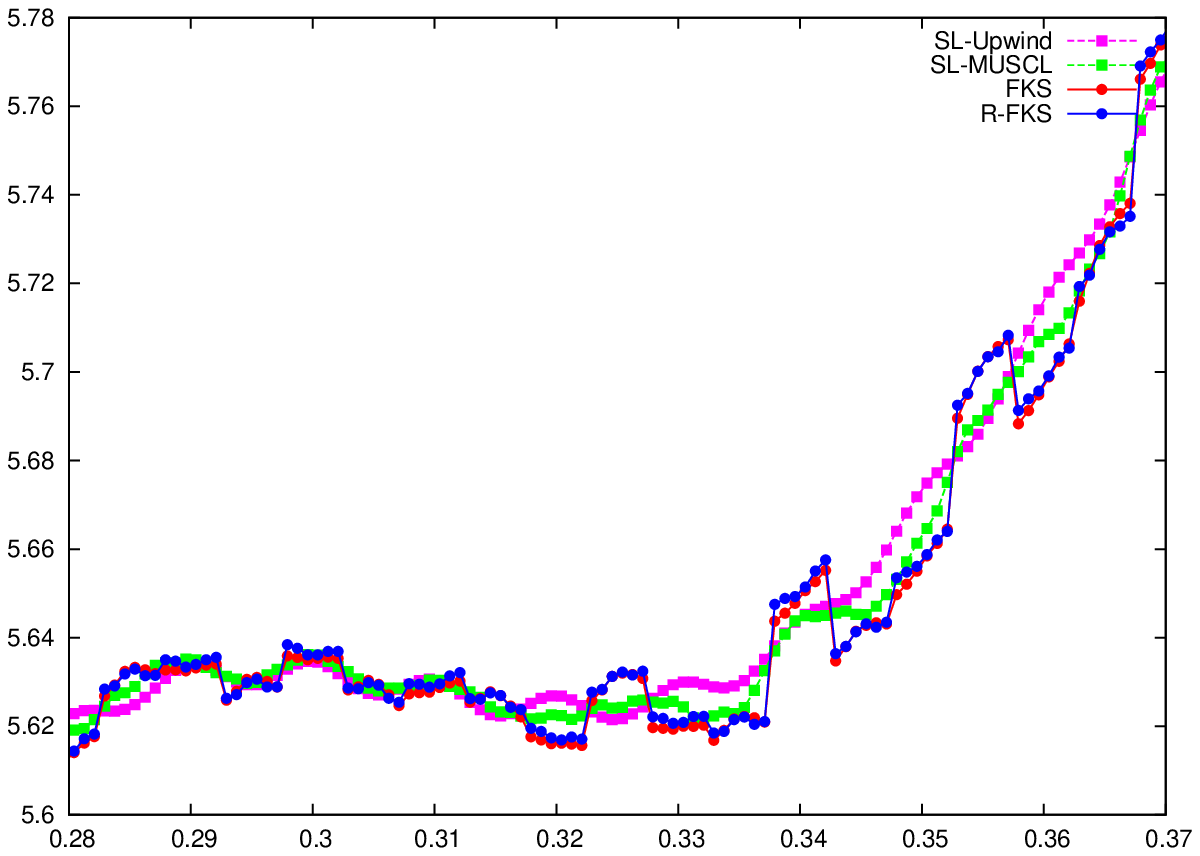}\\
    \hspace{-1.2cm}
    \includegraphics[scale=0.6]{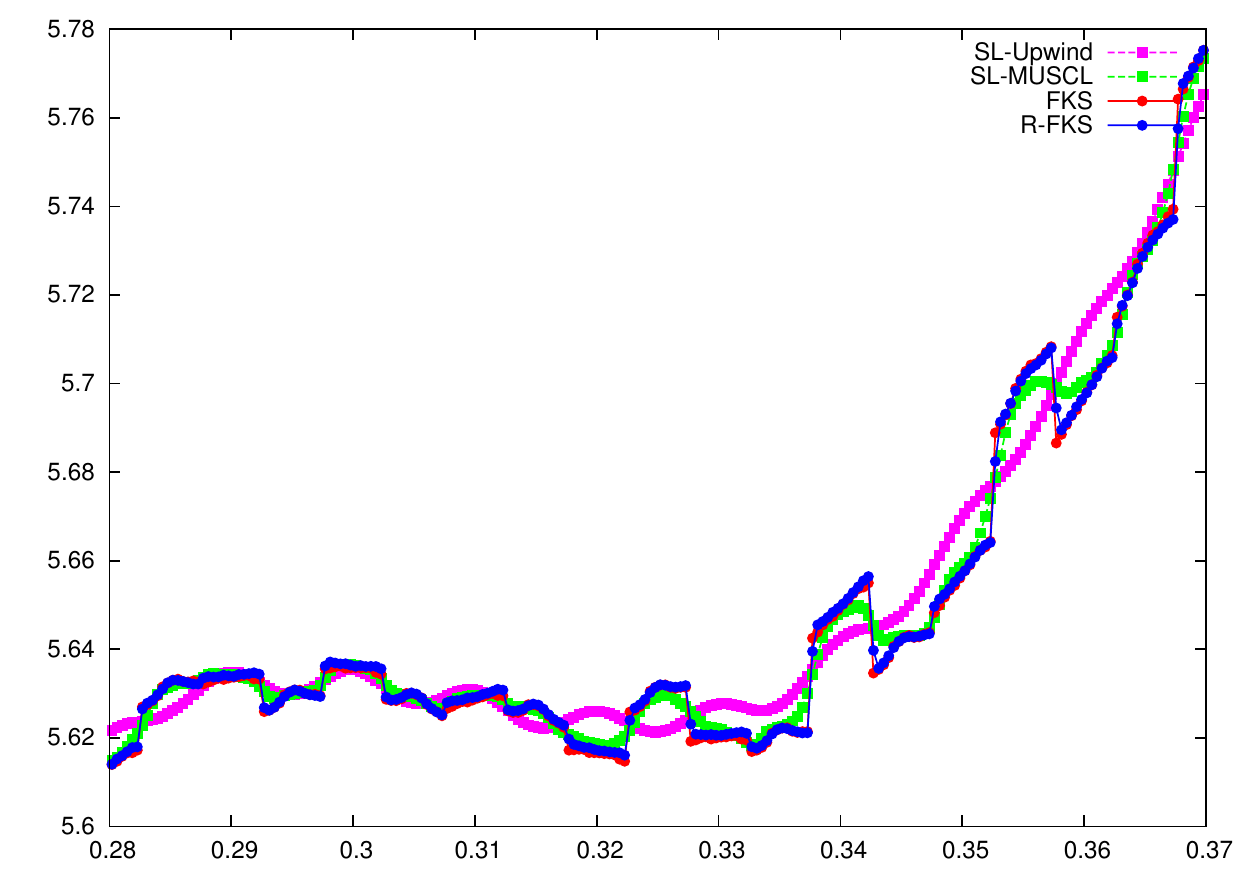}&
    \hspace{-0.75cm}
    \includegraphics[scale=0.6]{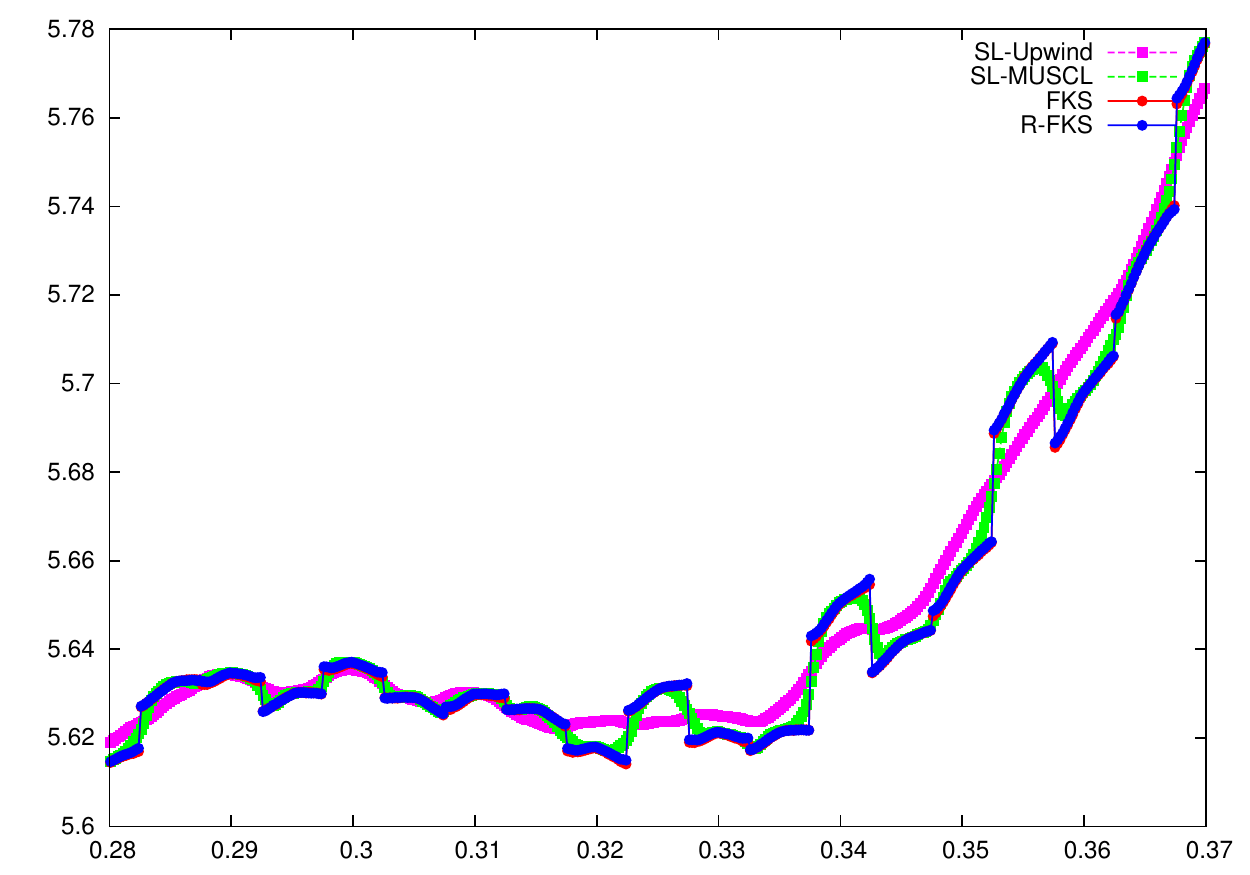}
  \end{tabular}
\caption{Oscillating problem --- Zooms on
  temperature for 
  $600$, $1200$, $2400$ and $4800$ cells (from top-left to bottom-right) --- 
  Solution at $t_{\text{final}}=0.025$  with $\nu=10^{2}$
  for the SL-Upwind, SL-MUSCL, FKS and R-FKS schemes.
}
\label{fig:oscillcomp_temperature}
\end{center}
\end{figure}

% Data
For the first set of simulations, the final time is fixed at $t_{\text{final}}=0.025$, while we choose $\delta=0.02$. 
Then, an increasing sequence of spacial meshes is used for the four schemes: $M=600$, $1200$, $2400$ and $4800$,
while both the number of velocity cells $N=50$ and the velocity domain $L_v=[-30:30]$ are kept fixed.
In figure~\ref{fig:oscill600} we present the density, velocity and temperature
for the four schemes when $600$ cells are used. Observing this figure the schemes seem
to capture the same main waves and they seem to overlap apart from the R-FKS
which seems to oscillate. But a closer view shows that this is not the case. Thus, in order to enhance the differences between the schemes
and investigate the small scale structures of the solution, in the next figures, we present magnifications
on particular locations for each of the macroscopic variables.
In figure~\ref{fig:oscillcomp_density} we present the mesh convergence study
for the density variable with four different meshes: $600, 1200$, $2400$ and $4800$ cells.
Velocity and temperature variable magnifications results are displayed in figure~\ref{fig:oscillcomp_velocity}
and \ref{fig:oscillcomp_temperature} respectively. 
From these figures we observe that although the SL-Upwind scheme can capture the main flow, most
of the small scale structures are damped out due to an excessive numerical diffusion.
The same but less pronounced effect can be observed for the SL-MUSCL scheme. 
Indeed its embedded limiter\footnote{Note that any other limiter will have
more or less the same tendency because they are all designed to fulfill
a Discrete Maximum Principle, and, as such it clips all $t^{n+1}$ extrema 
to min/max $t^n$ data no matter if it represents a physically justified new extrema or if this extrema results from 
a numerical spurious oscillation.} 
has a know tendency to clip extrema, which in this extreme case yields to excessive numerical diffusion.
When the mesh is refined the SL-Upwind scheme starts to develop a wavy behavior
which is a remembrance of the true physical oscillations. SL-MUSCL scheme, being second-order
accurate is capturing more oscillations, but, is nonetheless still destroying the stair case initial
profile. Remarkably the (R-)FKS scheme is able to maintain the small scale structures of the 
solution even with the coarser mesh. When the mesh is refined some more details 
of the flow are further captured. The spacial mesh refinement
shows that SL-MUSCL, SL-Upwind and (R-)FKS converge towards the same solution. 
However, thanks to its almost exact transport phase, the FKS scheme family
is drastically reducing the numerical diffusion in this kinetic regime.
From those curves we can clearly see that the FKS family is genuinely able to maintain the structure of the solution for longer times and with
much less mesh points whereas the intrinsic numerical diffusion of 
semi-Lagrangian schemes damps such structures.
%Note that it is nt surprising that the classical FKS is producing comparable results
%with the new R-FKS. 
%=== B E G I N    F I G U R E ================================================
%=== E N D   F I G U R E ================================================

%=== B E G I N    F I G U R E ================================================
%=== E N D   F I G U R E ================================================

%=== B E G I N    F I G U R E ================================================
%=== E N D   F I G U R E ================================================

%=== B E G I N    F I G U R E ================================================
%=== E N D   F I G U R E ================================================

\subsubsection{Cost and efficiency study}

Although the numerical results presented in this paper are obtained in the one
dimensional case, it is nonetheless interesting to monitor the cost
of each numerical scheme. Two situations are considered; the first one
reports the CPU time at fixed meshes independent of the obtained accuracy\footnote{Basically
it measures the time we have to wait when spacial/velocity meshes are chosen.}, whereas
the second situation measures the cost for a given accuracy (that is spacial mesh refinement
is pursued up to the point the given error is attained).
Because the exact solution does not exist for this problem we have to design a reference diagnostics
to decide when a given scheme has attained an ``acceptable'' accuracy.

\paragraph{Cost at fixed meshes.} 
We present in table~\ref{tab:oscill_cost}
the cost of each of the previous simulations (from the mesh convergence study).
For each simulation we report the number of time steps $N_{\text{cycle}}$ and the CPU time $\text{T}$ in second(s).
We also compute
the time per cycle by $T_{\text{cycle}} = \text{T}/N_{\text{cycle}}$ and
time per cycle per physical cell $N_c$ by $T_{\text{cell}} = \text{T}/N_{\text{cycle}}/N_c$.
Last we propose the ratio of CPU time with respect to the corresponding SL-Upwind scheme results.
From this table we can see that at \underline{fixed meshes} FKS is faster than SL-Upwind
and SL-MUSCL and R-FKS are about two times more expensive than FKS. When SL-MUSCL and R-FKS are compared then we can observe
about $20\%$ speed-up using R-FKS. Of course in a multidimensional setting the situation is completely different as we have shown in \cite{FKS, FKS_HO}, in this
case the gain in term of computational cost is much larger for the FKS method. We expect the R-FKS method to be more expensive than FKS in more than
one dimension, since we pass from a piecewise constant reconstruction to a linear reconstruction which leads to some additional operations
but with a comparable computational cost \cite{FKS_HO2}.
%== T A B L E =========================================================
\begin{table}
  \begin{tabular}{|p{0.35cm}|cc|c||c|c|c|c|c|}
    \hline
    \multirow{2}{*}{\begin{sideways}\textbf{Code.}\end{sideways}}
    & &  &\textbf{Cell \#}    & \textbf{Cycle} & \textbf{Time} & T\textbf{/cycle} & T\textbf{/cell} & \textbf{Time ratio}\\
    & $N_v$ & {\begin{sideways} Vel.\end{sideways}}&$N_c\times N_v$ &  $N_{\text{cycle}}$  &  $\text{T}$ (s)  & $T_{\text{cycle}}$ (s) & $T_{\text{cell}}$ (s) & vs SL-Upwind \\
    \hline
    \hline
    \hline
    \multirow{8}{*}{\begin{sideways}  \textbf{SL-Upwind} \end{sideways}}  
    & \multirow{8}{*}{ $50$} &
    \multirow{8}{*}{\begin{sideways} $[-15,15]$  \end{sideways}}
    &$600\times 50$  & \multirow{2}{*}{$220$}   & \multirow{2}{*}{$0.772$ s}   & \multirow{2}{*}{$0.0035$}  & \multirow{2}{*}{$5.86\times 10^{-6}$} & \multirow{2}{*}{---}\\
    &&&$=3 \times 10^4$ & &  & & &\\
    \cline{4-9}
    &&&$1200 \times 50$  & \multirow{2}{*}{$441$}   & \multirow{2}{*}{$3.317$ s}   & \multirow{2}{*}{$0.0075$}  & \multirow{2}{*}{$6.27\times 10^{-6}$} & \multirow{2}{*}{---}\\
    &&&$= 6\times 10^{4}$ & & & & &\\
    \cline{4-9}
    &&&$2400 \times 50$ & \multirow{2}{*}{$882$}  & \multirow{2}{*}{$14.233$ s}   & \multirow{2}{*}{$0.0161$} & \multirow{2}{*}{$6.72\times 10^{-6}$} &  \multirow{2}{*}{---} \\
    &&&$= 12\times 10^4$ & &  &  & &\\
    \cline{4-9}
    &&&$4800 \times 50$ & \multirow{2}{*}{$1764$} & \multirow{2}{*}{$61.008$ s}&  \multirow{2}{*}{$0.0346$} & \multirow{2}{*}{$7.21\times 10^{-6}$} & \multirow{2}{*}{---}\\
    &&&$= 24\times 10^{4}$ & &  & & &\\
    \hline
    \hline
    \hline
    \multirow{8}{*}{\begin{sideways}  \textbf{SL-MUSCL} \end{sideways}}  
    & \multirow{8}{*}{ $50$} &
    \multirow{8}{*}{\begin{sideways} $[-15,15]$  \end{sideways}}
    &$600\times 50$  & \multirow{2}{*}{$220$}   & \multirow{2}{*}{$1.633$ s}   & \multirow{2}{*}{$0.0074$}  & \multirow{2}{*}{$1.24\times 10^{-5}$} & \multirow{2}{*}{$2.12$}\\
    &&&$=3 \times 10^4$ & &  & & &\\
    \cline{4-9}
    &&&$1200 \times 50$  & \multirow{2}{*}{$441$}   & \multirow{2}{*}{$6.784$ s}   & \multirow{2}{*}{$0.0154$}  & \multirow{2}{*}{$1.28\times 10^{-5}$} & \multirow{2}{*}{$2.05$}\\
    &&&$= 6\times 10^{4}$ & & & & &\\
    \cline{4-9}
    &&&$2400 \times 50$ & \multirow{2}{*}{$882$}  & \multirow{2}{*}{$29.166$ s}   & \multirow{2}{*}{$0.0331$} & \multirow{2}{*}{$1.38\times 10^{-5}$} &  \multirow{2}{*}{$2.05$} \\
    &&&$= 12\times 10^4$ & &  &  & &\\
    \cline{4-9}
    &&&$4800 \times 50$ & \multirow{2}{*}{$1764$} & \multirow{2}{*}{$275.812$ s}&  \multirow{2}{*}{$0.1564$} & \multirow{2}{*}{$3.26\times 10^{-5}$} & \multirow{2}{*}{$4.52$}\\
    &&&$= 24\times 10^{4}$ & &  & & &\\
    \hline
    \hline
    \hline
    \multirow{8}{*}{\begin{sideways}  \textbf{FKS} \end{sideways}} 
    & \multirow{8}{*}{ $50$} &
    \multirow{8}{*}{\begin{sideways} $[-15,15]$  \end{sideways}}
    &$600\times 50$  & \multirow{2}{*}{$220$}   & \multirow{2}{*}{$0.518$ s}   & \multirow{2}{*}{$0.0024$}  & \multirow{2}{*}{$3.92\times 10^{-6}$} & \multirow{2}{*}{$0.67$}\\
    &&&$=3 \times 10^4$ & &  & & &\\
    \cline{4-9}
    &&&$1200 \times 50$  & \multirow{2}{*}{$441$}   & \multirow{2}{*}{$2.101$ s}   & \multirow{2}{*}{$0.0048$}  & \multirow{2}{*}{$3.97\times 10^{-6}$} & \multirow{2}{*}{$0.63$}\\
    &&&$= 6\times 10^{4}$ & & & & &\\
    \cline{4-9}
    &&&$2400 \times 50$ & \multirow{2}{*}{$882$}  & \multirow{2}{*}{$8.754$ s}   & \multirow{2}{*}{$0.0099$} & \multirow{2}{*}{$4.14\times 10^{-6}$} &  \multirow{2}{*}{$0.62$} \\
    &&&$= 12\times 10^4$ & &  &  & &\\
    \cline{4-9}
    &&&$4800 \times 50$ & \multirow{2}{*}{$1764$} & \multirow{2}{*}{$38.544$ s}&  \multirow{2}{*}{$0.0220$} & \multirow{2}{*}{$4.55\times 10^{-6}$} & \multirow{2}{*}{$0.63$}\\
    &&&$= 24\times 10^{4}$ & &  & & &\\
    \hline    
    \hline
    \hline
    \multirow{8}{*}{\begin{sideways}  \textbf{R-FKS} \end{sideways}} 
    & \multirow{8}{*}{ $50$} &
    \multirow{8}{*}{\begin{sideways} $[-15,15]$  \end{sideways}}
    &$600\times 50$  & \multirow{2}{*}{$220$}   & \multirow{2}{*}{$1.503$ s}   & \multirow{2}{*}{$0.0068$}  & \multirow{2}{*}{$1.14\times 10^{-5}$} & \multirow{2}{*}{$1.95$}\\
    &&&$=3 \times 10^4$ & &  & & &\\
    \cline{4-9}
    &&&$1200 \times 50$  & \multirow{2}{*}{$441$}   & \multirow{2}{*}{$5.981$ s}   & \multirow{2}{*}{$0.0136$}  & \multirow{2}{*}{$1.13\times 10^{-5}$} & \multirow{2}{*}{$1.80$}\\
    &&&$= 6\times 10^{4}$ & & & & &\\
    \cline{4-9}
    &&&$2400 \times 50$ & \multirow{2}{*}{$882$}  & \multirow{2}{*}{$23.781$ s}   & \multirow{2}{*}{$0.0270$} & \multirow{2}{*}{$1.12\times 10^{-5}$} &  \multirow{2}{*}{$1.67$} \\
    &&&$= 12\times 10^4$ & &  &  & &\\
    \cline{4-9}
    &&&$4800 \times 50$ & \multirow{2}{*}{$1764$} & \multirow{2}{*}{$95.430$ s}&  \multirow{2}{*}{$0.0541$} & \multirow{2}{*}{$1.13\times 10^{-5}$} & \multirow{2}{*}{$1.56$}\\
    &&&$= 24\times 10^{4}$ & &  & & &\\
    \hline
  \end{tabular}
  \caption{ \label{tab:oscill_cost}
    Oscillating test case results for $600$, $1200$, $2400$ and $4800$ cells for the SL-Upwind, SL-MUSCL, FKS and R-FKS schemes.
    Time per cycle is obtained by $T_{\text{cycle}} = \text{T}/N_{\text{cycle}}$ and
    time per cycle per cell by $T_{\text{cell}} = \text{T}/N_{\text{cycle}}/N_c$.
  }
\end{table}
%=== E N D   T A B L E ================================================

\paragraph{Cost at fixed accuracy.}
However, as already mentioned, the previous table is not sufficient to give an idea of the performances of this class of methods. 
We also need to have a measure of the computational cost at fixed accuracy. For the oscillating test case we do not dispose of 
an exact reference solution. However, we have observed that both SL-MUSCL and the FKS family seem to capture an equivalent
solution with discontinuities at the same locations when the mesh is fine enough. Moreover thanks to the mesh convergence study we estimated that the characteristics length between
two discontinuuities is of the order $\Delta l = 4.8 \times 10^{-3}$ (this corresponds to $23$ cells of size $\Delta x = 1/4800$). 
%because the highly oscillating behavior of the solution makes the choice of the numerical scheme
%crucial. 
In addition from the SL-MUSCL data we estimated that SL-MUSCL scheme spreads those discontinuities 
over $\mathcal{N}=15$ cells. Therefore for the SL-MUSCL scheme it would require about 
$N_c \geq 2 \mathcal{N} / \delta l$ cells to capture each of the ``plateaus'' between two discontinuities.
In short SL-MUSCL scheme requires at minima $N_c=6250$ cells to provide an accurate
'enough' numerical solution. 
For this amount of cells, according to  table~\ref{tab:oscill_cost}, the number
of timesteps becomes $N_{\text{cycle}}=(0.3675) \, N_c$, indeed the
the ratio between timesteps and number of cell is constant equal to 
 $0.3675 = 1764/4800 = 882/2400 = 441/1200 (\simeq= 220/600)$.
Finaly, given the time needed to update one cell during one timestep, $T_{\text{cell}}$ from table~\ref{tab:oscill_cost}, we can estimate the CPU time needed by the 
SL-MUSCL scheme via equation: $\text{T}=N_c \times N_{\text{cycle}} \times T_{\text{cell}}$, 
that is about $466\text{s}\simeq 7\text{mn} 45\text{s}$.
In table~\ref{tab:oscill_cost2} we report these CPU times along with the other parameters from
the four schemes. As expected from the numerical results, the numerical diffusion of SL-Upwind
and SL-MUSCL schemes drastically penalizes their efficiency as they would need about $3$ hours and $8$ minutes respectively
to compute a valid solution when the FKS family rounds about $2$ to $3$ seconds.
 % RAPH ADD ON
We conclude by observing that, as usual, the computational cost analysis is a difficult task to accomplish and thus the numbers reported in the tables are
are only meant to provide a rough idea of the amount of resources that can already be saved in 1D by this family of schemes.
%\footnote{Even if $\mathcal{N}$ for the SL schemes would be divided by a factor $3$ or $4$ the gain would still be in favor of FKS family.}. 
Vastly more interesting gains are expected when 2D$\times$2D and more important 3D$\times$3D simulations will be performed.
 
%== T A B L E =========================================================
\begin{table}
  \begin{tabular}{|c||c|c|c|c|}
    \hline
                    & \textbf{SL-Upwind} &  \textbf{SL-MUSCL} & \textbf{FKS} & \textbf{R-FKS}  \\
    \hline
    $\mathcal{N}$   & $150$ & $15$ & $2$ & $2$ \\
    $\delta l$      & $4.8 \times 10^{-3}$ & $4.8 \times 10^{-3}$ & $4.8 \times 10^{-3}$ & $4.8 \times 10^{-3}$ \\
    \hline   
    $N_c=2\mathcal{N}/\delta l$         
                    & $62500$ & $6250$ & $833$ & $833$ \\
    $N_{\text{cycle}}=(0.3675) \, N_c$ & $22969$ & $2297$ & $306$ & $306$ \\     
    $T_{\text{cell}}$ (s) & $7\times 10^{-6}$ & $3.25\times 10^{-5}$ & $4.55\times 10^{-6}$ & $1.13\times 10^{-5}$\\
    \hline
    \hline
    \multirow{2}{*}{$\text{T}=N_c \times N_{\text{cycle}} \times T_{\text{cell}}$}
                   & $10048.83$s     & $466.55$s     & \multirow{2}{*}{$1.16$s} & \multirow{2}{*}{$2.88$s} \\ 
                   & ($2$h$47$mn) & ($7$mn$45$s) & & \\
    \hline
    \hline
  \end{tabular}
  \caption{ \label{tab:oscill_cost2}
    Oscillating test case results --- Estimation of the CPU time $\text{T}$
    needed to capture the plateaus of estimated characteristics length $\delta l$.
    The schemes diffuse discontinuities over $\mathcal{N}$ cells.
    The required number of cells is $N_c=2\mathcal{N}/\delta l$, the number
    of time steps $N_{\text{cycle}}=(0.3675) \, N_c$, where $0.3675=1764/4800$ 
    the ratio time steps/$N_c$ observed in table~\ref{tab:oscill_cost},
    %as well as the
    %time needed to update one cell during one timestep, $T_{\text{cell}}$.
    %Finally 
    $\text{T}=N_c \times N_{\text{cycle}} \times T_{\text{cell}}$.
  }
\end{table}
%=== E N D   T A B L E ================================================

%------------------------------------

%------------------------------------
% CONCLUSION
%

%
% CONCLUSIONS AND PERSPECTIVES
% 
\section{Conclusion and perspectives} \label{sec:conclusion}

% Conclusion
In this paper we have presented a new class of semi-Lagrangian methods based on high order polynomial reconstructions
of the distribution function which extends the Fast Kinetic Scheme developed in \cite{FKS, FKS_HO}. The main idea
is based on updating the extreme points of the distribution function instead of updating the solution on the cells centers in order
not to lose maxima and minima of the solution because of the numerical diffusion.
In the first part of this paper, the method has been presented in all its generality and different approaches to compute
the collision operator on the extrema are proposed. A critical analysis leads to choose the one which possesses the best properties
in term of precision and efficiency. 
In a second part, details are then given for the specific case of a piece-wise linear polynomial reconstruction and a BGK collision operator.
The technique is based on the 
%observation
assumption that the distribution function and its associated Maxwellian share 
the same space topology for equal velocity of the lattice. Finally, 
a min/max choices is made in order
to treat the discontinuities of the Maxwellian function at the interface points.

% Numerics
In the third part of the paper, this new FKS scheme with reconstruction has been tested on two problems.
The first one, a Riemann like problem, shows that in the fluid regime the 
method shares similar behaviors with respect a second-order accurate limited semi-Lagrangian (SL) scheme
whereas the original FKS reproduces the results of a first-order accurate semi-Lagrangian scheme. 
For intermediate regimes, the new scheme is still better than the original FKS and the second order accurate
semi-Lagrangian scheme but differences are extremely small. 
In a full kinetic regime this test case 
%is not able 
does not help 
to make a clear difference between first/second order accurate SL or FKS schemes.
Conversely the second test case is a highly oscillatory kinetic test
 presenting many discontinuities which are damped by the classical semi-Lagrangian schemes. 
Contrarily the FKS family is able to maintain the
discontinuities of the solution over time due to its natural anti-dissipation formulation.
The results obtained by the FKS schemes are therefore superior in terms of precision and in terms of computational costs.

% 
% % Pespectives
In the near future we will investigate the extension of this approach to the 
multi-dimensional case; the same reconstruction technique may apply but an analysis related to the choices to be done at interface surfaces/lines (3D/2D)
should be carried out. 
We plan also to extend the new FKS scheme to the case of Boltzmann collisional operator.

%------------------------------------

%==========================================================================================
%
%    A C K N O W LE D G E M E N T S 
%
\section*{Acknowledgements}
%This work has been supported by the French 'Agence Nationale pour la Recherche (ANR)' in the frame of the contract JCJC ``ALE INC(ubator) 3D''.
This work was partly supported by Campus France under the PHC Gallileo 2015-2016 number 32272UL.
%
%    F I N    de   A C K N O W LE D G E M E N T S
%
%
%===========================================================================================

%=============================================================================================
%
%    B I B L I O
%
%%% Local Variables:
%%% TeX-master: "MOOD_FKS.tex"
%%% End:

%
%   F I N   de   B I B L I O
%
%=============================================================================================


\begin{thebibliography}{99}
\small


% \bibitem{AlaiaII} {\sc A. Alaia, G. Puppo},
% {\it A hybrid method for hydrodynamic-kinetic flow - Part II - Coupling of hydrodynamic and kinetic models,}
% J. Comput. Phys., 231 (2012) pp.~5217--5242

%\bibitem{AA}
%{\sc O.~Aktas, and N.R.~Aluru}, {\it A combined continuum/DSMC
%technique for multiscale analysis of microfluidic filters}. J. Comp.
%Phys., vol.178, (2002) pp.~342--72.

%\bibitem{Babovsky}
%{\sc H.~Babovsky}, {\it On a simulation scheme for the Boltzmann
%equation}, Math. Methods Appl. Sci., 8 (1986), pp.~223--233.


\bibitem{bird}
{\sc{G.A.Bird}}, {\em{Molecular gas dynamics and direct simulation
of gas flows}}, Clarendon Press, Oxford (1994).

% \bibitem{Mieussens_localV_14}
% {\sc{C. Baranger, J. Claudel, N. H{\'e}rouard, L. Mieussens},} 
% {\it Locally refined discrete velocity grids for stationary rarefied flow simulations}, 
% J. Comput. Phys., 257(15), 572-593 (2014)

\bibitem{birsdall}
{\sc{C.K. Birsdall, A.B. Langdon}}, {\em{Plasma Physics Via Computer
Simulation}}, Institute of Physics (IOP), Series in Plasma Physics
(2004).

\bibitem{bobylev}
{\sc{A.V. Bobylev, A. Palczewski, J. Schneider,}} {\it{ On
approximation of the Boltzmann equation by discrete velocity
models.}} C. R. Acad. Sci. Paris Ser. I. Math. 320, (1995), pp.
639-644.

% \bibitem{Mieussens_localV2_14}
% {\sc{S. Brull, L. Mieussens}}, 
% {\it  Local discrete velocity grids for deterministic rarefied flow simulations }, 
% to appear J. Comput. Phys., (2014).

%\bibitem{BLPQ}
%{\sc J.~F.~Bourgat, P.~LeTallec, B.~Perthame, and Y.~Qiu}, {\it
%Coupling Boltzmann and Euler equations without overlapping}, in
%Domain Decomposition Methods in Science and Engineering, Contemp.
%Math. 157, AMS, Providence, RI, (1994), pp.~377--398.


%\bibitem{BLT}
%{\sc J.~F.~Bourgat, P.~LeTallec, M.D.~Tidriri}, {\it Coupling
%Boltzmann and Navier-Stokes equations by friction}. J. Comput. Phys.
%127, vol. 2 (1996), pag.~227--245.

% \bibitem{Boyd}
% {\sc J. Burt, I. Boyd}, {\it A low diffusion particle method for
% simulating compressible inviscid flows}, J. Comput. Phys., Vol.
% 227, (2008), pp. 4653-4670.

%\bibitem{Boyd2}
%{\sc J. Burt, I. Boyd}, {\it A hybrid particle approach for
%continuum and rarefied flow simulation}, J. Comput. Phys., Vol. 228,
%(2009), pp. 460-475.


% \bibitem{CPima}
% {\sc R.~E.~Caflisch, L.~Pareschi}, {\it Towards an hybrid method for
% rarefied gas dynamics}, IMA Vol. App. Math., vol. 135 (2004),
% pp.~57--73.

\bibitem{Cf}
{\sc R.~E.~Caflisch}, {\it Monte Carlo and Quasi-Monte Carlo
Methods}, Acta Numerica (1998)  pp.~1--49.

\bibitem{cercignani}
{\sc C.~Cercignani}, {\em The Boltzmann Equation and Its
Applications}, Springer-Verlag, New York, (1988).

\bibitem{Chacon}
{\sc L. Chacón, D. del-Castillo-Negrete, C.D. Hauck} ,
{\em  An asymptotic-preserving semi-Lagrangian algorithm for the time-dependent anisotropic heat transport equation},
J. Comp. Phys. (2014), vol. 272, pp. 719-746.

% \bibitem{Chen}
% {\sc Chen,S. and Doolen, G.D.} ,
% {\em  Lattice Boltzmann Method For Fluid Flows},
% Annu. Rev. Fluid Mech. 1998. %30:32964

\bibitem{Cheng}
{\sc C.Z. Cheng, G. Knorr}, {\it The integration of the Vlasov
equation in configuration space}, J. Comput. Phys., vol. 22
(1976), pp. 330-351.

% \bibitem{Clain}
% {\sc S. Clain, S. Diot, R. Loub{\`e}re} ,
% {\em A high-order finite volume method for systems of conservation laws—Multi-dimensional Optimal Order Detection (MOOD)}, 
% J. Comput. Phys. 230 (2011), pp. 4028–4050.
% 
% \bibitem{Clain1}
% {\sc S. Clain, S. Diot, R. Loub{\`e}re} ,
% {\em The multidimensional optimal order detection method in the three-dimensional case: very high-order finite volume method for hyperbolic systems},  
% Internat. J. Numer. Methods Fluids 73 (2013), pp. 362–392.

\bibitem{CrSon}
{\sc N. Crouseilles, T. Respaud, E. Sonnendrucker}, {\it A Forward
semi-Lagrangian Method for the Numerical Solution of the Vlasov
Equation}, Comp. Phys. Comm. 180, 10 (2009) 1730--1745.

\bibitem{CrSon1}
{\sc N. Crouseilles, M. Mehrenberger, E. Sonnendrucker}, {\it
Conservative semi-Lagrangian schemes for Vlasov equations}, J. Comp.
Phys., (2010) pp.~1927--1953.

%\bibitem{degond}
%{\sc N. Crouseilles, P.Degond, M.Lemou}, {\it A hybrid kinetic-fluid
%model for solving the gas-dynamics Boltzmann BGK equation}, J.
%Comput. Phys., vol. 199 (2004), pp. 776-808.

% \bibitem{DeVuyst1}
% {\sc F. De Vuyst, F. Salvarani}, {\it
% GPU-accelerated numerical simulations of the Knudsen gas on time-dependent domains,}
% Comput. Phys. Comm., 184, 3 (2013) pp. 532-536.
% 
% \bibitem{DeVuyst2}
% {\sc F. De Vuyst, C. Labourdette, C. Rey}, {\it
% GPU-accelerated real time visualization and interaction for coupled fluid dynamics},
% Submitted, preprint (\verb|http://hal.archives-ouvertes.fr/hal-00837555|)
% 


\bibitem{dimarco4}
{\sc P. Degond, G.~Dimarco, L. Pareschi}, {\it The Moment Guided
Monte Carlo Method}, Int. J. Num. Meth. Fluids, Vol.67, (2011), pp.
189-213.


\bibitem{Des}
{\sc L. Desvillettes, S. Mischler}, {\it About the splitting
algorithm for Boltzmann and BGK equations.} Math. Mod. \& Meth. in
App. Sci. bf 6, (1996), pp. 1079--1101.

\bibitem{FKS_HO2}
{\sc G.~Dimarco, C. Hauck, R. Loub{\`e}re}, {\it Multidimensional high order FKS schemes for the Boltzmann equation}, in progress.

\bibitem{dimarco-review}
{\sc G.~Dimarco, L. Pareschi}, {\it Numerical methods for kinetic equations}, ACTA Numerica, Vol. 23 (2014), pp. 369-520.

\bibitem{dimarco-IMEX}
{\sc G.~Dimarco, L. Pareschi}, {\it Asymptotic preserving implicit-explicit Runge-Kutta methods for non linear kinetic equations}, 
SIAM J. Num. Anal., Vol. 49, (2011), pp. 2057-2077 .. 

\bibitem{dimarco-Exp}
{\sc G.~Dimarco, L. Pareschi}, {\it Exponential Runge-Kutta methods for stiff kinetic equations}, 
SIAM J. Num. Anal., Vol. 51, (2013),  pp. 1064-1087. 

% \bibitem{dimarcoMC}
% {\sc P. Degond, G.~Dimarco}, {\it Fluid simulations with localized Boltzmann upscaling by direct Monte Carlo.},
% J. Comput. Phys., vol. 231 (2012), pp. 2414-2437.

%\bibitem{degond2}
%{\sc P. Degond, S. Jin, L. Mieussens}, {\it A smooth transition
%between kinetic and hydrodynamic equations.}, J. Comput. Phys., vol.
%209 (2005), pp. 665-694.

%\bibitem{degond3}
%{\sc P.Degond, J.-G. Liu, L. Mieussens}, {\it Macroscopic fluid
%models with localized kinetic upscaling effects}, SIAM MMS, vol. 5
%(2006), pp. 940-979.

%\bibitem{dimarco3}
%{\sc P. Degond, G. Dimarco, L. Mieussens.}, {\it A moving interface
%method for dynamic kinetic-fluid coupling.} J. Comput. Phys., Vol.
%227, pp. 1176-1208, (2007).

%\bibitem{dimarco4}
%{\sc P. Degond, G. Dimarco, L. Mieussens.}, {\it A multiscale
%kinetic-fluid solver with dynamic localization of kinetic effects.}
%J. Comput. Phys., Vol. 229, pp. 4907-4933, (2010).

\bibitem{dimarco2}
{\sc G.~Dimarco, L.~Pareschi}, {\it A Fluid Solver Independent
Hybrid method for Multiscale Kinetic Equations}, SIAM J. Sci.
Comput. Vol. 32, (2010), pp. 603--634.

\bibitem{dimarco3}
{\sc G.~Dimarco},
{\it The hybrid moment guided Monte Carlo method for the Boltzmann equation.}
Kin. Rel. Models, Vol 6, pp. 291-315 (2013).

%\bibitem{dimarco3}
%{\sc G.~Dimarco, L.~Pareschi}, {\it Hybrid multiscale methods II.
%Kinetic equations}, SIAM Mult. Model. and Simul. Vol 6., (2007), pp.
%1169--1197.

\bibitem{FKS}
{\sc{G. Dimarco, R. Loub{\`e}re}},
{\it Towards an ultra efficient kinetic scheme. Part I: basics on the BGK equation},
J. Comput. Phys., Vol. 255, 2013, pp 680-698.

\bibitem{FKS_HO}
{\sc{G. Dimarco, R. Loub{\`e}re}},
{\it Towards an ultra efficient kinetic scheme. Part II: the high order case},
J. Comput. Phys., Vol. 255, 2013, pp 699-719.

% \bibitem{FKS_DD}
% {\sc{G. Dimarco, R. Loub{\`e}re, V. Rispoli}},
% {\it A multiscale fast semi-Lagrangian method for rarefied gas dynamics},
% submitted  2014.


%\bibitem{DiMiRi13} {\sc G. Dimarco, L. Mieussens, V. Rispoli},
%{\it An asymptotic preserving automatic domain decomposition method for the
%Vlasov-Poisson-BGK system with applications to plasmas,} submitted to J. Comput. Phys. (2013)

%\bibitem{dimarco5}
%{\sc G.~Dimarco and L.~Pareschi}, {\it Domain decomposition
%techniques and hybrid multiscale methods for kinetic equations},
%Proceedings of the 11th International Conference on Hyperbolic
%problems: Theory, Numerics, Applications, pp. 457-464.

\bibitem{Filbet2} {\sc F.~Filbet, G.~Russo}, {\it High order numerical methods for the space non-homogeneous Boltzmann equation}.
J. Comput. Phys., 186  (2003), 457--480.

\bibitem{Filbet}
{\sc F. Filbet, E. Sonnendr\"ucker, P. Bertrand}, {\it Conservative
Numerical schemes for the Vlasov equation.} J. Comput. Phys. 172,
(2001) pp.~166--187.

%\bibitem{Frish}
%{\sc U. Frisch,  B. Hasslacher, and Y. Pomeau},
%{\it Lattice gas automata for the Navier-Stokes equations} Phys. Rev. Lett. 56, 1505 (1986).

% \bibitem{Frezzotti}
% {\sc A. Frezzotti, G. Ghiroldi, and L. Gibelli}, {\it Direct solution of the
% Boltzmann equation for a binary mixture on GPUs}, Proceedings of the
% 27th International Symposium on Rarefied Gas Dynamics, 2011, pp. 884–
% 889
% 
% \bibitem{Frezzotti2}
% {\sc A. Frezzotti, G. Ghiroldi, and L. Gibelli},{\it
% Solving model kinetic equations on GPUs}, Computers and Fluids, 50
% (2011), pp. 136–146.
% 
% \bibitem{Frezzotti3}
% {\sc A. Frezzotti, G. Ghiroldi, and L. Gibelli}, {\it Solving the Boltzmann equation on GPUs}, Comput. Phys. Comm.,
% 182 (2011), pp. 2445–2453.

\bibitem{Gamba}
{\sc  I.M. Gamba, S. H. Tharkabhushaman}, {\it Spectral - Lagrangian
based methods applied to computation of Non - Equilibrium
Statistical States.} J. Comput. Phys. 228, (2009) pp.~2012--2036.

% \bibitem{Gamba2}
% {\sc I. Gamba and J. Haack}, {\it
% High performance computing with a conservative spectral Boltzmann
% solver}, Proceedings of the 28th International Symposium on Rarefied
% Gas Dynamcis, 2012, pp. 334–341.

\bibitem{Gross}
{\sc E.P. Gross P.L. Bathnagar, M. Krook}, {\it A model for
collision processes in gases. {I}. small amplitude processes in
charged and neutral one-component systems} Phys. Rev. 94 (1954), pp.
511--525.

\bibitem{Gu}
{\sc Y. G\"{u}\c{c}l\"{u}, W.N.G. Hitchon}, {\it
A high order cell-centered semi-Lagrangian scheme for
multi-dimensional kinetic simulations of neutral gas flows.} J. Comput. Phys. 231, (2012) pp.~3289-3316.

\bibitem{Gu2}
{\sc Y. G\"{u}\c{c}l\"{u}, A.J. Christlieb, W.N.G. Hitchon}, {\it
Arbitrarily high order Convected Scheme solution of the Vlasov-Poisson system} J. Comput. Phys. 270, (2014) pp.~711-752.

% \bibitem{Haack}
% {\sc J. R. Haack}, {\it A hybrid OpenMP and MPI implementation of a
% conservative spectral method for the Boltzmann
% equation}. Submitted.	arXiv:1301.4195.

%\bibitem{HashHassan}
%{\sc D.B.~Hash and H.A.~Hassan},
%\newblock {\it Assessment of schemes for coupling Monte Carlo and Navier-Stokes solution
%methods}, J. Thermophys. Heat Transf., 10, (1996), pp.~242--249.

\bibitem{Hadji}
{\sc T. Homolle, N. Hadjiconstantinou}, {\it A low-variance
deviational simulation Monte Carlo for the Boltzmann equation.} J.
Comput. Phys., Vol 226 (2007), pp 2341--2358.

\bibitem{Hadji1}
{\sc T. Homolle, N. Hadjiconstantinou}, {\it Low-variance
deviational simulation Monte Carlo.} Phys. Fluids, Vol 19 (2007),
041701.

\bibitem{Jin}
{\sc  S. Jin}, {\it Efficient asymptotic-preserving (ap) schemes for some multiscale
kinetic equations}, SIAM J. Sci. Comput. Vol. 21  (1999), 441–454.

%\bibitem{Klar}
%{\sc A. Klar and C. Schmeiser}, {\it Numerical passage from
%radiative heat transfer to nonlinear diffusion models}, Math. Models
%Methods Appl. Sci., Vol. 11, pp.749-767, 2001.

%\bibitem{Kolobov12}
%{\sc V.I. Kolobov, R.R. Arslanbekov},
%{\it Towards adaptive kinetic-fluid simulations of weakly ionized plasmas},
%J. Comput. Phys. 231 (2012) 839--869

%\bibitem{Kobolov}
%{\sc V.I. Kolobov, R.R. Arslanbekov,V.V. Aristov, A.A. Frolova, S.A.
%Zabelok}, {\it Unified Solver for Rarefied and Continuum Flows with
%Adaptive Mesh and Algorithm Refinement}, J. Comput. Phys., Vol. 223
%No. 2, pp. 589--608 (2007).


%\bibitem{MUSCL}
%{\sc B. van Leer},
%{\it Towards the Ultimate Conservative Difference Scheme, V. A Second Order Sequel to Godunov's Method,}
%J. Comput. Phys. 32 (1979) 101--136

%\bibitem{Letallec}
%{\sc P.~LeTallec and F.~Mallinger},
%\newblock {\it Coupling Boltzmann and Navier-Stokes by half
%fluxes}
%\newblock {J. Comput. Phys.}, vol .136  (1997), pp.~51--67.

\bibitem{leveque} {\sc R.~J.~LeVeque},
   {\em Numerical Methods for Conservation Laws},
   Lectures in Mathematics, Birkhauser Verlag, Basel (1992).

%\bibitem{Levermore} {\sc D. Levermore, W.J.Morokoff, B.T. Nadiga}, {\it Moment realizability and the validity of the
%Navier Stokes equations for rarefied gas dynamics}, Physics of
%Fluids, Vol. 10, No. 12, (1998).

%\bibitem{sliu}
%{\sc S.~Liu}, {\em Monte Carlo strategies in scientific computing},
%Springer, (2004).

%\bibitem{ALE_INC}
%{\sc R. Loub{\`e}re},
%{\it First steps into ALE INC(ubator).
%A 2D arbitrary-Lagrangian-Eulerian code on general polygonal mesh for compressible flows: Version
%1.0.0}, Los Alamos National Laboratory, Report, LA-UR-04-8840, 2004.

% \bibitem{Malkov}
% {\sc E. Malkov and M. Ivanov}, {\it Parallelization of algorithms for solving the
% boltzmann equation for GPU-based computations}, Proceedings of the 27th
% International Symposium on Rarefied Gas Dynamics, 2011, pp. 946–951.

\bibitem{Mieussens}
{\sc L. Mieussens}, {\it Discrete Velocity Model and Implicit Scheme
for the BGK Equation of Rarefied Gas Dynamic}, Math. Models Meth.
App. Sci., Vol. 10, (2000), 1121--1149.

%\bibitem{Mouhot}
%{\sc C. Mouhot and L. Pareschi}, {\it Fast algorithms for computing the Boltzmann
%collision operator}, Math. Comp. Vol. 75, (2006), 1833-1852.


% \bibitem{Nanbu80}
% {\sc K.~Nanbu}, {\it Direct simulation scheme derived from the
% Boltzmann equation}, J. Phys. Soc. Japan, vol. 49 (1980),
% pp.~2042--2049.

%\bibitem{Noh87}
%{\sc W.F. Noh,} {\it Errors for calculations of strong shocks using an artificial viscosity and
%an artificial heat flux.}, J. Comput. Phys. 72, (1987), pp 78-120

\bibitem{Pal}
{\sc A. Palczewski, J. Schneider, A.V. Bobylev,} {\it A consistency
result for a discrete-velocity model of the Boltzmann equation.}
SIAM J. Numer. Anal. 34, (1997) pp. 1865-1883.

\bibitem{Pal1}
{\sc A. Palczewski, J. Schneider,} {\it Existence, stability, and
convergence of solutions of discrete velocity models to the
Boltzmann equation.} J. Statist. Phys. 91, (1998) pp. 307-326.

%\bibitem{PR}
%{\sc L.~Pareschi, G.~Russo}, {\it Time Relaxed Monte Carlo methods
%for the Boltzmann equation}, SIAM J. Sci. Comput. 23 (2001),
%pp.~1253--1273.

\bibitem{Pareschi}
{\sc L.~Pareschi, G.~Russo}, {\it Numerical solution of the Boltzmann equation I:
Spectrally accurate approximation of the collision operator}, SIAM J. Numer.
Anal. Vol. 37, pp. 1217–1245 (2000).

\bibitem{Pareschi2}
{\sc L. Pareschi, G. Russo and G. Toscani}, {\it Fast spectral methods for the Fokker–
Planck–Landau collision operator}, J. Comput. Phys. Vol 165, pp. 216–236 (2000).

\bibitem{PaTo}
  L. Pareschi, G. Toscani, Interacting multi-agent systems. Kinetic equations and Monte Carlo methods, Oxford University Press, USA, (2013).


%\bibitem{pa-inria}
%{\sc L.~Pareschi}, {\it Hybrid multiscale methods for hyperbolic and
%kinetic problems}, Esaim Proceedings, Vol. 15, T. Goudon, E.
%Sonnendrucker \& D. Talay Editors (2005), pp.87-120.

% \bibitem{Pieraccini}
% {\sc S. Pieraccini, G. Puppo}, {\it Implicit-explicit schemes for
% BGK kinetic equations}, J. Sci. Comp. (2007), pp. 1-28.


%\bibitem{LB}
%{\sc K.N. Premnath and J. Abraham.}
%{\it Three-dimensional multi-relaxation time (MRT) lattice-Boltzmann models
%for multiphase flow},
%J. Comput. Phys. ,Vol. 224( 2), 10 June 2007, 539-559


% \bibitem{Pullin78}
% {\sc D.~I.~Pullin}, {\it Direct simulation methods for compressible
% inviscid ideal gas flow}, J. Comput. Phys., 34 (1980), pp.~231--244.

%\bibitem{RGV1}
%{\sc R.~Roveda, D.B.~Goldstein, and P.L.~Varghese}, {\it Hybrid
%Euler/Particle Approach for Continuum/ Rarefied Flows}, AIAA J.
%Spacecraft Rockets 35, (1998), pp.~258--265.

%\bibitem{RGV2}
%{\sc R.~Roveda, D.B.~Goldstein, and P.L.~Varghese}, {\it Hybrid
%Euler/Direct Simulation Monte Carlo Calculation of Unsteady Slit
%Flow}, AIAA J. Spacecraft Rockets, vol. 37 (2000), pp.~753-760.

% \bibitem{Schwartz}
% {\sc T.E. Schwartzentruber, L.C. Scalabrin, I.D. Boyd}, {\it A
% modular particle-continuum numerical method for hypersonic
% non-equilibrium gas flows} J. Comput. Phys., Vol. 225 (2007), pp.
% 1159-1174.

%\bibitem{SSE}
%{\sc T. Schulze, P. Smereka and W. E}, {\it Coupling kinetic
%Monte-Carlo and continuum models with application to epitaxial
%growth}, J. Comput. Phys., vol. 189 (2003), pp.~197--211.

%\bibitem{Sedov}
%{\sc L.I. Sedov},
%%Sedov, L. I., "Propagation of strong shock waves," Journal of Applied Mathematics and Mechanics, Vol. 10, pages 241 - 250 (1946).
%{\em Similarity and Dimensional Methods in Mechanics} (New York: Academic), 1959

\bibitem{Shoucri}
{\sc M. Shoucri, G. Knorr}, {\it Numerical integration of the
Vlasov equation}. J. Comput. Phys., vol. 14 (1974), pp. 8492.

%\bibitem{capsule}
%{\sc U. Shiva Prasad, G. Srinivas}
%{\it
%Flow Simulation over Re-Entry Bodies at Supersonic \& Hypersonic Speeds,}
%International Journal of Engineering Research and Development, Vol. 2, pp. 29-34 (2012).

\bibitem{Strang}
{\sc G.~Strang}, {\it On the construction and the comparison of
    difference schemes}. SIAM J. Numer. Anal., (1968), pp.~506--517.

%\bibitem{SunBoyCan04} {\sc Q. Sun, I.D. Boyd, G.V. Candler},
%{\it A hybrid continuum-particle approach for modeling subsonic, rarefied gas flows},
%J. Comput. Phys., 194 (2004) pp.~256--277.

%\bibitem{TK}
%{\sc S. Tiwari, and A. Klar}, {\it An adaptive domain
%decomposition procedure for Boltzmann and Euler equations}. J.
%Comp. Appl. Math., vol.90, (1998) pp.~223--37.

%\bibitem{tiwari_JCP}
%  {\sc S. Tiwari}, {\it Coupling of the Boltzmann and Euler equations with automatic domain decomposition},
%  J. Comput. Phys., vol. 144, 1998, 710--726.

%\bibitem{tiwari_JCP1}
%{\sc S. Tiwari, A. Klar, S. Hardt}, {\it A Particle-Particle Hybrid
%Method for Kinetic and Continuum Equations}, J. Comput. Phys., Vol.
%228, (2009) pp. 7109-7124.

% \bibitem{Titarev}
% {\sc V. A. Titarev}, {\it Efficient deterministic modelling of three-dimensional rarefied gas flows.}
%  Commun. Comput. Phys. Vol., pp. 162192, (2012).

% \bibitem{Titarev2}
% {\sc V. A. Titarev, M. Dumbser, S. Utyuzhnikov}, {\it
% Construction and comparison of parallel implicit kinetic solvers in three spatial dimensions,}
% J. Comp. Phys., Vol 256, pp. 17-33 (2014).



 %\bibitem{Wad}
%{\sc D.C. Wadsworth, D.A. Erwin} {\it Two dimensional hybrid
%continuum/particle simulation approach for rarefied hypersonic
%flows.} AIAA Paper 92-2975, 1992.

%\bibitem{WanBoy03} {\sc W.L. Wang, I.D. Boyd},
%{\it Predicting Continuum breakdown in Hypersonic Viscous Flows},
%Phys. Fluids, 15 (2003) pp.~91--100.

%\bibitem{WijHad04} {\sc H.S. Wijesinghe, N.G. Hadjiconstantinou},
%{\it Discussion of Hybrid Atomistic-Continuum Methods for Multiscale Hydrodynamics},
%Int. J. Multi. Comp. Eng., 2 (2004) pp.~189--202.

%\bibitem{Wij}
%{\sc H.S. Wijesinghe, R.D. Hornung, A. L. Garcia, N. G.
%Hadjiconstantinou}, {\it Three-dimensional hybrid
%continuum-atomistic simulations for multiscale hydrodynamics.} J.
%Fluids Eng., Vol. 126, pp. 768-777, 2004.

%\bibitem{Wolf00}
%{\sc D.A. Wolf-Gladrow,} {\em Lattice-Gas Cellular Automata and Lattice Boltzmann Models: An Introduction},
%Springer-Verlag, Berlin-Heidelberg, 2000



\end{thebibliography}
\end{document}